%% file: CompactSurfaces_050111.tex
\input amstex
\let\myfrac=\frac
\input eplain
\let\frac=\myfrac
\input epsf




\loadeufm \loadmsam \loadmsbm
\message{symbol names}\UseAMSsymbols\message{,}

\font\myfontdefault=cmr10

\font\mytdmchapfont=cmb10 at 14pt
\font\mytdmheadfont=cmb10 at 10pt
\font\mytdmsubheadfont=cmr10

\magnification 1200
\newif\ifinappendices
\newif\ifundefinedreferences
\newif\ifchangedreferences
\newif\ifloadreferences
\newif\ifmakebiblio
\newif\ifmaketdm

\undefinedreferencesfalse
\changedreferencesfalse


\loadreferencestrue
\makebibliofalse
\maketdmfalse

\def\headpenalty{-400}     
\def\proclaimpenalty{-200} 

%
%

\def\alphanum#1{\ifcase #1 _\or A\or B\or C\or D\or E\or F\or G\or H\or I\or J\or K\or L\or M\or N\or O\or P\or Q\or R\or S\or T\or U\or V\or W\or X\or Y\or Z\fi}
\def\gobbleeight#1#2#3#4#5#6#7#8{}

\newwrite\references
\newwrite\tdm
\newwrite\biblio

\newcount\chapno
\newcount\headno
\newcount\subheadno
\newcount\procno
\newcount\figno
\newcount\citationno

\def\setcatcodes{%
\catcode`\!=0 \catcode`\\=11}%

\ifloadreferences
    {\catcode`\@=11 \catcode`\_=11%
    \input references.tex %
    }%
\else
    \openout\references=references.tex
\fi

\newcount\newchapflag 
\newcount\showpagenumflag 

\global\chapno = -1 
\global\citationno=0
\global\headno = 0
\global\subheadno = 0
\global\procno = 0
\global\figno = 0

\def\resetcounters{%
\global\headno = 0%
\global\subheadno = 0%
\global\procno = 0%
\global\figno = 0%
}

\global\newchapflag=0 
\global\showpagenumflag=0 

\def\chinfo{\ifinappendices\alphanum\chapno\else\the\chapno\fi}%
\def\headinfo{\the\headno}
\def\subheadinfo{\the\headno.\the\subheadno}
\def\procinfo{\the\headno.\the\procno}
\def\figinfo{\the\headno.\the\figno}
\def\citationinfo{\the\citationno}%
\def\nextheadno{\global\advance\headno by 1 \global\subheadno = 0 \global\procno = 0}
\def\nextsubheadno{\global\advance\subheadno by 1}
\def\nextprocno{\global\advance\procno by 1 \procinfo}
\def\nextfigno{\global\advance\figno by 1 \figinfo}

{\global\let\noe=\noexpand%
%
%
\catcode`\@=11%
\catcode`\_=11%
\setcatcodes%
!global!def!_@@internal@@makeref#1{%
!global!expandafter!def!csname #1ref!endcsname##1{%
!csname _@#1@##1!endcsname%
!expandafter!ifx!csname _@#1@##1!endcsname!relax%
    !write16{#1 ##1 not defined - run saving references}%
    !undefinedreferencestrue%
!fi}}%
!global!def!_@@internal@@makelabel#1{%
!global!expandafter!def!csname #1label!endcsname##1{%
!edef!temptoken{!csname #1info!endcsname}%
!ifloadreferences%
    !expandafter!ifx!csname _@#1@##1!endcsname!relax%
        !write16{#1 ##1 not hitherto defined - rerun saving references}%
        !changedreferencestrue%
    !else%
        !expandafter!ifx!csname _@#1@##1!endcsname!temptoken%
        !else
            !write16{#1 ##1 reference has changed - rerun saving references}%
            !changedreferencestrue%
        !fi%
    !fi%
!else%
    !expandafter!edef!csname _@#1@##1!endcsname{!temptoken}%
    !edef!textoutput{!write!references{\global\def\_@#1@##1{!temptoken}}}%
    !textoutput%
!fi}}%
!global!def!makecounter#1{!_@@internal@@makelabel{#1}!_@@internal@@makeref{#1}}%
!unsetcatcodes%
}
\makecounter{ch}%
\makecounter{head}%
\makecounter{subhead}%
\makecounter{proc}%
\makecounter{fig}%
\makecounter{citation}%
\def\newref#1#2{%
\def\temptext{#2}%
\edef\bibliotextoutput{\expandafter\gobbleeight\meaning\temptext}%
\global\advance\citationno by 1\citationlabel{#1}%
\ifmakebiblio%
    \edef\fileoutput{\write\biblio{\noindent\hbox to 0pt{\hss$[\the\citationno]$}\hskip 0.2em\bibliotextoutput\medskip}}%
    \fileoutput%
\fi}%
\def\cite#1{%
$[\citationref{#1}]$%
\ifmakebiblio%
    \edef\fileoutput{\write\biblio{#1}}%
    \fileoutput%
\fi%
}%
%
%
%

\let\mypar=\par


\def\raggedleft{\leftskip=0pt plus 1fil \parfillskip=0pt}


\font\lettrinefont=cmr10 at 28pt
\def\lettrine #1[#2][#3]#4%
{\hangafter -#1 \hangindent #2
\noindent\hskip -#2 \vtop to 0pt{
\kern #3 \hbox to #2 {\lettrinefont #4\hss}\vss}}

\font\mylettrinefont=cmr10 at 28pt
\def\mylettrine #1[#2][#3][#4]#5%
{\hangafter -#1 \hangindent #2
\noindent\hskip -#2 \vtop to 0pt{
\kern #3 \hbox to #2 {\mylettrinefont #5\hss}\vss}}


\edef\Pagetitle={Blank}

\headline={\hfil\Pagetitle\hfil}

\footline={\hfil\myfontdefault\folio\hfil}

\def\nextoddpage
{
\newpage%
\ifodd\pageno%
\else%
    \global\showpagenumflag = 0%
    \null%
    \vfil%
    \eject%
    \global\showpagenumflag = 1%
\fi%
}


\def\newchap#1#2%
{%
%
%
\global\advance\chapno by 1%
\resetcounters%
%
%
\newpage%
\ifodd\pageno%
\else%
    \global\showpagenumflag = 0%
    \null%
    \vfil%
    \eject%
    \global\showpagenumflag = 1%
\fi%
\global\newchapflag = 1%
\global\showpagenumflag = 1%
%
%
{\font\chapfontA=cmsl10 at 30pt%
\font\chapfontB=cmsl10 at 25pt%
\null\vskip 5cm%
{\chapfontA\raggedleft\hfil%
{%
\ifnum\chapno=0
    \phantom{%
    \ifinappendices%
        Annexe \alphanum\chapno%
    \else%
        \the\chapno%
    \fi}%
\else%
    \ifinappendices%
        Annexe \alphanum\chapno%
    \else%
        \the\chapno%
    \fi%
\fi%
}%
\par}%
\vskip 2cm%
{\chapfontB\raggedleft%
\lineskiplimit=0pt%
\lineskip=0.8ex%
\hfil #1\par}%
\vskip 2cm%
}%
\edef\Pagetitle{#2}%
%
%
\ifmaketdm%
    \def\temp{#2}%
    \def\tempbis{\nobreak}%
    \edef\chaptitle{\expandafter\gobbleeight\meaning\temp}%
    \edef\mynobreak{\expandafter\gobbleeight\meaning\tempbis}%
    \edef\textoutput{\write\tdm{\bigskip{\noexpand\mytdmchapfont\noindent\chinfo\ - \chaptitle\hfill\noexpand\folio}\par\mynobreak}}%
\fi%
\textoutput%
}


\def\newhead#1%
{%
\ifhmode%
    \mypar%
\fi%
\ifnum\headno=0%
\else%
    \nobreak\vskip -\lastskip%
    \nobreak\vskip .5cm%
    \line{\hfil$\diamond$\hfil}%
    \penalty\headpenalty\vskip .5cm%
\fi%
\nextheadno%
\ifmaketdm%
    \def\temp{#1}%
    \edef\sectiontitle{\expandafter\gobbleeight\meaning\temp}%
    \edef\textoutput{\write\tdm{\noindent{\noexpand\mytdmheadfont\quad\headinfo\ - \sectiontitle\hfill\noexpand\folio}\par}}%
    \textoutput%
\fi%
\font\headfontA=cmbx10 at 14pt%
{\headfontA\noindent\headinfo\ -\ #1.\hfil}%
\nobreak\vskip .5cm%
}%


\def\newsubhead#1%
{%
\ifhmode%
    \mypar%
\fi%
\ifnum\subheadno=0%
\else%
    \penalty\headpenalty\vskip .4cm%
\fi%
\nextsubheadno%
\ifmaketdm%
    \def\temp{#1}%
    \edef\subsectiontitle{\expandafter\gobbleeight\meaning\temp}%
    \edef\textoutput{\write\tdm{\noindent{\noexpand\mytdmsubheadfont\quad\quad\subheadinfo\ - \subsectiontitle\hfill\noexpand\folio}\par}}%
    \textoutput%
\fi%
\font\subheadfontA=cmsl10 at 12pt
{\subheadfontA\noindent\subheadinfo\ #1.\hfil}%
\nobreak\vskip .25cm%
}%

%
%


\font\mathromanten=cmr10
\font\mathromanseven=cmr7
\font\mathromanfive=cmr5
\newfam\mathromanfam
\textfont\mathromanfam=\mathromanten
\scriptfont\mathromanfam=\mathromanseven
\scriptscriptfont\mathromanfam=\mathromanfive
\def\mathroman{\fam\mathromanfam}


\font\sansseriften=cmss10
\font\sansserifseven=cmss7
\font\sansseriffive=cmss5
\newfam\sansseriffam
\textfont\sansseriffam=\sansseriften
\scriptfont\sansseriffam=\sansserifseven
\scriptscriptfont\sansseriffam=\sansseriffive
\def\mathsf{\fam\sansseriffam}


\font\boldten=cmb10
\font\boldseven=cmb7
\font\boldfive=cmb5
\newfam\mathboldfam
\textfont\mathboldfam=\boldten
\scriptfont\mathboldfam=\boldseven
\scriptscriptfont\mathboldfam=\boldfive
\def\mathbf{\fam\mathboldfam}


\font\mycmmiten=cmmi10
\font\mycmmiseven=cmmi7
\font\mycmmifive=cmmi5
\newfam\mycmmifam
\textfont\mycmmifam=\mycmmiten
\scriptfont\mycmmifam=\mycmmiseven
\scriptscriptfont\mycmmifam=\mycmmifive

\def\hexa#1{\ifcase #1 0\or 1\or 2\or 3\or 4\or 5\or 6\or 7\or 8\or 9\or A\or B\or C\or D\or E\or F\fi}
\mathchardef\mathi="7\hexa\mycmmifam7B
\mathchardef\mathj="7\hexa\mycmmifam7C


\font\mymsbmten=msbm10 at 8pt
\font\mymsbmseven=msbm7 at 5.6pt
\font\mymsbmfive=msbm5 at 4pt
\newfam\mymsbmfam
\textfont\mymsbmfam=\mymsbmten
\scriptfont\mymsbmfam=\mymsbmseven
\scriptscriptfont\mymsbmfam=\mymsbmfive

\mathchardef\mybeth="7\hexa\mymsbmfam69
\mathchardef\mygimmel="7\hexa\mymsbmfam6A
\mathchardef\mydaleth="7\hexa\mymsbmfam6B


\def\placelabel[#1][#2]#3{{%
\setbox10=\hbox{\raise #2cm \hbox{\hskip #1cm #3}}%
\ht10=0pt%
\dp10=0pt%
\wd10=0pt%
\box10}}%

\def\placefigure#1#2#3{%
\medskip%
\midinsert%
\vbox{\line{\hfil#2\epsfbox{#3}#1\hfil}%
\vskip 0.3cm%
\line{\noindent\hfil\sl Figure \nextfigno\hfil}}%
\medskip%
\endinsert%
}


\newif\ifinproclaim%
\global\inproclaimfalse%
\def\proclaim#1{%
\medskip%
%
%
\bgroup%
\inproclaimtrue%
\setbox10=\vbox\bgroup\leftskip=0.8em\noindent{\bf #1}\sl%
}

\def\endproclaim{%
\egroup%
\setbox11=\vtop{\noindent\vrule height \ht10 depth \dp10 width 0.1em}%
\wd11=0pt%
\setbox12=\hbox{\copy11\kern 0.3em\copy11\kern 0.3em}%
\wd12=0pt%
\setbox13=\hbox{\noindent\box12\box10}%
\noindent\unhbox13%
\egroup%
\medskip\ignorespaces%
}

\def\proclaim#1{%
\medskip%
\bgroup%
\inproclaimtrue%
\noindent{\bf #1}%
\nobreak\medskip%
\sl%
}

\def\endproclaim{%
\mypar\egroup\penalty\proclaimpenalty\medskip\ignorespaces%
}

\def\noskipproclaim#1{%
\medskip%
\bgroup%
\inproclaimtrue%
\noindent{\bf #1}\nobreak\sl%
}

\def\endnoskipproclaim{%
\mypar\egroup\penalty\proclaimpenalty\medskip\ignorespaces%
}


\def\ninn{{n\in\Bbb{N}}}

\def\colon{\hphantom{0}:}
\def\proof{{\noindent\bf Proof:\ }}

\def\remark{{\noindent\sl Remark:\ }}
\def\suite#1#2{({#1}_{#2})_{#2\in\Bbb{N}}}
\def\mlim{\mathop{{\mathroman Lim}}}
\def\mlimsup{\mathop{{\mathroman LimSup}}}

\def\msup{\mathop{{\mathroman Sup}}}
\def\minf{\mathop{{\mathroman Inf}}}
\def\msf#1{{\mathsf #1}}

\def\qed{~$\square$}
\def\munion{\mathop{\cup}}
\def\minter{\mathop{\cap}}
\def\myitem#1{%
\ifinproclaim%
    \item{#1}%
\else%
    \noindent\hbox to .5cm{\hfill#1\hss}
\fi}%

\catcode`\@=11
\def\Eqalign#1{\null\,\vcenter{\openup\jot\m@th\ialign{%
\strut\hfil$\displaystyle{##}$&$\displaystyle{{}##}$\hfil%
&&\quad\strut\hfil$\displaystyle{##}$&$\displaystyle{{}##}$%
\hfil\crcr #1\crcr}}\,}
\catcode`\@=12

\def\makeop#1{%
\global\expandafter\def\csname op#1\endcsname{{\mathroman #1}}}%

\def\makeopsmall#1{%
\global\expandafter\def\csname op#1\endcsname{{\mathroman{\lowercase{#1}}}}}%

\makeopsmall{ArcTan}%
\makeopsmall{ArcCos}%
\makeop{Arg}%
\makeop{Det}%
\makeop{Log}%
\makeop{Re}%
\makeop{Im}%
\makeop{Dim}%
\makeopsmall{Tan}%
\makeop{Ker}%
\makeopsmall{Cos}%
\makeopsmall{Sin}%
\makeop{Exp}%
\makeopsmall{Tanh}%
\makeop{Tr}%
\makeop{End}%
\makeop{Long}%
\makeop{Ch}%
\makeop{Exp}%
\makeop{Int}%
\makeop{Ext}%
\makeop{Aire}%
\makeop{Im}%
\makeop{Conf}%
\makeop{Exp}%
\makeop{Mod}%
\makeop{Log}%
\makeop{Ext}%
\makeop{Int}%
\makeop{Dist}%
\makeop{Aut}%
\makeop{Id}%
\makeop{GL}%
\makeop{SO}%
\makeop{Homeo}%
\makeop{Vol}%
\makeop{Ric}%
\makeop{Hess}%
\makeop{Euc}%
\makeop{Isom}%
\makeop{Max}%
\makeop{Long}%
\makeop{Fixe}%
\makeop{Wind}%
\makeop{Mush}%
\makeop{Ad}%
\makeop{loc}%
\makeop{Len}%
\makeop{Area}%

\let\emph=\bf

\hyphenation{quasi-con-formal}

%
%

\ifmakebiblio%
    \openout\biblio=biblio.tex %
    {%
        \edef\fileoutput{\write\biblio{\bgroup\leftskip=2em}}%
        \fileoutput
    }%
\fi%

\newref{BallGromSch}{Ballman W., Gromov M., Schroeder V., {\sl Manifolds of nonpositive curvature}, Progress in Mathematics, {\bf 61}, Birkh\"auser, Boston, (1985)}
\newref{Grom}{Gromov M., Foliated plateau problem, part I : Minimal varieties, {\sl GAFA} {\bf 1}, no. 1, (1991), 14--79}%
\newref{LabB}{Labourie F., Probl\`emes de Monge-Amp\`ere, courbes holomorphes et laminations, {\sl GAFA} {\bf 7}, no. 3, (1997), 496--534}
\newref{LabA}{Labourie F., Un lemme de Morse pour les surfaces convexes, {\sl Invent. Math.} {\bf 141} (2000), 239--297}
\newref{LV}{Lehto O., Virtanen K. I., {\sl Quasiconformal mappings in the plane}, Die Grund\-lehren der mathemathischen Wissenschaften, {\bf 126}, Springer-Verlag, New York-Heidelberg, (1973)}
\newref{Muller}{Muller M.P., Gromov's Schwarz lemma as an estimate of the gradient for holomorphic curves, In: {\sl Holomorphic curves in symplectic geometry},
Progress in Mathematics, {\bf 117}, Birkh\"auser, Basel, (1994)}
\newref{RosSpruck}{Rosenberg H., Spruck J. On the existence of convex hyperspheres of constant Gauss curvature in hyperbolic space, {\sl J. Diff. Geom.} {\bf 40} (1994), no. 2,
379--409} 
\newref{SmiA}{Smith G., Special Legendrian structures and Weingarten problems, Preprint, Orsay (2005)}%
\newref{SmiB}{Smith G., Hyperbolic Plateau problems, Preprint, Orsay (2005)}%
\newref{SmiE}{Smith G., Th\`ese de doctorat, Paris (2004)}%

\ifmakebiblio%
    {\edef\fileoutput{\write\biblio{\egroup}}%
    \fileoutput}%
\fi%

%
%
%
\document
\myfontdefault
\global\chapno=1
\global\showpagenumflag=1
\def\Pagetitle{}
\null
\vfill
\def\centre{\rightskip=0pt plus 1fil \leftskip=0pt plus 1fil \spaceskip=.3333em \xspaceskip=.5em \parfillskip=0em \parindent=0em}%
\def\textmonth#1{\ifcase#1\or January\or Febuary\or March\or April\or May\or June\or July\or August\or September\or October\or November\or December\fi}
\font\abstracttitlefont=cmr10 at 14pt
{\abstracttitlefont\centre Pointed $k$-surfaces\par}
\bigskip
{\centre Graham Smith\par}
\bigskip
{\centre \the\day\ \textmonth\month\ \the\year\par}
\bigskip
{\centre Equipe de topologie et dynamique,\par
Laboratoire des math\'ematiques,\par
B\^atiment 425,\par
UFR des sciences d'Orsay,\par
91405 Orsay CEDEX, FRANCE\par}
\bigskip
\noindent{\emph Abstract:\ }Following on from \cite{SmiB}, we provide a complete geometric description of solutions to the Plateau
problem $(S,\varphi)$ when $S$ is a compact Riemann surface with a finite number of points removed.
\vfill
\nextoddpage
\def\Pagetitle{\sl Pointed $k$-surfaces}
\newhead{Introduction}
\noindent In this paper, by establishing a result permitting us to describe the behaviour ``at infinity'' of surfaces of constant
Gaussian curvature immersed in three dimensional hyperbolic space, we obtain a complete geometric description of solutions to
the Plateau problem for compact Riemann surfaces with marked points.
\medskip
\noindent Let $\Bbb{H}^3$ be three dimensional hyperbolic space, and let $\partial_\infty\Bbb{H}^3$ be its ideal boundary (see, for
example \cite{BallGromSch}). The ideal boundary of $\Bbb{H}^3$ may be identified canonically with the Riemann sphere $\hat{\Bbb{C}}$. In
this context, following \cite{LabA} and \cite{SmiB}, we define a {\emph Plateau problem} to be a pair $(S,\varphi)$ where $S$ is a Riemann
surface and $\varphi:S\rightarrow\partial_\infty\Bbb{H}^3$ is a locally conformal mapping (i.e. a locally homeomorphic holomorphic mapping). The Plateau problem $(S,\varphi)$ is said to be of hyperbolic, parabolic or elliptic type depending on whether $S$ is hyperbolic,
parabolic or elliptic respectively.
\medskip
\noindent Let $U\Bbb{H}^3$ be the unitary bundle over $\Bbb{H}^3$. For $i:S\rightarrow\Bbb{H}^3$ an immersion, using the canonical
orientation of $S$, we may define the unit normal exterior vector field $\msf{N}$ over $S$. This field is a section of $U\Bbb{H}^3$ over $i$.
We define the {\emph Gauss lifting} $\hat{\mathi}$ of $i$ by $\hat{\mathi}=\msf{N}$. We define a {\emph $k$-surface} to be an immersed
surface $\Sigma=(S,i)$ in $\Bbb{H}^3$ of constant Gaussian curvature $k$ whose Gauss lifting $\hat{\Sigma}=(S,\hat{\mathi})$ is a complete
immersed surface in $U\Bbb{H}^3$. For $k\in (0,1)$, a solution to the Plateau problem $(S,\varphi)$ is a $k$-surface $\Sigma=(S,i)$ such
that, if we denote by $\overrightarrow{n}$ the Gauss-Minkowski mapping of $\Bbb{H}^3$, then the Gauss lifting $\hat{\mathi}$ of $i$
satisfies:
$$
\varphi = \overrightarrow{n}\circ\hat{\mathi}.
$$
\noindent In \cite{SmiB} we show that, if $(S,\varphi)$ is a hyperbolic Plateau problem, then, for all $k\in (0,1)$ there exists a unique
solution $i$ to the Plateau problem $(S,\varphi)$ with constant Gaussian curvature $k$. Moreover, we show that $i$ depends continuously on
$\varphi$. In this paper, following on from these ideas, we study the structure of solutions to the Plateau problem $(S,\varphi)$ when
$S$ is a compact Riemann surface with isolated marked points. The following result, which provides the key to the rest of the paper,
describes the behaviour ``at infinity'' of solutions to the Plateau problem:
\proclaim{Theorem \nextprocno\ {\sl Boundary Behaviour Theorem}}
\noindent let $S$ be a hyperbolic Riemann surface and let $\varphi:S\rightarrow\hat{\Bbb{C}}$ be a locally conformal mapping. For $k\in (0,1)$, let
$i:S\rightarrow U\Bbb{H}^3$ be an immersion such that $(S,i)$ is the unique solution to the Plateau problem $(S,\varphi)$
with constant Gaussian curvature $k$.
\medskip
\noindent Let $K$ be a compact subset of $S$ and let $\Omega$ be a connected component of $S\setminus K$. Let
$q$ be an arbitrary point in the boundary of $\varphi(\Omega)$ that is not in $\varphi(\overline{\Omega}\minter K)$.
If $\suite{p}{n}\in\Omega$ is a sequence of points such that $(\varphi(p_n))_{n\in\Bbb{N}}$ tends towards $q$, then the sequence $(i(p_n))_{n\in\Bbb{N}}$ also tends towards $q$.
\endproclaim
\proclabel{PresentationChIIILimites}
\remark This theorem confirms our intuition concerning solutions to the Plateau problem. In particular, if $S$ is a Jordan domain in
$\partial_\infty\Bbb{H}^3$, if $\varphi$ is the canonical embedding and if $i:S\rightarrow\Bbb{H}^3$ is a solution to the Plateau
problem $(S,\varphi)$, then the ideal boundary of the immersed surface $(S,i)$ coincides with $\partial S$.
\medskip
\noindent We use this theorem to study the behaviour of solutions to the Plateau problem near to isolated singularities. We begin by
a series of definitions concerning tubes about geodesics. For $\Gamma$ a geodesic in $\Bbb{H}^3$, we define $N_\Gamma$ to be the
normal bundle over $\Gamma$ in $U\Bbb{H}^3$:
$$
N_\Gamma = \left\{ n_p\in U\Bbb{H}^3\text{ s.t. }p\in\Gamma, n_p\perp T_p\Gamma\right\}.
$$
\noindent A {\emph tube} about $\Gamma$ is a pair $T=(S,\hat{\mathi})$ where $S$ is a complete surface and
$\hat{\mathi}:S\rightarrow N_\Gamma$ is a covering map. Since $N_\Gamma$ is conformally equivalent to $S^1\times\Bbb{R}$, where $S^1$ is the
circle of radius $1$ in $\Bbb{C}$, we may assume either that $S=S^1\times\Bbb{R}$ or that $S=\Bbb{R}\times\Bbb{R}$. In the former case,
$\hat{\mathi}$ is a covering map of finite order, and, if $k$ is the order of $\hat{\mathi}$, then we say that the tube $T$ is a
{\emph tube of order $k$}. The application $\hat{\mathi}$ is then unique up to vertical translations and horizontal rotations of
$S^1\times\Bbb{R}$. In the latter case, we say that $T$ is a {\emph tube of infinite order}. The application $\hat{\mathi}$ is then
unique up to translations of $\Bbb{R}\times\Bbb{R}$. In the sequel, we will only be interested in tubes of finite order.
\medskip
\noindent Let $S$ be a compact surface and let $\Cal{P}$ be a finite set of points in $S$. Let
$\hat{\mathi}:S\setminus\Cal{P}\rightarrow U\Bbb{H}^3$ be an immersion. Let $p$ be an arbitrary point in $\Cal{P}$. We say that
$(S\setminus\Cal{P},\hat{\mathi})$ is {\emph asymtotically tubular} of order $k$ about $p$ if and only if it is a bounded graph over
a half tube of order $k$ in $U\Bbb{H}^3$, which tends towards the tube itself as one tends towards infinity. More precisely, let $\opExp:TU\Bbb{H}^3\rightarrow U\Bbb{H}^3$ be the exponential mapping and let $NN_\Gamma$ be the normal bundle of $N_\Gamma$.
$(S\setminus\Cal{P},\hat{\mathi})$ is asymptotically tubular of order $k$ about $p$ if and only if there exists:
\medskip
\myitem{(i)} a geodesic $\Gamma$ and a tube $T=(S^1\times\Bbb{R},\hat{\mathj})$ of order $k$ about $\Gamma$,
\medskip
\myitem{(ii)} a section $\lambda$ of $\hat{\mathj}^*NN_\Gamma$ over $S^1\times (0,\infty)$,
\medskip
\myitem{(iii)} a neighbourhood $\Omega$ of $p$ in $S$ such that $\Cal{P}\minter\Omega=\left\{p\right\}$, and
\medskip
\myitem{(iv)} a diffeomorphism $\alpha:S^1\times (0,\infty)\rightarrow\Omega\setminus\left\{p\right\}$,
\medskip
\noindent such that:
\medskip
\myitem{(i)} $\hat{\mathi}\circ\alpha = \opExp\circ\lambda$,
\medskip
\myitem{(ii)} $\alpha(e^{i\theta},t)\rightarrow p$ as $t\rightarrow\infty$, and
\medskip
\myitem{(iii)} for all $p\in\Bbb{N}$, the derivative $D^p\lambda(e^{i\theta},t)$ tends to zero as $t$ tends to $+\infty$.
\medskip
\noindent We now obtain the following result:
\proclaim{Theorem \nextprocno}
\noindent Let $S$ be a Riemann surface. Let $\Cal{P}$ be a discrete subset of $S$ such that $S\setminus\Cal{P}$ is hyperbolic. Let $\varphi:S\rightarrow\hat{\Bbb{C}}$ be a ramified covering having critical points in $\Cal{P}$. Let $\kappa$ be a real number in $(0,1)$. Let
$i:S\setminus\Cal{P}\rightarrow\Bbb{H}^3$ be the unique solution to the Plateau problem $(S\setminus\Cal{P},\varphi)$ with constant Gaussian
curvature $\kappa$. Let $\hat{\Sigma}=(S\setminus\Cal{P},\hat{\mathi})$ be the Gauss lifting of $\Sigma$.
\medskip
\noindent Let $p_0$ be an arbitrary point in $\Cal{P}$. If $\varphi$ has a critical point of order $k$ at $p_0$, then $\hat{\Sigma}$ is asymptotically tubular of order
$k$ at $p_0$.
\endproclaim
\proclabel{PresentationChIIIRevetementsRamifiees}
\remark This means that if the mapping $\varphi$ has a critical point of order $k$ at $p_0$, and is thus equivalent to $z\mapsto z^k$, then
the immersed surface $(S\setminus\Cal{P},i)$ wraps $k$ times about a geodesic which terminates at $\varphi(p_0)$. We observe that
critical points of order $1$ are admitted, even though they are not, strictly speaking, critical points.
\medskip
\noindent We also obtain a converse to this result:
\proclaim{Theorem \nextprocno}
\noindent let $S$ be a surface and let $\Cal{P}\subseteq S$ be a discrete subset of $S$. Let $i:S\setminus\Cal{P}\rightarrow\Bbb{H}^3$ be
an immersion such that $\Sigma=(S\setminus\Cal{P},i)$ is a $k$-surface (and is thus the solution to a Plateau problem). Let
$\overrightarrow{n}:U\Bbb{H}^3\rightarrow\partial_\infty\Bbb{H}^3$ be the Gauss-Minkowski mapping which sends $U\Bbb{H}^3$ to
$\partial_\infty\Bbb{H}^3$. Let $\hat{\mathi}$ be the Gauss lifting of $i$ so that $\varphi=\overrightarrow{n}\circ\hat{\mathi}$ defines
the Plateau problem to which $i$ is the solution. Let $\Cal{H}$ be the holomorphic structure generated over $S\setminus\Cal{P}$ by
the local homeomorphism $\varphi$.
\medskip
\noindent Let $p_0$ be an arbitrary point in $\Cal{P}$, and suppose that $\Sigma$ is asymptotically tubular of order $k$ about $p_0$. Then
there exists a unique holomorphic structure $\tilde{\Cal{H}}$ over $(S\setminus\Cal{P})\munion\left\{p_0\right\}$ and a unique
holomorphic mapping $\tilde{\varphi}:(S\setminus\Cal{P})\munion\left\{p_0\right\}\rightarrow\hat{\Bbb{C}}$ such that $\tilde{\Cal{H}}$ and
$\tilde{\varphi}$ extend $\Cal{H}$ and $\varphi$ respectively. Moreover, $\tilde{\varphi}$ has a critical point of order $k$ at
$p_0$.
\endproclaim
\proclabel{PresentationChIIISurfacesAsymptotiquementTubulaires}
\remark Together, these two theorems provide a complete geometric descrition of solutions to the Plateau problem $(S,\varphi)$ when
$S$ is a compact Riemann surface with a finite number of marked points.
\medskip
\noindent In the first section of this paper, we provide an overview of the definitions and notations that will be used in the sequel.
In the second section, we study the differential geometry of the unitary bundle of a Riemannian manifold, focusing, in particular, on the
canonical contact and complex structures of this bundle. In the third section, we define the Plateau problem, providing various auxiliary
definitions and recalling existing results of \cite{LabA} and \cite{SmiB} which will be required in the sequel. In the fifth section, we
prove theorem \procref{PresentationChIIILimites}. In the sixth section, we study the geometry of the Plateau problem
$(\Bbb{D}^*,z\mapsto z)$, which provides a model for the study of all other cases. In the seventh section, we prove theorem
\procref{PresentationChIIIRevetementsRamifiees}, and in the final section we prove theorem
\procref{PresentationChIIISurfacesAsymptotiquementTubulaires}.
\medskip
\noindent These results provoke the following reflections concerning potential future avenues of research: first, we obtain a homeomorphism between the space of meromorphic mappings over compact Riemann surfaces with a finite number of marked points on the one hand and complete positive pseudo-holomorphic curves immersed in $U\Bbb{H}^3$ with cylindrical ends on the other. These pseudo-holomorphic curves project down to surfaces of constant Gaussian curvature immersed into $\Bbb{H}^3$. Such an equivalence may well permit us to better understand the
structure of either one or both of these two spaces. Second, by integrating primatives of the canonical volume form of $\Bbb{H}^3$ over
these immersed surfaces, one obtains a ``volume'' bounded by these surfaces. If this volume can be shown to be finite, then we would obtain a
new function over the Teichm\"uller space of compact Riemann surfaces with marked points. We would then be interested in the properties of
such a function.
\newhead{Immersed Surfaces - Definitions and Notations}
\newsubhead{Definitions}
\noindent In this section we will review basic definitions from the theory of immersed submanifolds and establish the notations that will be used throughout this article.
\medskip
\noindent Let $M$ be a smooth manifold. An {\emph immersed submanifold\/} is a pair $\Sigma=(S,i)$ where $S$ is a smooth manifold
and $i:S\rightarrow M$ is a smooth immersion. An {\emph immersed hypersurface\/} is an immersed submanifold of codimension $1$.
\medskip
\noindent Let $g$ be a Riemannian metric on $M$. We give $S$ the unique Riemannian metric $i^*g$ which makes $i$
into an isometry. We say that $\Sigma$ is {\emph complete\/} if and only if the Riemannian manifold $(S,i^*g)$ is.
\newsubhead{Normal Vector Fields, Second Fundamental Form, Convexity}
\noindent Let $\Sigma$ be a hypersurface immersed in the Riemannian manifold $M$. There exists a canonical embedding $i_*$ of the
tangent bundle $TS$ of $S$ into the pull-back $i^*TM$ of the tangent bundle of $M$. This embedding may be considered as a section
of $\opEnd(TS, i^*TM)$. We denote by $T\Sigma$ the image of $TS$ under the action of this embedding.
\medskip
\noindent Let us suppose that both $M$ and $S$ are oriented. We define $N\Sigma\subseteq i^*TM$, the normal bundle of $\Sigma$, by:
$$
N\Sigma = T\Sigma^\perp.
$$
\noindent $N\Sigma$ is a one dimensional subbundle of $i^*TM$ from which it inherits a canonical Riemannian metric. Using the orientations of $S$ and $M$, we define the
{\emph exterior unit normal vector field\/}, $\msf{N}_\Sigma\in\Gamma(S,N\Sigma)$, over $\Sigma$ in $M$. This is a global section of $N\Sigma$ which consequently trivialises
this bundle. We define the {\emph Weingarten operator\/}, $A_\Sigma$, which is a section of $\opEnd(TS, T\Sigma)$, by:
$$
A_\Sigma(X) = (i^*\nabla)_X \msf{N}_\Sigma.
$$
\noindent Since there exists a canonical isomorphism (being $i_*$) between $T\Sigma$ and $TS$, we may equally well view $A_\Sigma$ as a section of $\opEnd(TS)$. This section is
self-adjoint with respect to the canonical Riemannian metric over $S$. We thus define the {\emph second fundamental form\/}, $II_\Sigma$, which is a symmetric bilinear form
over $TS$ by:
$$
II_\Sigma(X,Y) = \langle A_\Sigma X, Y \rangle.
$$
\noindent $\Sigma$ is said to be {\emph convex\/} at $p\in S$ if and only if the bilinear form $II_\Sigma$ is either positive or negative
definite at $p$. $\Sigma$ is then said to be {\emph locally convex\/} if and only if it is convex at every point. Through a slight
abuse of language, we will say that $\Sigma$ is {\emph convex\/} in this case. Bearing in mind that the sign of $II_\Sigma$ depends on the sign of $\msf{N}_\Sigma$, which in
turn depends on the choice of orientation of $S$, if $\Sigma$ is convex, then we may choose the orientation of $S$ such that $II_\Sigma$ is positive definite. Consequently, in
the sequel, if $\Sigma$ is convex, then we will assume that $II_\Sigma$ is positive definite.
\newsubhead{Curvature}
\noindent For $i\in\left\{0,...,n\right\}$, let us denote by $\xi_i:\opEnd(\Bbb{R}^n)\rightarrow\Bbb{R}$ the {\emph $i$'th symmetric function\/} over $\opEnd(\Bbb{R}^n)$. Thus,
for all matrices $A$ and for all real $t$:
$$
\opDet(I + tA) = \sum_{i=0}^n\xi_i(A)t^i.
$$
\noindent The group $\opGL(\Bbb{R}^n)$ acts on $\opEnd(\Bbb{R}^n)$ by the {\emph adjoint action\/} which we by $\opAd$. Thus, for all
$A,B\in\opEnd(\Bbb{R}^n)$:
$$
\opAd(A)\cdot B = ABA^{-1}.
$$
\noindent For all $i$, the function $\xi_i$ is invariant under the adjoint action of $\opGL(\Bbb{R}^n)$ on $\opEnd(\Bbb{R}^n)$. For this reason, for any manifold $M$, and for all $i$, the
function $\xi_i$ defines a unique mapping $\tilde{\xi}_i:\Gamma(\opEnd(TM))\rightarrow C^\infty(M)$. For simplicity, we will denote $\tilde{\xi}_i$ also by $\xi_i$.
\medskip
\noindent Let $\Sigma=(S,i)$ be an oriented hypersurface immersed in an oriented Riemannian manifold $M$ of dimension $n+1$. For all $i\in\left\{0,...,n\right\}$, we refer to
the function $\xi_i(A_\Sigma)$ as the {\emph $i$'th higher principal curvature\/} of $\Sigma$. In particuler, $\xi_1(A_\Sigma)$ is the {\emph mean curvature\/} of $\Sigma$,
which is also denoted by $H_\Sigma$, and $\xi_n(A_\Sigma)$ is the {\emph Gaussian curvature\/} of $\Sigma$, which is also denoted by $k_\Sigma$.
\medskip
\noindent In this paper, we study oriented surfaces of constant Gaussian curvature immersed into three dimensional hyperbolic space.
\medskip
\noindent Let $p$ be an arbitrary point in $S$. Let $\Sigma'=(S',i')$ be another oriented immersed hypersurface in $M$. We say that $\Sigma'$ is {\emph tangent} to $\Sigma$ at
$p$ if and only if there exists $p'\in S'$ such that $i(p)=i'(p')$ and:
$$
T_p\Sigma = T_{p'}\Sigma'.
$$
\noindent We call $p'$ a {\emph point of tangency} of $\Sigma'$ on $\Sigma$. For such a pair of points $(p,p')$, we may show that there exists:
\medskip
\myitem{(i)} a neighbourhood $U$ of $p$ in $S$ and a neighbourhood $U'$ of $p'$ in $S'$,
\medskip
\myitem{(ii)} a diffeomorphism $\varphi:(U,p)\rightarrow (U',p')$, and
\medskip
\myitem{(iii)} a function $\lambda:U\rightarrow \Bbb{R}$,
\medskip
\noindent such that, if $\msf{N}:S\rightarrow TM$ is the exterior unit normal vector field over $\Sigma$ in $M$, and if $\opExp:TM\rightarrow M$ is the exponential mapping of
$M$, then, for all $x\in U$:
$$
(i'\circ\varphi)(x) = \opExp(\lambda(x)\msf{N}(x)).
$$
\noindent In otherwords $\Sigma'$ is locally a graph over $\Sigma$ near $p$. Moreover, since $\Sigma'$ is tangent to $\Sigma$ at $p$, we obtain:
$$
d\lambda(p) = 0.
$$
\noindent If $\varphi'$ is another diffeomorphism defined in a neighbourhood of $p$ such that $\varphi'(p)=p'$ and if $\lambda'$ is another function defined in a neighbourhood of
$p$ such that $(i'\circ\varphi')(x) = \opExp(\lambda'(x)\msf{N}(x))$ for all $x$ in a neighbourhood of $p$, then the pairs of functions $(\varphi,\lambda)$ and
$(\varphi',\lambda')$ coincide in a neighbourhood of $p$.
\medskip
\noindent If $\Sigma'$ is tangent to $\Sigma$ at $p$, then we say that $\Sigma'$ is an {\emph exterior tangent} (resp. {\emph interior tangent}) to $\Sigma$ at $p$ if and only
if $\lambda\geqslant 0$ (resp. $\lambda\leqslant 0$) in a neighbourhood of $p$. We now obtain the following {\emph weak geometric maximum principal}:
\proclaim{Lemma \nextprocno {\sl Weak Geometric Maximum Prinicipal.}}
\noindent Let $M$ be an oriented manifold. Let $\Sigma=(S,i)$ and $\Sigma'=(S',i')$ be two oriented immersed hypersurfaces in $M$. Let $p$ be a point in $S$ and
suppose that $\Sigma'$ is an exterior tangent to $\Sigma$ at $p$. Let $p'\in S'$ be a point of exterior tangency of $\Sigma'$ on $\Sigma$. If $k_\Sigma(p)$ and
$k_{\Sigma'}(p')$ denote the Gaussian curvatures of $\Sigma$ at $p$ and $\Sigma'$ at $p'$ respectively, we obtain:
$$
k_{\Sigma}(p) \leqslant k_{\Sigma'}(p')
$$
\endproclaim
\proclabel{LemmaGeometricMaximumPrincipal}
\noindent A proof of this result may be found in \cite{LabA}. An analogous result exists when $\Sigma'$ is an interior tangent to $\Sigma$ at $p$.
\newsubhead{Haussdorf Convergence}
\noindent In the sequel, we will make use of the notion of Haussdorf convergence of sequences of compact sets contained withing a given metrisable space. The following lemmata
will permit us to better understand the nature of the Haussdorf topology. First, we recall a classical result which tells us that the Haussdorf topology of a compact metrisable
space does not depend on the metric chosen over that space:
\proclaim{Lemma \nextprocno}
\noindent Let $X$ be a compact metrisable space. Let $g_1$ and $g_2$ be two metrics over $X$ compatible with the topology of $X$. Let $\suite{A}{n}, A_0$ be compact subsets
of $X$. The sequence $\suite{A}{n}$ converges to $A_0$ in the $g_1$-Haussdorf topology if and only if it converges to $A_0$ in the
$g_2$-Haussdorf topology.
\endproclaim
\noindent In particular, the Haussdorf topology of $\Bbb{H}^3\munion\partial_\infty\Bbb{H}^3$ is well defined. Next, we have a result concerning the relationship between the
Haussdorf topology and the topology of uniform convergence for homoemorphisms of a given compact metric space:
\proclaim{Lemma \nextprocno}
\noindent Let $(X,d)$ be a compact metric space. Let $\suite{Y}{n}, Y_0\subset X$ be subsets of $X$ such that $\suite{Y}{n}$ converges to $Y_0$ in the Haussdorf topology. Let
$\suite{\alpha}{n},\alpha_0$ be homeomorphisms of $X$ such that $(\alpha_n)_{n\in\Bbb{N}}$ converges to $\alpha_0$ in the compact-open topology (i.e. the topology of uniform
convergence). The sequence $(\alpha_n(Y_n))_{n\in\Bbb{N}}$ converges to $\alpha_0(Y_0)$ in the Haussdorf topology.
\endproclaim
\proclabel{LemmeConvergenceDesEnsemblesEtDesFonctions}
\noindent Finally, we have a result concerning the intersections of two sequences of compact sets that converge:
\goodbreak
\proclaim{Lemma \nextprocno}
\noindent Let $(X,d)$ be a compact metric space. Let $\suite{A}{n},A_0\subseteq X$ and $\suite{B}{n},B_0\subseteq X$ be compact sets such that $\suite{A}{n}$ and $\suite{B}{n}$
converge to $A_0$ and $B_0$ respectively in the Haussdorf topology. If, for all $n$:
$$
A_n\minter B_n\neq\emptyset.
$$
\noindent then:
$$
A_0\minter B_0\neq\emptyset.
$$
\endproclaim
\proclabel{LemmeLIntersectionDesLimites}
\noindent Proofs of lemmata \procref{LemmeConvergenceDesEnsemblesEtDesFonctions} and \procref{LemmeLIntersectionDesLimites} may be found in appendix $A$ of \cite{SmiE}.
\newsubhead{Pointed Manifolds, Convergence}
\noindent In the sequel, we will use the concept of Cheeger/Gromov convergence for complete pointed immersed submanifolds.
\medskip
\noindent A {\emph pointed Riemannian manifold\/} is a pair $(M,p)$ where $M$ is a Riemannnian manifold and $p$ is a point in $M$. If $(M,p)$ and $(M',p')$ are pointed
manifolds then a {\emph morphism\/} (or {\emph mapping\/}) from $(M,p)$ to $(M',p')$ is a (not necessarily even continuous) function from $M$ to $M'$ which sends $p$ to $p'$ and
is of type $C^\infty$ in a neighbourhood of $p$. The family of pointed manifolds along with these morphisms forms a category. In this section, we will discuss a notion of convergence
for this familly. It should be borne in mind that even though this family is not a set, we may consider it as such. Indeed, since every manifold may be plunged into an infinite
dimensional real vector space, we may discuss, instead, the equivalent family of pointed finite dimensional submanifolds of this vector space, and this is a set.
\medskip
\noindent Let $(M_n,p_n)_{n\in\Bbb{N}}$ be a sequence of complete pointed Riemannian manifolds. For all $n$, we denote by $g_n$ the Riemannian metric over $M_n$. We say that the
sequence $(M_n,p_n)_{n\in\Bbb{N}}$ {\emph converges\/} to the complete pointed manifold $(M_0,p_0)$ in the {\emph Cheeger/Gromov topology\/} if and only if, for all $n$, there
exists a mapping $\varphi_n:(M_0,p_0)\rightarrow (M_n,p_n)$ such that, for every compact subset $K$ of $M_0$, there exists $N\in\Bbb{N}$ such that for all $n\geqslant N$:
\medskip
\myitem{(i)} the restriction of $\varphi_n$ to $K$ is a $C^\infty$-diffeomorphism onto its image, and
\medskip
\myitem{(ii)} if we denote by $g_0$ the Riemannian metric over $M_0$, then the sequence of metrics $(\varphi_n^*g_n)_{n\geqslant N}$ converges to $g_0$ in the $C^\infty$ topology
over $K$.
\medskip
\noindent We refer to the sequence $(\varphi_n)_{n\in\Bbb{N}}$ as a sequence of {\emph convergence mappings\/} of the sequence $(M_n)_{n\in\Bbb{N}}$ with respect to the
limit $(M_0,p_0)$. The convergence mappings are trivially not unique. However, two sequences of convergence mappings $(\varphi_n)_{n\in\Bbb{N}}$ and $(\varphi'_n)_{n\in\Bbb{N}}$
are equivalent in the sense that there exists an isometry $\phi$ of $(M_0,p_0)$ such that, for every compact subset $K$ of $M_0$, there exists $N\in\Bbb{N}$ such that for
$n\geqslant N$:
\medskip
\myitem{(i)} the mapping $(\varphi_n^{-1}\circ\varphi_n')$ is well defined over $K$, and
\medskip
\myitem{(ii)} the sequence $(\varphi_n^{-1}\circ\varphi_n')_{n\geqslant N}$ converges to the $\phi$ in the $C_\oploc^\infty$ topology over $K$.
\medskip
\noindent One may verify that this mode of convergence does indeed arise from a topological structure over the space of complete pointed manifolds. Moreover, this topology
is Haussdorf (up to isometries).
\medskip
\noindent Most topological properties are unstable under this limiting process. For example, the limit of a sequence of simply connected manifolds is not necessarily
simply connected. On the other hand, the limit of a sequence of surfaces of genus $k$ is a surface of genus at most $k$ (but quite possibly with many holes).
\medskip
\noindent Let $M$ be a complete Riemannian manifold. A {\emph pointed immersed submanifold\/} in  $M$ is a pair $(\Sigma,p)$ where $\Sigma=(S,i)$ is an immersed submanifold in
$M$ and $p$ is a point in $S$.
\medskip
\noindent Let $(\Sigma_n,p_n)_{n\in\Bbb{N}}=(S_n,p_n,i_n)_{n\in\Bbb{N}}$ be a sequence of complete pointed immersed submanifolds in $M$. We say that
$(\Sigma_n,p_n)_{n\in\Bbb{N}}$ {\emph converges\/} to $(\Sigma_0,p_0)=(S_0,p_0,i_0)$ in the {\emph Cheeger/Gromov topology\/} if and only if $(S_n,p_n)_{n\in\Bbb{N}}$ converges
to $(S_0,p_0)$ in the Cheeger/Gromov topology, and, for every sequence $(\varphi_n)_{n\in\Bbb{N}}$ of convergence mappings with respect to this limit, and for every compact
subset $K$ of $S_0$, the sequence of functions $(i_n\circ\varphi_n)_{n\geqslant N}$ converges to the function $(i_0\circ\varphi_0)$ in the $C^\infty$ topology over $K$.
\medskip
\noindent As before, this mode of convergence arises from a topological structure over the space of complete immersed submanifolds. Moreover, this topology is
Haussdorf (up to isometries).
\newsubhead{``Common Sense'' Lemmata}
\noindent In order to make good use of the concept of Cheeger/Gromov convergence, it is helpful to recall some basic lemmas concerning the topological properties of functions
acting on open subsets of $\Bbb{R}^n$. The results that follow are essentially formal expressions of ``common sense''. To begin with, we recall a result concerning the inverses
of a sequence of functions that converges:
\proclaim{Lemma \nextprocno}
\noindent Let $\Omega\subseteq\Bbb{R}^n$ be an open set. Let $\suite{f}{n},f_0:\Omega\rightarrow\Bbb{R}^n$ be such that,
for every $n$, the function $f_i$ is a homeomorphism onto its image. Let $\Omega'$ be the image of $\Omega$ under the action of $f_0$.
If $\suite{f}{n}$ converges towards $f_0$ locally uniformly in $\Omega$, then the sequence
$(f_n^{-1})_{n\in\Bbb{N}}$ converges towards $f_0^{-1}$ locally uniformly in $\Omega'$.
\medskip
\noindent To be precise, for every compact subset $K$ in $\Omega'$, there exists $N\in\Bbb{N}$ such that, for every $n\geqslant N$,
the set $K$ is contained within $f_n(\Omega)$ and $(f_n^{-1})_{n\geqslant N}$ converges towards $f_0^{-1}$ uniformly over $K$.
\medskip
\noindent Moreover, if every $f_n$ is of type $C^m$ and if $\suite{f}{n}$ converges to $f_0$ in the $C^m_{\oploc}$ topology, then $(f_n^{-1})_{n\in\Bbb{N}}$ also converges towards
$f_0^{-1}$ in the $C^m_{\oploc}$ topology.
\endproclaim
\proclabel{LemmeLimitesDInverses}
\noindent We recall a result concerning the injectivity of the limit:
\proclaim{Lemma \nextprocno}
\noindent Let $\Omega\subseteq\Bbb{R}^n$ be an open set. Let $\suite{f}{n},f_0:\Omega\rightarrow\Bbb{R}^n$ be such that for every $n>0$, the function $f_n$ is a homeomorphism
onto its image. If $(f_n)_{n\in\Bbb{N}}$ tends towards $f_0$ locally uniformly, and if, moreover, $f_0$ is a local homeomorphism, then $f_0$ is injective.
\endproclaim
\proclabel{LemmeInjectiviteDeLaLimite}
\noindent We recall a converse of this result for $C^2$ functions:
\proclaim{Lemma \nextprocno}
\noindent Let $\Omega\subseteq\Bbb{R}^n$ be an open set. Let $\suite{f}{n},f_0:\Omega\rightarrow\Bbb{R}^n$ be $C^2$ functions such that
$\suite{f}{n}$ converges towards $f_0$ in the $C^2_{\oploc}$ topology. If $f_0$ is a diffeomorphism onto its image, then for every compact subset
$K$ in $\Omega$, there exists $N\in\Bbb{N}$ such that, for $n\geqslant N$, the restriction of $f_n$ to $K$ is injective.
\endproclaim
\proclabel{LemmeInjectiviteLorsquOnSApprocheDuLimite}
\noindent Finally, we have a result concerning the images of a sequence of functions:
\proclaim{Lemma \nextprocno}
\noindent Let $\Omega\subseteq\Bbb{R}^n$ be an open set. Let $\suite{f}{n}:\Omega\rightarrow\Bbb{R}^n$ be such that for every $n$, the function $f_n$ is
a homeomorphism onto its image. If there exists a local homeomorphism $f_0:\Omega\rightarrow\Bbb{R}^n$ such that $\suite{f}{n}$ converges to $f_0$ locally uniformly, then, for
every compact subset $K\subseteq f_0(\Omega)$, there exists $N\in\Bbb{N}$ such that for $n\geqslant N$:
$$
K\subseteq f_n(\Omega).
$$
\endproclaim
\proclabel{LemmeCompacteDansLImageDeLaLimite}
\noindent The interested reader may find a discussion and proofs of these results in the appendix A of \cite{SmiE}.
\newhead{The Unitary Bundle of a Riemannian Manifold}
\newsubhead{Geometric Structures Over $TM$}
\noindent Let $M$ be a Riemannian manifold. We define $\pi:TM\rightarrow M$ to be the canonical projection of the tangent space of $M$ onto $M$. We denote by
$HTM\subseteq TTM$ the {\emph horizontal bundle\/} of the Levi-Civita covariant derivative of $M$. We denote by $VTM\subseteq TTM$ the {\emph vertical bundle\/}
over $TM$. To be precise, $VTM$ is defined to be the kernel of the projection $\pi$ within $TTM$. The tangent bundle of $TM$ is the direct sum of these two subbundles:
$$
TTM = HTM\oplus VTM.
$$
\noindent Each of $HTM$ and $VTM$ is canonically isomorphic to $\pi^*TM$. We denote by $i_H$ (resp. $i_V$), which is a section of $\opEnd(HTM,\pi^*TM)$ (resp. $\opEnd(VTM,\pi^*TM)$), the canonical isomorphism sending $HTM$ (resp. $VTM$) to $\pi^*TM$. We obtain the following isomorphism:
$$
i_H\oplus i_V:TTM\rightarrow\pi^*TM\oplus\pi^*TM.
$$
\noindent For every pair of vector fields $X,Y\in\Gamma(M,TM)$ over $M$ we define the vector field $\left\{X,Y\right\}$ over $TM$ by:
$$
(i_H\oplus i_V)(\left\{X,Y\right\}) = (\pi^* X, \pi^* Y).
$$
\noindent Trivially, every vector field over $TM$ may be expressed (at least locally) in terms of a linear combination of such vector fields. In the same way, for
a given point $p\in M$ in $M$ and for a given triplet of vectors $X,Y,q\in T_pM$ over $p$, we may define $\left\{X,Y\right\}_q\in T_qTM$ by:
$$
(i_H\oplus i_V)_q\left\{X,Y\right\}_q = (\pi^*_qX, \pi^*_q Y).
$$
\noindent Finally, for a given vector field $X$ over $M$, we may define $X^H$ and $X^V$ by:
$$\matrix
X^H \hfill&=\left\{X,0\right\}, \hfill\cr
X^V \hfill&=\left\{0,X\right\}. \hfill\cr
\endmatrix$$
\newsubhead{Geometric Structures Over $UM$}
\noindent For $M$ a Riemannian manifold, we define $UM$, the {\emph unitary bundle\/} over $M$ by:
\headlabel{HeadDifferentialGeometryInTTM}
$$
UM = \left\{X\in TM\ |\ \|X\| = 1\right\}.
$$
\noindent We define the {\emph tautological vector fields} $T^H$ and $T^V$ over the tangent space $TM$ to $M$ by:
$$\matrix
T^H(q) \hfill&=\left\{q,0\right\}_q, \hfill\cr
T^V(q) \hfill&=\left\{0,q\right\}_q. \hfill\cr
\endmatrix$$
\noindent Let $i:UM\rightarrow TM$ be the canonical embedding. Let $HUM$ (resp. $VUM$) be the restriction of $HTM$ (resp. $VTM$) to $UM$:
$$\matrix
HUM \hfill&=i^*HTM, \hfill\cr
VUM \hfill&=i^*VTM. \hfill\cr
\endmatrix$$
\noindent The section $i^*T^H$ (resp. $i^*T^V$) is nowhere vanishing. It consequently defines a one dimensional subbundle of $HUM$ (resp. $VUM$). In order
to simplify the notation we will also denote this section by $T^H$ (resp. $T^V$). We denote the one dimensional subbundle that it generates by $\langle q\rangle_H$
(resp. $\langle q\rangle_V$). We define the subbundles $\langle q\rangle_H^\perp$ and $\langle q\rangle_V^\perp$ of $HUM$ and $VUM$ respectively by:
$$\matrix
\langle q\rangle_H^\perp \hfill&=\langle T^H \rangle^\perp, \hfill\cr
\langle q\rangle_V^\perp \hfill&=\langle T^V \rangle^\perp. \hfill\cr
\endmatrix$$
\noindent Since parallel transport preserves the length of vectors and thus sends $UM$ onto itself, the immersion $i$ induces the following
isomorphism of vector bundles:
$$
i_*:TUM\rightarrow HUM\oplus \langle q\rangle_V^\perp.
$$
\noindent In order to simplify our notation, we consider $HUM$, $\langle q\rangle_H^\perp$ and $\langle q\rangle_V^\perp$ as subbundles of $TUM$. In particular, we define $WUM$
by:
$$
WUM = \langle q\rangle_H^\perp\oplus\langle q\rangle_V^\perp.
$$
\noindent The subbundle $WUM$ defines, in fact, a contact structure over $UM$, and we will consequently refer to it as the {\emph contact bundle\/} over $UM$.
In the sequel we will denote the bundles $HUM$, $VUM$, $WUM$ etc. by $H$, $V$, $W$ etc.
\noindent For $k>0$, we write $\nu=\sqrt{k}$ and we define the metric $g^\nu$ over $TTM$ such that, for every pair of vector fields $X,Y\in\Gamma(M,TM)$ over $M$ we have:
$$
g^\nu(\left\{X,Y\right\},\left\{X,Y\right\}) = \langle X, X\rangle + \nu^{-2}\langle Y, Y\rangle.
$$
\noindent We denote also by $g^\nu$ the metric induced over $UM$ by $g^\nu$ and the canonical embedding of $UM$ into $TM$.
\medskip
\noindent From now on, we will suppose that $M$ is oriented and three dimensional. This allows us to canonically identify $TM$ and $TM \wedge TM$ and consequently to define a
vector product $\wedge$ over $TM$. We then define the canonical {\emph complex structures\/} $J^H$ (resp. $J^V$) over $\langle q\rangle_H^\perp$ (resp.
$\langle q\rangle_V^\perp$) such that for every $q\in UM$ and for every vector $X$ orthogonal to $q$:
$$\matrix
J^H_q\left\{X,0\right\}_q \hfill&=\left\{q\wedge X, 0\right\}_q, \hfill\cr
J^V_q\left\{0,X\right\}_q \hfill&=\left\{0, q\wedge X\right\}_q. \hfill\cr
\endmatrix$$
\noindent In order to simplify notation we refer to both $J^H$ and $J^V$ by $J$. We define the isomorphism $j:HTM\rightarrow VTM$ by:
$$
j = i_V^{-1}\circ i_H.
$$
\noindent This isometry sends $\langle q\rangle_H$ onto $\langle q\rangle_V$ and consequently $\langle q\rangle_H^\perp$ onto $\langle q\rangle_V^\perp$. Moreover, we
trivially obtain the following commutative diagram:
$$
\commdiag{
\langle q\rangle_H^\perp &\mapright^j &\langle q\rangle_V^\perp \cr
\mapdown_J& &\mapdown_J\cr
\langle q\rangle_H^\perp &\mapright^j &\langle q\rangle_V^\perp}
$$
\noindent We define the {\emph complex structure\/} $J^\nu$ over $W=\langle q\rangle_H^\perp\oplus\langle q\rangle_V^\perp$ by:
$$
J^\nu = \pmatrix 0 & j^{-1}\nu^{-1}J \cr j\nu J & 0 \cr \endpmatrix.
$$
\noindent $J^\nu$ is compatible with the metric $g^\nu$ over $W$. In the sequel, in order to simplify the notation, we identify $\langle q\rangle_H^\perp$ and
$\langle q\rangle_V^\perp$ through the isomorphism $j$, and we thus write:
$$
J^\nu = \pmatrix 0 & \nu^{-1}J \cr \nu J & 0 \cr \endpmatrix.
$$
\noindent By composing this form with the orthogonal projection of $TUM$ onto $W$, we may extend this form to one defined on $TUM$.
\medskip
\noindent Let $q$ be a point in $UM$. Let $\Sigma\subseteq W_q$ be a plane in $W_q$. We say that $\Sigma$ is the graph of the matrix $A$ over $\langle q\rangle_H^\perp$
if and only if:
$$
\Sigma = \left\{ \left\{V,AV\right\} | V\in q^\perp \right\}.
$$
\noindent We say that the plane $\Sigma$ is {\emph k-complex\/} if and only if it is stable under the action of $J^\nu$. In this case, if it is the graph of a matrix $A$,
we may trivially show that $A$ is symmetric. The plane $\Sigma$ is then said to be {\emph positive\/} if and only if it is the graph of a positive definite matrix.
\newsubhead{Holomorphic Curves, k-Surfaces}
\noindent Let $M$ be a compact oriented three dimensional Riemannian manifold and let $\Sigma=(S,i)$ be a convex hypersurface immersed in $M$. Let $\mathsf{N}_\Sigma$ be
the normal exterior vector field to $\Sigma$. We define the {\emph Gauss lifting\/} $\hat{\Sigma}=(S,\hat{\mathi})$ of $\Sigma$ by:
$$
(S,\hat{\mathi}) = (S,\mathsf{N}_\Sigma).
$$
\noindent For $k>0$, we say that $\Sigma$ is a {\emph k-surface\/} if and only if $\hat{\Sigma}$ is complete and the Gaussian curvature of $\Sigma$ is always equal to $k$.
\medskip
\noindent We say that $\hat{\Sigma}$ is a {\emph k-holomorphic curve\/} if and only if all its tangent planes are k-complex planes, and we say that it is
{\emph positive\/} if and only if all its tangent planes are positive.
\medskip
\noindent These concepts are related by the following elementary result:
\proclaim{Lemma \nextprocno}
\noindent Let $M$ be an oriented three dimensional Riemannian manifold. Let $\Sigma=(S,i)$ be a convex hypersurface immersed in $M$. $\Sigma$ is a k-surface if and only if $\hat{\Sigma}$ use a complete positive k-holomorphic curve.
\endproclaim
\proof See, for example \cite{LabB}.\qed
\medskip
\noindent We now consider the case where $M=\Bbb{H}^3$. Let $\overrightarrow{n}$ be the Gauss-Minkowski mapping which sends $U\Bbb{H}^3$ to
$\partial_\infty\Bbb{H}^3\cong\hat{\Bbb{C}}$. For $\hat{\Sigma}=(S,\hat{\mathi})$ a k-holomorphic curve in $U\Bbb{H}^3$, we define $\varphi:S\rightarrow\hat{\Bbb{C}}$ by:
$$
\varphi = \overrightarrow{n}\circ\hat{\mathi}.
$$
\noindent Let $\Cal{H}$ be the canonical holomorphic structure over $\hat{\Bbb{C}}$. We obtain the following result:
\proclaim{Lemma \nextprocno}
\noindent Let $\hat{\Sigma}=(S,\hat{\mathi})$ be a positive k-holomorphic curve in $U\Bbb{H}^3$. Let $\overrightarrow{n}$ be the Gauss-Minkowski mapping which sends $U\Bbb{H}^3$
to $\partial_\infty\Bbb{H}^3\cong\hat{\Bbb{C}}$. Let $\Cal{H}$ be the canonical conformal structure over $\hat{\Bbb{C}}$. Let $\Cal{H}'$ be the conformal structure
generated over $S$ by $\hat{\mathi}^*g^\nu$ and the canonical orientation of $S$. The two structures $\Cal{H}'$ and $\varphi^*\Cal{H}'$ are quasiconformally equivalent.
\endproclaim
\proclabel{LemmaCSQCEquivalence}
\proof See \cite{SmiB}. \qed
\newhead{The Plateau Problem}
\newsubhead{Definitions}
\noindent A {\emph Hadamard manifold} is a complete, connected and simply connected manifold of negative sectional curvature. The manifold $\Bbb{H}^3$ is an example of a
$3$-dimensional Hadamard manifold. In \cite{LabA}, Labourie studies the {\emph Plateau problem} for constant Gaussian curvature
hypersurfaces immersed in a three dimensional Hadamard manifold $M$. In the language of this paper, a Plateau problem is a pair $(S,\varphi)$
where $S$ is a Riemann surface and $\varphi:S\rightarrow\partial_\infty\Bbb{H}^3$ is a local conformal mapping. A solution to
this Plateau problem is an immersion $i:S\rightarrow\Bbb{H}^3$ such that the immersed hypersurface $(S,i)$ is a
k-surface and, if we denote by $\hat{\mathi}$ the Gauss lifting of $i$ and by $\overrightarrow{n}$ the Gauss-Minkowski mapping, then:
$$
\overrightarrow{n}\circ\hat{\mathi} = \varphi.
$$
\newsubhead{Tubes, Tubular Surfaces, Asymptotically Tubular Surfaces}
\noindent In this section we will define tubes about geodesics which, as will be shown in the following sections, play a special role in the study of k-surfaces.
\medskip
\noindent For $\Gamma$ a geodesic in $\Bbb{H}^3$, we define $N_\Gamma$ to be the normal bundle over $\Gamma$ in $U\Bbb{H}^3$:
$$
N_\Gamma = \left\{ n_p\in U\Bbb{H}^3{ s.t. }p\in\Gamma, n_p\perp T_p\Gamma\right\}.
$$
\noindent A {\emph tube} about $\Gamma$ is a pair $T=(S,\hat{\mathi})$ where $S$ is a complete surface and $\hat{\mathi}:S\rightarrow N_\Gamma$ is a covering map.
Since $N_\Gamma$ is conformally equivalent to $S^1\times\Bbb{R}$, where $S^1$ is the circle of radius $1$ in $\Bbb{C}$, we may assume either that
$S=S^1\times\Bbb{R}$ or that $S=\Bbb{R}\times\Bbb{R}$. In the former case $\hat{\mathi}$ is a covering map of finite order, and, if $k$ is the order of $\hat{\mathi}$, then
we say that the tube $T$ is a {\emph tube of order $k$}. The application $\hat{\mathi}$ is then unique up to vertical translations and horizontal rotations of $S^1\times\Bbb{R}$.
In the latter case, we say that the tube $T$ is a {\emph tube of infinite order}. The application $\hat{\mathi}$ is then unique up to translations of $\Bbb{R}\times\Bbb{R}$. In
either case, we call the point $(0,0)$ the {\emph origin} of the tube $T$. In the sequel, we will only be interested in tubes of finite order.
\medskip
\noindent Let $T=(S^1\times\Bbb{R},\hat{\mathi})$ be a tube of order $k$. We define the fields $\partial_\theta$ and $\partial_t$ over $T$ by:
$$\matrix
\partial_\theta(e^{i\theta},t) \hfill&= [\phi\mapsto(e^{i\theta+i\phi},t)],\hfill\cr
\partial_t(e^{i\theta},t) \hfill&= [s\mapsto(e^{i\theta},t+s)].\hfill\cr
\endmatrix$$
\noindent Using the definition of $g^\nu$, we find that every fibre of $N_\Gamma$ is a circle of length $2\pi\nu^{-1}$. Consequently, since $\hat{\mathi}$ is a covering
map of order $k$, and since $S^1$ is of length $2\pi$, it follows by homogeneity that:
$$
\|T\hat{\mathi}\cdot\partial_\theta\| = k\nu^{-1}.
$$
\noindent Since $\hat{\mathi}$ is locally conformal, we obtain:
$$
\|T\hat{\mathi}\cdot\partial_t\| = k\nu^{-1}.
$$
\noindent Let $\opExp:TU\Bbb{H}^3\rightarrow U\Bbb{H}^3$ be the exponential mapping over $U\Bbb{H}^3$. Let $NN_\Gamma$ be the normal bundle over $N_\Gamma$ in $TU\Bbb{H}^3$.
Let $T=(S^1\times\Bbb{R},\hat{\mathi})$ be a tube of order $k$ about $\Gamma$. We define the normal bundle $NT$ over $T$ by:
$$
NT = \hat{\mathi}^*NN_\Gamma.
$$
\noindent For $r\in\Bbb{R}$ we define $T_r$ by:
$$
T_r = (S^1\times (-r,r),\hat{\mathi}).
$$
\noindent We define $NT_r$ to be the restriction of $NT$ to the set $S^1\times (-r,r)$. Let $(\hat{\Sigma},p)$ be a pointed immersed surface in $U\Bbb{H}^3$. We say that
$(\hat{\Sigma},p)$ is a {\emph graph} over $T$ of half length $r$ if and only if there exists:
\medskip
\myitem{(i)} a neighbourhood $\Omega$ of $S$ about $p$,
\medskip
\myitem{(ii)} a diffeomorphism $\varphi:S^1\times (-r,r)\rightarrow\Omega$, and
\medskip
\myitem{(iii)} a section $\lambda\in\Gamma(S^1\times(-r,r),NT_r)$,
\medskip
\noindent such that $\varphi(0,0)=p$ and:
$$
\opExp\circ\lambda = \hat{\mathi}\circ\varphi.
$$
\noindent We call $\varphi$ a {\emph graph diffeomorphism} of $(\hat{\Sigma},p)$ over $T_r$ and we call $\lambda$ a {\emph graph function} of $(\hat{\Sigma},p)$ over $T_r$.
\medskip
\noindent For $\epsilon\in\Bbb{R}^+$, we define $N_\epsilon N_\Gamma$ by:
$$
N_\epsilon N_\Gamma = \left\{ v_p\in NN_\Gamma\text{ s.t. } \|v_p\|\leqslant\epsilon\right\}.
$$
\noindent Since $U\Bbb{H}^3$ is homogeneous, there exists $\epsilon\in\Bbb{R}^+$ independant of $\Gamma$ such that the restriction of $\opExp$ to $N_\epsilon N_\Gamma$ is a
diffeomorphism onto its image.
\medskip
\noindent It follows that if $S$ is a graph over a tube of order $k$ of half length $r$ with graph diffeomorphism $\varphi$ and graph function $\lambda$, and if
$\|\lambda\|<\epsilon$, then $\lambda$ and $\varphi$ are unique.
\medskip
\noindent We define the {\emph upper half tube\/} $T_+$ of $T$ by:
$$
T_+ = (S^1\times (0,\infty),\hat{\mathi}).
$$
\noindent We define $NT_+$ to be the restriction of $NT$ to the set $S^1\times (0,\infty)$.
\medskip
\noindent Let $S$ be surface and let $p$ be a point in $S$. Let $\hat{\mathi}:S\setminus\left\{p\right\}\rightarrow U\Bbb{H}^3$ be an immersion. We define the immersed
surface $\hat{\Sigma}$ by:
$$
\hat{\Sigma} = (S\setminus\left\{p\right\},\hat{\mathi}).
$$
\noindent We say that $\hat{\Sigma}$ is a {\emph graph} over $T_+$ near $p$ if and only if there exists:
\medskip
\myitem{(i)} a neighbourhood $\Omega$ of $S$ about $p$,
\medskip
\myitem{(ii)} a diffeomorphism $\varphi:S^1\times(0,\infty)\rightarrow\Omega\setminus\left\{p\right\}$, and
\medskip
\myitem{(iii)} a section $\lambda\in\Gamma(S^1\times(0,\infty),NT_+)$,
\medskip
\noindent such that $\varphi(e^{i\theta},t)$ tends to $p$ as $t$ tends to $\infty$ and:
$$
\opExp\circ\lambda = \hat{\mathi}\circ\varphi.
$$
\noindent As before, we call $\varphi$ a {\emph graph diffeomorphism} of $\hat{\Sigma}$ over $T_+$ and we call $\lambda$ a {\emph graph function} of $\hat{\Sigma}$
over $T_+$. Similarly, if $\|\lambda\|<\epsilon$, then $\lambda$ and $\varphi$ are unique up to composition with an affine transformation of $S^1\times\Bbb{R}$.
\medskip
\noindent There exists a canonical trivialisation $\tau:NT\rightarrow (S_1\times\Bbb{R})\times\Bbb{R}^3$ which is unique up to composition by an endomorphism in $S0(3)$.
Consequently, we may interpret a graph funtion $\lambda$ as a function on a subset of $S^1\times\Bbb{R}$ taking values in $\Bbb{R}^3$. We say that $\hat{\Sigma}$ is
{\emph asymptotically tubular of order $k$} about $p$ if and only if there exists a tube $T$ of order $k$ such that:
\medskip
\myitem{(i)} $\hat{\Sigma}$ is a graph over $T_+$, and
\medskip
\myitem{(ii)} if $\lambda$ is a graph function of $\hat{\Sigma}$ over $T_+$, then, for all $p\in\Bbb{N}$:
$$
\|D^p\lambda(e^{i\theta},t)\|\rightarrow 0\text{ as }t\rightarrow +\infty.
$$
\newsubhead{The Space of Solutions}
\noindent We define $\Cal{L}$ to be the set of Gauss liftings of pointed k-surfaces in $\Bbb{H}^3$:
$$
\Cal{L} = \left\{(\hat{\Sigma},p)|\Sigma\text{ is a k-surface in }\Bbb{H}^3, p\in\hat{\Sigma} \right\}.
$$
\noindent We define $\Cal{L}_\infty$ to be the set of pointed tubes in $U\Bbb{H}^3$:
$$
\Cal{L}_\infty = \left\{(T,p)| T\text{ is a tube about a geodesic }\gamma\text{ in }\Bbb{H}^3,p\in T\right\}.
$$
\noindent We define $\overline{\Cal{L}}$ to be the union of these two sets:
$$
\overline{\Cal{L}} = \Cal{L}\munion\Cal{L}_\infty.
$$
\noindent The justification for this notation will become clear presently. Following \cite{LabA}, we will refer to the surfaces contained in $\Cal{L}_\infty$ as {\emph curtain
surfaces\/}.
\medskip
\noindent In \cite{SmiB}, we proved the existence of solutions to Plateau problems of hyperbolic type. We quote theorem $1.1$ of this paper:
\proclaim{Theorem \nextprocno\ {\bf [Smith, 2004]}\ {\sl Hyperbolic Existence Theorem.}}
\noindent Let $\varphi:\Bbb{D}\rightarrow\hat{\Bbb{C}}$ be locally conformal. Then, for every $k\in (0,1)$, there exists a unique solution $i_k:\Bbb{D}\rightarrow\Bbb{H}^3$
to the Plateau problem $(\Bbb{D},\varphi)$.
\endproclaim
\proclabel{ThmSmithExistence}
\noindent We will also require the following result of Labourie concerning the relationship between two solutions to the Plateau problem. We quote the
theorem $7.2.1$ of \cite{LabA}:
\proclaim{Theorem \nextprocno\ {\bf [Labourie, 2000]}\ {\sl Structure of Solutions.}}
\noindent There exists at most one solution to the Plateau problem. Moreover, let $\Sigma=(S,i)$ be a k-surface
in $\Bbb{H}^3$ and let $\hat{\Sigma}=(S,\hat{\mathi})$ be its Gauss lifting. Define $\varphi$ by
$\varphi=\overrightarrow{n}\circ\hat{\mathi}$ such that $i$ is the solution to the Plateau problem $(S,\varphi)$.
Let $\Omega$ by an open subset of $S$. There exists a solution to the Plateau problem $(\Omega,\varphi)$ which is
a graph over $\Omega$ in the extension of $\Sigma$. In otherwords, there exists $f:\Omega\rightarrow\Bbb{R}^+$ such
that the solution coincides (up to reparametrisation) with $(\Omega,\opExp(fN_\Sigma))$.
\endproclaim
\proclabel{LemmeLabNature}
\noindent Finally, as a compactness result, we use the principal result of \cite{LabB}, which translates into
our framework as follows:
\proclaim{Theorem \nextprocno\ {\bf [Labourie, 2000]}\ {\sl Compactness.}}
\noindent Let $(\Sigma_n)_{n\in\Bbb{N}}=(S_n,i_n)_{n\in\Bbb{N}}$ be a sequence of k-surfaces in $\Bbb{H}^3$ and, for every $n$, let
$\hat{\Sigma}_n=(S_n,\hat{\mathi}_n)$ be the Gauss lifting of $\Sigma_n$. For every $n$, let $p_n\in S_n$ be an
arbitrary point of $S_n$. If there exists a compact subset $K\subseteq U\Bbb{H}^3$ such that $\hat{\mathi}_n(p_n)\in K$ for every $n$, then there exists
$(\hat{\Sigma}_0,p_0)\in\overline{\Cal{L}}=\Cal{L}\munion\Cal{L}_\infty$ such that, after extraction of a
subsequence, $(\hat{\Sigma}_n,p_n)_{n\in\Bbb{N}}$ converges to $(\hat{\Sigma}_0,p_0)$ in the Cheeger/Gromov topology.
\endproclaim
\proclabel{ThmLabCompacite}
\remark One may also obtain this theorem as a special case of \cite{SmiA}.
\newhead{The Behaviour at Infinity}
\newsubhead{The Key Result}
\noindent The following theorem provides the key to the rest of this paper:
\proclaim{Theorem \procref{PresentationChIIILimites}\ {\sl Boundary Behaviour Theorem}}
\noindent let $S$ be a hyperbolic Riemann surface and let $\varphi:S\rightarrow\hat{\Bbb{C}}$ be a locally conformal mapping. For $k\in (0,1)$, let
$i:S\rightarrow U\Bbb{H}^3$ be an immersion such that $(S,i)$ is the unique solution to the Plateau problem $(S,\varphi)$
with constant Gaussian curvature $k$.
\medskip
\noindent Let $K$ be a compact subset of $S$ and let $\Omega$ be a connected component of $S\setminus K$. Let
$q$ be an arbitrary point in the boundary of $\varphi(\Omega)$ that is not in $\varphi(\overline{\Omega}\minter K)$.
If $\suite{p}{n}\in\Omega$ is a sequence of points such that $(\varphi(p_n))_{n\in\Bbb{N}}$ tends to $q$, then the sequence $(i(p_n))_{n\in\Bbb{N}}$ also tends towards $q$.
\endproclaim
\proof We will assume the contrary in order to obtain a contradiction. Let us denote by $\hat{\mathi}$ the Gauss lifting
of $i$. We identify $\Bbb{H}^3$ with $\Bbb{C}\times (0,\infty)$ and $\partial_\infty\Bbb{H}^3$ with $\hat{\Bbb{C}}$.
After applying an isometry of $\Bbb{H}^3$ if necessary, we may identify $q$ with $\infty$. For all $n$, we define
$q_n$ by $q_n=\varphi(p_n)$. Let $(z_n,\lambda_n)_{n\in\Bbb{N}}$ be a sequence such that, for all $n$:
$$
i(p_n) = (z_n,\lambda_n).
$$
\noindent Since $(i(p_n))_\ninn$ does not tend towards infinity, we may assume that there exists $R>0$ such that, for all
$n$:
$$
\|z_n,\lambda_n\|<R.
$$
\noindent Let us define the sequence of isometries $\suite{A}{n}$ of $\Bbb{H}^3$ such that, for all $n$:
$$
A_n(z,\lambda) = \frac{1}{\lambda_n}(z-z_n,\lambda).
$$
\noindent In particular, for all $n$, we have:
$$
A_n(z_n,\lambda_n) = (1,0).
$$
\noindent For all $n$, we note also by $A_n$ the automorphism of $\partial_\infty\Bbb{H}^3=\hat{\Bbb{C}}$ induced by the
action of $A_n$, and, if we denote by $\|\cdot\|$ the Euclidean norm over $\Bbb{C}$, we obtain:
$$
\|A_n(q_n)\| \geqslant \frac{1}{R}(\|q_n\| - R).
$$
\noindent In particular, since $\suite{q}{n}$ tends to infinity, the sequence of points $(A_n(q_n))_{n\in\Bbb{N}}$ also tends to infinity. For all $n$, we define
$i_n:S\rightarrow\Bbb{H}^3$ by:
$$
i_n = A_n\circ i.
$$
\noindent For every $n$, we denote the Gauss lifting of $i_n$ by $\hat{\mathi}_n$. Since $i_n(p_n)=(0,1)$, by Labourie's compactness theorem (theorem
\procref{ThmLabCompacite}), there exists an immersed surface $(S_0,\hat{\mathi}_0,p_0)$ in $U\Bbb{H}^3$ (which is possibly
a tube) such that $(S,\hat{\mathi}_n,p_n)_{n\in\Bbb{N}}$ tends towards this surface. Let $\overrightarrow{n}$ be the Gauss-Minkowksi mapping which sends $U\Bbb{H}^3$ to $\partial_\infty\Bbb{H}^3=\hat{\Bbb{C}}$. We obtain:
$$\matrix
(\overrightarrow{n}\circ\hat{\mathi}_0)(p_0) \hfill&=\mlim_{n\rightarrow\infty}(\overrightarrow{n}\circ\hat{\mathi}_n)(p_n) \hfill\cr
&=\mlim_{n\rightarrow\infty}(\overrightarrow{n}\circ A_n\circ\hat{\mathi})(p_n)\hfill\cr
&=\mlim_{n\rightarrow\infty}(A_n\circ\overrightarrow{n}\circ\hat{\mathi})(p_n)\hfill\cr
&=\mlim_{n\rightarrow\infty}(A_n\circ\varphi)(p_n) \hfill\cr
&=\mlim_{n\rightarrow\infty}A_n(q_n) \hfill\cr
&=\infty. \hfill\cr
\endmatrix$$
\noindent For $p\in S$ an arbitary point, $\epsilon\in (0,\infty)$ a positive real number, and $g$ a metric over $S$,
we define $B_\epsilon(p;g)$ to be the ball of radius $\epsilon$ about $p$ in $S$ with respect to the metric $g$. Let us
furnish $S$ with the metric $\hat{\mathi}^*g^\nu$. For all $\epsilon\in (0,\infty)$, since the surface
$(S,\hat{\mathi}^*g^\nu)$ is complete, there exists $N\in\Bbb{N}$ such that, for all $n\geqslant N$:
$$
B_\epsilon(p_n;\hat{\mathi}^*g^\nu) \subseteq \Omega.
$$
\noindent Indeed, otherwise, since these balls are connected, we may assume that, for all $n$:
$$
B_\epsilon(p_n;\hat{\mathi}^*g^\nu)\minter K \neq \emptyset.
$$
\noindent It thus follows that the sequence $\suite{p}{n}$ is contained in the ball of radius $\epsilon$ about the compact
set $K$. Since this ball is also compact, we may assume that there exists $p_0'\in\overline{\Omega}$ such that
the sequence $\suite{p}{n}$ converges to $p_0'$. By continuity $\varphi(p_0')=q$ and consequently $q$ is either
in the image of $\Omega$ or in the image of $\overline{\Omega}\minter K$. In either case, this contradicts the
hypotheses on $q$.
\medskip
\noindent For all $n$, let us define the metric $g_n$ over $S$ by:
$$
g_n = \hat{\mathi}_n^*g^\nu = \hat{\mathi}^* A_n^*g^\nu.
$$
\noindent For all $n$, since $A_n$ is an isometry, the metric $g_n$ coincides with $g$. For all $n$, let us define $B_n$
by:
$$
B_n = B_\epsilon(p_n,g_n).
$$
\noindent We may thus assume that $B_n$ is contained in $\Omega$ for all $n$, and consequently that $\infty$ is not in $\varphi(B_n)$. Since $A_n$ preserves $\infty$, we obtain:
$$
\infty \notin (A_n\circ\varphi)(B_n) = (\overrightarrow{n}\circ\hat{\mathi}_n)(B_n).
$$
\noindent By choosing $\epsilon$ to be sufficiently small, we may assume that the restriction of
$\overrightarrow{n}\circ\hat{\mathi}_0$ to $B_0=B_\epsilon(p_0,g_0)$ is a homeomorphism onto its image. Consequently, by common sense
lemma \procref{LemmeInjectiviteLorsquOnSApprocheDuLimite}, we may assume that, for all $n$, the restriction of the mapping
$(\overrightarrow{n}\circ\hat{\mathi}_n)$ to $B_n$ is a homeomorphism onto its image. Thus, by common sense lemma
\procref{LemmeCompacteDansLImageDeLaLimite}, $\infty$ is not in $(\overrightarrow{n}\circ\hat{\mathi}_0)(B_0)$. We
thus obtain the desired contradiction and the result follows.\qed
\newhead{The Geometry of The Problem $(\Bbb{D}^*,z\mapsto z)$}
\newsubhead{Overview of Geometric Properties of the Solution}
\noindent The solution to the problem $(\Bbb{D},z\mapsto z)$ will serve as a model for the study of the general case. In this section we will establish some of its geometric
properties.
\medskip
\noindent We identify $\Bbb{H}^3$ with $\Bbb{C}\times (0,\infty)$ and $\partial_\infty\Bbb{H}^3$ with $\hat{\Bbb{C}}$.
We define $\Bbb{D}^*$ by:
$$
\Bbb{D}^* = \left\{z\in\Bbb{C} | 0<\left|z\right|<1\right\}.
$$
\noindent We define $\varphi:\Bbb{D}^*\rightarrow\hat{\Bbb{C}}$ by:
$$
\varphi(z)=z.
$$
\noindent Since $\Bbb{D}^*$ is hyperbolic, by the hyperbolic existence theorem (theorem \procref{ThmSmithExistence}), there exists a unique solution
to the Plateau problem $(\Bbb{D}^*,\varphi)$. Let us denote this solution by $i:\Bbb{D}^*\rightarrow\Bbb{H}^3$, and let
$\hat{\mathi}$ be its Gauss lifting. We define the immersed surface $\Sigma$ by $\Sigma=(\Bbb{D}^*,i)$.
\medskip
\noindent The following result gives us a better idea of the shape of $\Sigma$:
\proclaim{Lemma \nextprocno\ {\sl First Structure Lemma.}}
\noindent There exists $f:\Bbb{D}^*\rightarrow (0,\infty)$ which only depends on $r=\left|z\right|$ such that $\Sigma$
coincides with the graph of $f$ over $\Bbb{D}^*$. Moreover $f(r)$ tends towards $0$ as $r$ tends towards $0$ and $1$.
\endproclaim
\proclabel{LemmeSolutionEstgraphe}
\noindent This result will be proven in section \subheadref{SubHeadGraph}. We denote by $h$ the metric on $\Bbb{H}^3$ and
we define the metric $g$ over $\Sigma$ by $g=i^*h$. By the uniqueness of solutions to the Plateau problem, $g$ is invariant
under rotations and reflections of $\Bbb{D}^*$. In section \subheadref{SubHeadBehaviourOfMetric}, we will prove the
following result concerning $g$:
\proclaim{Lemma \nextprocno\ {\sl Second Structure Lemma.}}
\noindent The Riemannian manifold $(\Bbb{D}^*,g)$ is complete. Moreover, if we define the vector field $\partial_\theta$
over $\Bbb{D}^*$ by:
$$
\partial_\theta(re^{i\theta}) = [t\rightarrow re^{i(\theta + t)}]_{t=0},
$$
\noindent then $g(\partial_\theta,\partial_\theta)$ tends towards $0$ as $r$ tends to $0$.
\endproclaim
\proclabel{LemmeMetriqueEstMince}
\remark For $r\in ]0,1[$ we define the curve $c_r:]0,2\pi[\rightarrow\Bbb{D}^*$ by:
$$
c_r(\theta) = re^{i\theta}.
$$
\noindent We define $\opLen(c_r,g)$ to be the length of $c_r$ with respect to the metric $g$. Since $g$ is symmetric, we obtain:
$$\matrix
\opLen(c_r,g) \hfill&=\int_0^{2\pi}\sqrt{g(\partial_\theta,\partial_\theta)}d\theta \hfill\cr
&=2\pi\sqrt{g(\partial_\theta,\partial_\theta)}. \hfill\cr
\endmatrix$$
\noindent It follows that $g(\partial_\theta,\partial_\theta)$ tends to $0$ as $r$ tends to $0$ if and only if
$\opLen(c_r,g)$ tends to $0$ as $r$ tends to $0$.
\medskip
\noindent We define $T=i(S)$ to be the image in $\Bbb{H}^3$ of the mapping $i$. Let $\overline{T}$ be the closure of $T$
in $\Bbb{H}^3\munion\partial_\infty\Bbb{H}^3$. We will obtain the following result:
\proclaim{Lemma \nextprocno\ {\sl Third Structure Lemma.}}
\noindent Let $K$ be a compact subset of $\Bbb{H}^3$. Let $\suite{p}{n}\in\Bbb{H}^3$ be a sequence which converges towards $0\in\partial_\infty\Bbb{H}^3\cong\hat{\Bbb{C}}$. Let
$\suite{A}{n}$ be a sequence of isometries of $\Bbb{H}^3$ such that, for all $n$:
$$
A_np_n\in K.
$$
\noindent If, for every $n$, we define $\overline{T}_n$ by $\overline{T}_n=A_n\overline{T}$, then there exists $\overline{T}_0\subseteq\Bbb{H}^3\munion\partial_\infty\Bbb{H}^3$
which is either the closure in $\Bbb{H}^3\munion\partial_\infty\Bbb{H}^3$ of a geodesic in $\Bbb{H}^3$, or a point in $\partial\Bbb{H}^3$, such that, after extraction of a subsequence, $\suite{\overline{T}}{n}$ converges to $\overline{T}_0$ in the Haussdorf topology.
\medskip
\noindent We define $\overline{\Gamma}_{0,\infty}$ to be the closure in $\Bbb{H}^3\munion\partial_\infty\Bbb{H}^3$ of the geodesic joining $0$ and $\infty$. If, for every $n$,
we define $\overline{\Gamma}_n$ by $\overline{\Gamma}_n=A_n\overline{\Gamma}_{0,\infty}$, then the sequence $\suite{\overline{\Gamma}}{n}$ converges to $\overline{T}_0$ in the
Haussdorf topology.
\endproclaim
\proclabel{LemmeConvergenceDeT}
\noindent This result will be proven in section \subheadref{SubHeadHaussdorfConvergence}.
\newsubhead{A Graph Over $\Bbb{D}^*$}
\noindent By symmetry with respect to reflections and rotations, and by the uniqueness of solutions to the Plateau problem, there exist functions $i_1:(0,1)\rightarrow\Bbb{R}$
and $i_2:(0,1)\rightarrow (0,\infty)$ such that:
\subheadlabel{SubHeadGraph}%
$$
i(re^{i\theta}) = (i_1(r)e^{i\theta}, i_2(r)).
$$
\noindent Using the boundary behaviour theorem (theorem \procref{PresentationChIIILimites}), we obtain the following result:
\proclaim{Lemma \nextprocno}
\noindent Let $S$ be a surface and let $\varphi:S\rightarrow\hat{\Bbb{C}}$ be a local diffeomorphism. Let $i:S\rightarrow\Bbb{H}^3$ be an immersion such that $(S,i)$ is the
unique solution to the Plateau problem $(S,\varphi)$.
\medskip
\noindent Let $q$ be a point in the boundary of $\varphi(S)$. If $\suite{p}{n}$ is a sequence of points in $S$ such that $(\varphi(p_n))_{n\in\Bbb{N}}$ converges to $q$, then
$(i(p_n))_{n\in\Bbb{N}}$ also converges to $q$.
\endproclaim
\proclabel{Limites}
\proof We take $K=\emptyset$ in theorem \procref{PresentationChIIILimites} and the result follows.\qed
\medskip
\noindent This permits us to establish certain properties of $i_1$ and $i_2$:
\noskipproclaim{Lemma \nextprocno}
$$\matrix
i_1(r)\rightarrow 0\hfill&\text{ as }r\rightarrow 0,\hfill\cr
i_1(r)\rightarrow 1\hfill&\text{ as }r\rightarrow 1,\hfill\cr
i_2(r)\rightarrow 0\hfill&\text{ as }r\rightarrow 0,1\hfill.\cr
\endmatrix$$
\endnoskipproclaim
\proclabel{LemmeBehaviourOfIAtLimites}
\proof By the preceeding lemma, $i(z)\rightarrow(0,0)$ as $z$ tends to $0$ and, for all $\theta\in [0,2\pi]$,
$i(z)$ converges to $(e^{i\theta},0)$ as $z$ tends to $e^{i\theta}$. The result now follows.\qed
\medskip
\noindent Next we have the following result:
\proclaim{Lemma \nextprocno}
\noindent The function $i_1$ is strictly increasing.
\endproclaim
\proof Since $i$ is smooth, the function $i_1$ is also. We recall that the intersection of a strictly convex
surface with a geodesic consists of isolated points. Moreover, in the canonical identification of $\Bbb{H}^3$ with
$\Bbb{C}\times (0,\infty)$, the vertical lines are geodesics. Thus, since $\Sigma$ is a strictly convex
surface of revolution, the critical points of $i_1$ are either strict local maxima or strict local minima (see figure
\figref{FigureConvexiteStricteDeSigma}), and, in particular, they are isolated. We denote the Gauss-Minkowksi application that sends $U\Bbb{H}^3$ to
$\partial_\infty\Bbb{H}^3$ by $\overrightarrow{n}$. We identify $\partial_\infty\Bbb{H}^3$ with $\hat{\Bbb{C}}$ and,
for $t\in (0,1)$, we find that $\overrightarrow{n}\circ\hat{\mathi}(t)\in\hat{\Bbb{R}}\subseteq\hat{\Bbb{C}}$. Since $\Sigma$
is strictly convex, figure \figref{FigureConvexiteStricteDeSigma} illustrates how, if $t\in (0,1)$ is a local
minimum of $i_1$, then:
$$
\overrightarrow{n}\circ\hat{\mathi}(t) < i_1(t).
$$
\placefigure{\null}{%
\placelabel[-0.3][0]{$\overrightarrow{n}\circ\hat{\mathi}(t)$}%
\placelabel[1.3][0]{$i_1(t)$}%
\placelabel[3.9][0]{$t$}%
\placelabel[3.4][2.6]{$i$}%
}{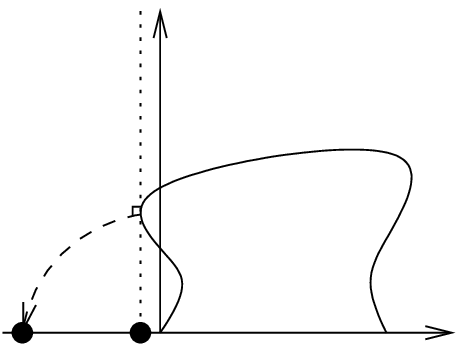}
\figlabel{FigureConvexiteStricteDeSigma}
\noindent Likewise, if $t\in (0,1)$ is a local maximum of $i_1$, then:
$$
\overrightarrow{n}\circ\hat{\mathi}(t) > i_1(t).
$$
\noindent Since $\overrightarrow{n}\circ\hat{\mathi}(t)=t\in ]0,1[$ and since $i_1(t)$ tends to $0$ and $1$ as $t$ tends to
$0$ and $1$ respectively, the function $i_1$ takes values in the interval $[0,1]$. Indeed, otherwise, by compactness, there
exists $t_0\in (0,1)$ such that, either $i_1(t_0)<0$ and $t_0$ is a minimum of $i_1$ or $i_1(t_0)>1$ and $t_0$ is a maximum
of $i_1$. In the first instance, we obtain:
$$
t_0 = (\overrightarrow{n}\circ\hat{\mathi})(t_0) < i_1(t_0) < 0.
$$
\noindent This is absurd. Likewise, the second possibility is absurd, and we thus obtain the desired
contradiction.
\medskip
\noindent Suppose that $t_0\in (0,1)$ is a strict local maximum of $i_1$. Since $i_1(t_0)\leqslant 1$, and since
$i_1(t)$ tends to $1$ as $t$ tends to $1$, there exists a strict local minimum $t_1$ of $i_1$ in the open interval
$(t_0,1)$ such that $i_1(t_1)<i_1(t_0)$. However:
$$\matrix
t_1 = \overrightarrow{n}\circ\hat{\mathi}(t_1) < i_1(t_1) \hfill\cr
\qquad < i_1(t_0) < \overrightarrow{n}\circ\hat{\mathi}(t_0) = t_0. \hfill\cr
\endmatrix$$
\noindent This is absurd. Consequently, there are no strict local maxima of $i_1$ in $(0,1)$. For the same reasons, there
are no strict local minima of $i_1$ in $(0,1)$. Consequently $i_1$ does not have any critical points in $(0,1)$ and
the result follows.\qed
\medskip
\noindent We may now prove the first structure lemma:
\medskip
\proclaim{Lemma \procref{LemmeSolutionEstgraphe} {\sl First Structure Lemma.}}
\noindent There exists $f:\Bbb{D}^*\rightarrow (0,\infty)$ which only depends on $r=\left|z\right|$ such that $\Sigma$
coincides with the graph of $f$ over $\Bbb{D}^*$. Moreover $f(r)$ tends towards $0$ as $r$ tends towards $0$ and $1$.
\endproclaim
\proof By the preceeding lemma, the application $i_1$ is strictly increasing and thus invertible. We define $\alpha:\Bbb{D}^*\rightarrow\Bbb{D}^*$ by:
$$
\alpha(r e^{i\theta}) = i_1^{-1}(r) e^{i\theta}.
$$
\noindent We define $f:\Bbb{D}^*\rightarrow (0,\infty)$ by:
$$
f(re^{i\theta})=i_2(i_1^{-1}(r)).
$$
\noindent The mapping $\alpha$ is a homeomorphism of $\Bbb{D}^*$ and:
$$\matrix
(i\circ\alpha)(re^{i\theta}) \hfill&= i(i_1^{-1}(r)e^{i\theta}) \hfill\cr
&= (i_1(i_1^{-1}(r))e^{i\theta}, i_2(i_1^{-1}(r))) \hfill\cr
&= (re^{i\theta}, f(re^{i\theta})). \hfill\cr
\endmatrix$$
\noindent The surface $\Sigma$ thus coincides with the graph of $f$ above $\Bbb{D}^*$. By definition, the function $f$
is independant of $\theta$, and by lemma \procref{LemmeBehaviourOfIAtLimites}, $f(r)$ tends to $0$ as $r$ tends to $0$ and $1$.
The result now follows.\qed
\newsubhead{The Properties of $i^*g$}
\noindent For all $\theta\in[0,2\pi]$, we define $D_\theta$ by:
\subheadlabel{SubHeadBehaviourOfMetric}%
$$
D_\theta = \left\{ z\in\Bbb{C} | \left|z-\frac{e^{i\theta}}{2}\right| < \frac{1}{2} \right\}.
$$
\noindent For all $\theta$, let $i_\theta$ be the unique solution to the Plateau problem
$(D_\theta, z\mapsto z)$ with constant Gaussian curvature equal to $k$. For all $\theta$, we define $\Sigma_\theta$ by
$\Sigma_\theta = (D_\theta, i_\theta)$. The surface $\Sigma_\theta$ is a surface equidistant from the (Euclidian)
hemisphere of radius $1/2$ centred on $e^{i\theta}/2$. In fact, $\Sigma_\theta$ is the intersection with the
upper half space of a Euclidean sphere whose centre is in the lower half space. Let $\Omega_\theta$ be the region exterior to this surface
(i.e. $\Omega_\theta$ is the intersection of the interior of this sphere with the upper half space).
We define $\Omega$ by:
$$
\Omega=\munion_{\theta\in[0,2\pi]}\Omega_\theta.
$$
\noindent We obtain the following result:
\proclaim{Lemma \nextprocno}
\noindent The surface $\Sigma$ est contained in the complement of $\Omega$.
\endproclaim
\proof Let $\theta\in[0,2\pi]$ be arbitrary. For $t\in ]0,1[$ we define $D_{t,\theta}$ by\colon
$$
D_{t,\theta} = \left\{z\in\Bbb{C} | \left|w-\frac{e^{i\theta}}{2}\right|<\frac{t}{2} \right\}.
$$
\noindent For all $t$, let $i_{t,\theta}$ be the unique solution to the Plateau problem $(D_{t,\theta},z\mapsto z)$
with constant Gaussian curvature equal to $k_t = (1-t) + tk > k$. For all $t$, we define $\Sigma_{t,\theta}$ by
$\Sigma_{t,\theta}=(D_{t,\theta}, i_{t,\theta})$. By considering the foliation of $\Bbb{H}^3$ defined by solutions of constant Gaussian
curvature $k_t$ to the Plateau problems given by all the discs centred on $e^{i\theta}/2$, using the weak geometric maximum principal
(lemma \procref{LemmaGeometricMaximumPrincipal}), we may show that the familly $(\Sigma_{t,\theta})_{t\in (0,1)}$ is a foliation of
$\Omega_\theta$.
\medskip
\noindent Let us denote by $\overline{\Sigma}_{t,\theta}$ the closure of the image of $\Sigma_{t,\theta}$ in
$\Bbb{H}^3\munion\partial_\infty\Bbb{H}^3$. Since $\partial_\infty\Sigma=\partial\Bbb{D}^*$, there exists
$\epsilon\in (0,1)$ such that, for all $t<\epsilon$:
$$
\overline{\Sigma}_{t,\theta}\minter\overline{\Sigma} = \emptyset.
$$
\noindent Let us define $t_0$ by:
$$
t_0 = \msup\left\{t | \overline{\Sigma}_{s,\theta}\minter\overline{\Sigma} = \emptyset\quad\forall s\in ]0,t[\right\}.
$$
\noindent We aim to show that $t_0=1$. We will assume the contrary in order to obtain a contradiction. Since $\overline{\Sigma}_\theta$ and $\overline{\Sigma}_{t_0,\theta}$
are compact, and since the foliation is continuous, we obtain:
$$
\overline{\Sigma}_{t_0,\theta}\minter\overline{\Sigma}\neq\emptyset.
$$
\noindent Now:
$$\matrix
&\partial_\infty\Sigma_{t_0,\theta} \hfill&=\partial D_{t_0,\theta} \hfill&\subseteq\hat{\Bbb{C}},\hfill\cr
&\partial_\infty\Sigma \hfill&=\partial \Bbb{D}^* \hfill&\subseteq\hat{\Bbb{C}}.\hfill\cr
\endmatrix$$
\noindent Thus:
$$\matrix
& \partial_\infty\Sigma_{t_0,\theta}\minter\partial_\infty\Sigma \hfill&=\emptyset \hfill& \cr
\Rightarrow \hfill& \Sigma_{t_0,\theta}\minter\Sigma \hfill&\neq\emptyset.\hfill& \cr
\endmatrix$$
\noindent Let $p_0$ be in the intersection of $\Sigma_{t_0,\theta}$ and $\Sigma$, and let us denote by $\opExt(\Sigma_{t_0,\theta})$ the exterior of $\Sigma_{t_0,\theta}$. We
obtain:
$$\matrix
&\opExt(\Sigma_{t_0,\theta})\hfill&= \munion_{0<t<t_0}\Sigma_{t,\theta} \hfill\cr
\Rightarrow &\Sigma_R\minter \opExt(\Sigma_{t_0,\theta}) \hfill&= \emptyset. \hfill\cr
\endmatrix$$
\noindent It follows that $\Sigma$ is tangent to $\Sigma_{t_0,\theta}$ within the interior of this surface at $p_0$.
However, since the Gaussian curvature of $\Sigma_{t_0,\theta}$ is strictly greater than that of $\Sigma$, we obtain
a contradiction by the weak geometric maximum principal (lemma \procref{LemmaGeometricMaximumPrincipal}). It thus follows that $t_0=1$ and we obtain:
$$
\Sigma\subseteq \Omega^c_\theta.
$$
\noindent Since $\theta\in [0,2\pi]$ is arbitrary, the result follows.\qed
\medskip
\noindent We now obtain the following result concerning the behaviour of $i$ and $f$:
\proclaim{Corollary \nextprocno}
\noindent There exists $B\in ]0,\infty[$ such that:
$$
\mlimsup_{r\rightarrow 0}\frac{i_1(r)}{i_2(r)} \leqslant B.
$$
\noindent In other words:
$$
\mlimsup_{r\rightarrow 0}\frac{r}{f(r)} \leqslant B.
$$
\endproclaim
\proof Since $f(r)=(i_2\circ i_1^{-1})(r)$ and $i_1(0)=0$, these two results are equivalent. Let $\Sigma_{1,0}$ be as in the proof of the
preceeding lemma. Let
$\tilde{f}:D_0\rightarrow ]0,\infty[$ be such that $\Sigma_{1,0}$ is the graph of $\tilde{f}$ over $D_0$. By the preceeding lemma:
$$
f \geqslant \tilde{f}.
$$
\noindent However, by considering the restriction of $\tilde{f}$ to $(0,1)$, there exists $B\in (0,\infty)$ such that:
$$
\mlimsup_{r\rightarrow 0}\frac{r}{\tilde{f}(r)} \leqslant B.
$$
\noindent The result now follows.\qed
\medskip
\noindent Let $\Gamma\subseteq\Bbb{H}^3$ be the geodesic going from $0$ to $\infty$. There exists a function $\delta:(0,\infty)\rightarrow (0,\infty)$ such that for all
$r$,$\theta$ and $\lambda$:
$$
d((re^{i\theta},\lambda),\Gamma) = \delta(r/\lambda).
$$
\noindent By the preceeding result, there exists $B\in (0,\infty)$ such that:
$$
\mlimsup_{z\rightarrow 0}\delta(i(z)) \leqslant B.
$$
\noindent We may now refine this estimate as follows:
\noskipproclaim{Lemma \nextprocno}
$$
\mlim_{z\rightarrow 0}\delta(i(z)) = 0.
$$
\endnoskipproclaim
\proof Let $\suite{p}{n}\in\Bbb{D}^*$ be a sequence of points that converges to $0$ such that:
$$
\delta(i(p_n))\rightarrow\mlimsup_{z\rightarrow 0}\delta(i(z)).
$$
\noindent For all $n$, let us define $(z_n,\lambda_n)_{n\in\Bbb{N}}\in\Bbb{H}^3$ by:
$$
(z_n,\lambda_n) = i(p_n).
$$
\noindent For all $n$, we define the isometry $A_n:\Bbb{H}^3\rightarrow\Bbb{H}^3$ by:
$$
A_n(z,\lambda) = \frac{1}{\lambda_n}(z,\lambda).
$$
\noindent In particular:
$$\matrix
(A_n\circ i)(p_n) \hfill&= A_n(z_n,\lambda_n) \hfill\cr
&=(\frac{z_n}{\lambda_n}, 1). \hfill\cr
\endmatrix$$
\noindent Since $\mlimsup_{n\rightarrow\infty}\left|z_n/\lambda_n\right|\leqslant B$, there exists a compact subset $K$
of $\Bbb{H}^3$ such that, for all $n$:
$$
(A_n\circ i)(p_n) \in K.
$$
\noindent For all $n$, we define $i_n$ by $i_n = A_n\circ i$ and we denote the Gauss lifting of $i_n$ by $\hat{\mathi}_n$.
For all $n$, we define the immersed surface $\Sigma_n$ by $\Sigma_n=(\Bbb{D}^*,i_n)$ and we denote the Gauss lifting of
$\Sigma_n$ by $\hat{\Sigma}_n$. By Labourie's compactness theorem (theorem \procref{ThmLabCompacite}), there exists a (possibly tubular) pointed immersed
surface $(\hat{\Sigma}_0,p_0)=(S_0,\hat{\mathi}_0,p_0)$ in $U\Bbb{H}^3$ such that $(\hat{\Sigma}_n,p_n)_{n\in\Bbb{N}}$ converges to
$(\hat{\Sigma}_0,p_0)$ in the Cheeger/Gromov topology.
\medskip
\noindent We define $\daleth$ by:
$$
\daleth=\limsup_{z\rightarrow 0}\delta(i(z)).
$$
\noindent Let $\eta\in(0,\infty)$ be an arbitrary positive real number. By definition, there exists a positive real number
$\epsilon\in(0,\infty)$ such that, for $0<\left|z\right|<\epsilon$:
$$
\delta(i(z)) \leqslant \daleth + \eta.
$$
\noindent For $g$ an arbitrary metric over $\Bbb{D}^*$, for $p\in\Bbb{D}^*$ an arbitrary point and for $R\in (0,\infty)$
an arbitrary positive real number, we define $d_g$ to be the metric (i.e. the distance function) generated over $\Bbb{D}^*$
by $g$ and we define $B_R(p,g)$ by:
$$
B_R(p,g) = \left\{q\in\Bbb{D}^* | d_g(p,q) < R \right\}.
$$
\noindent For all $n$, since $A_n$ is an isometry, we obtain:
$$\matrix
& \hat{\mathi}_n^*g^\nu \hfill&=\hat{\mathi}^*g^\nu \hfill\cr
\Rightarrow \hfill& B_R(p_n,\hat{\mathi}_n^*g^\nu) \hfill&= B_R(p_n,\hat{\mathi}^*g^\nu). \hfill\cr
\endmatrix$$
\noindent We fix $R$. Since $(\Bbb{D}^*, \hat{\mathi}^*g^\nu)$ is complete, if we define $D_\epsilon^*$ by:
$$
D_\epsilon^* = \left\{z\in\Bbb{D}^* | \left|z\right| < \epsilon \right\},
$$
\noindent then there exists a positive integer $N\in\Bbb{N}$ such that for all $n\geqslant N$:
$$
B_R(p_n,\hat{\mathi}^*g^\nu) \subseteq D_\epsilon^*.
$$
\noindent Consequently, for $n\geqslant N$:
$$\matrix
& B_R(p_n, \hat{\mathi}_n^*g^\nu) \hfill&\subseteq D_\epsilon^* \hfill\cr
\Rightarrow \hfill& \msup\left\{\delta(i_n(q)) | q\in B_R(p_n,\hat{\mathi}_n^*g^\nu) \right\} \hfill&\leqslant\daleth + \eta.\hfill\cr
\endmatrix$$
\noindent Thus, after taking limits, we obtain:
$$
\msup\left\{\delta(i_0(q)) | q\in B_R(p_0, \hat{\mathi}_0^*g^\nu) \right\} \leqslant \daleth + \eta.
$$
\noindent Since $\eta, R\in (0,\infty)$ are both arbitrary, we have:
$$
\msup\left\{\delta(i_0(q)) | q\in S_0\right\} \leqslant \daleth.
$$
\noindent However, by definition:
$$
\delta(i_n(p_n))\rightarrow \daleth\text{ as }n\rightarrow\infty.
$$
\noindent Consequently:
$$
\delta(i_0(p_0)) = \daleth.
$$
\noindent We will show that $\hat{\Sigma}_0$ is a tube. Indeed, suppose the contrary, in which case it is the Gauss
lifting of a k-surface $\Sigma_0=(S,i_0)$. The surface $\Sigma_0$ is an interior tangent at the point $p_0$ to the
surface $\delta^{-1}(B)$. However, the Gaussian curvature of $\delta^{-1}(B)$ is equal to $1$ and is thus strictly greater than $k$. The desired contradiction now follows by
the weak geometric maximum principal (lemma \procref{LemmaGeometricMaximumPrincipal}). Consequently $\hat{\Sigma}_0$ must be a tube about a geodesic $\Delta$.
\medskip
\noindent For all $p\in\Delta$:
$$
\delta(p) \leqslant \daleth.
$$
\noindent Consequently, $\Delta$ remains within a fixed distance of $\Gamma$. It thus follows that $\Delta$ and $\Gamma$ coincide and that, for all $p\in\Delta$:
$$
\delta(p) = 0.
$$
\noindent In particular, we obtain:
$$
\daleth = \delta(i_0(p_0)) = 0.
$$
\noindent The desired result now follows.\qed
\medskip
\noindent In particular, we obtain:
\noskipproclaim{Corollary \nextprocno}
$$
\mlim_{r\rightarrow 0}\frac{i_1(r)}{i_2(r)} = \mlim_{r\rightarrow 0}\frac{r}{f(r)} = 0.
$$
\endnoskipproclaim
\proclabel{CorLaFormeDeLagraphe}
\noindent We may now prove the second structure lemma:
\proclaim{Lemma \procref{LemmeMetriqueEstMince} {\sl Second Structure Lemma.}}
\noindent The Riemannian manifold $(\Bbb{D}^*,g)$ is complete. Moreover, if we define the vector field $\partial_\theta$
over $\Bbb{D}^*$ by:
$$
\partial_\theta(re^{i\theta}) = [t\rightarrow re^{i(\theta + t)}]_{t=0},
$$
\noindent then $g(\partial_\theta,\partial_\theta)$ tends towards $0$ as $r$ tends to $0$.
\endproclaim
\proof Let $h$ be the Riemannian metric over $\Bbb{H}^3$. By the boundary behaviour theorem (theorem \procref{PresentationChIIILimites}), the
immersion $i$ is proper, and thus the metric $g=i^*h$ is complete. For $r\in (0,1)$, we define the curve
$c_r:(0,2\pi)\rightarrow\Bbb{D}^*$ such that, for all $\theta$:
$$
c_r(\theta)=re^{i\theta}.
$$
\noindent We have:
$$\matrix
& (i\circ c_r)(\theta) \hfill&=(i_1(r)e^{i\theta}, i_2(r)) \hfill\cr
\Rightarrow \hfill& \opLen_h(i\circ c_r) \hfill&= 2\pi\frac{i_1(r)}{i_2(r)} \hfill\cr
\Rightarrow \hfill& \lim_{r\rightarrow 0}\opLen_h(i\circ c_r) \hfill&= 0. \hfill\cr
\endmatrix$$
\noindent However, for all $r$:
$$
\opLen_h(i\circ c_r) = \opLen_g(c_r).
$$
\noindent The result now follows.\qed
\newsubhead{Convergence in the Haussdorf Topology}
\noindent We now prove the third structure lemma:
\subheadlabel{SubHeadHaussdorfConvergence}%
\proclaim{Lemma \procref{LemmeConvergenceDeT} {\sl Third Structure Lemma.}}
\noindent Let $K$ be a compact subset of $\Bbb{H}^3$. Let $\suite{p}{n}\in\Bbb{H}^3$ be a sequence which converges towards $(0,0)\in\partial_\infty\Bbb{H}^3$. Let
$\suite{A}{n}$ be a sequence of isometries of $\Bbb{H}^3$ such that, for all $n$:
$$
A_np_n\in K.
$$
\noindent If, for every $n$, we define $\overline{T}_n$ by $\overline{T}_n=A_n\overline{T}$, then there exists $\overline{T}_0\subseteq\Bbb{H}^3\munion\partial_\infty\Bbb{H}^3$
which is either the closure in $\Bbb{H}^3\munion\partial_\infty\Bbb{H}^3$ of a geodesic in $\Bbb{H}^3$, or a point in $\partial\Bbb{H}^3$, such that $\suite{\overline{T}}{n}$
converges to $\overline{T}_0$ in the Haussdorf topology.
\medskip
\noindent We define $\overline{\Gamma}_{0,\infty}$ to be the closure in $\Bbb{H}^3\munion\partial_\infty\Bbb{H}^3$ of the geodesic joining $0$ and $\infty$. If, for every $n$,
we define $\overline{\Gamma}_n$ by $\overline{\Gamma}_n=A_n\overline{\Gamma}_{0,\infty}$, then the sequence $\suite{\overline{\Gamma}}{n}$ converges to $\overline{T}_0$ in the
Haussdorf topology.
\endproclaim
\proof Let $\suite{p}{n}\in\Bbb{H}^3$ be a sequence in $\Bbb{H}^3$ that converges to $0$. For all $n$, let us define $(w_n,\lambda_n)$ by:
$$
p_n = (w_n,\lambda_n).
$$
\noindent For all $n$, we define the isometry $M_n$ of $\Bbb{H}^3$ by:
$$
M_n(z,\lambda) = \frac{1}{\lambda_n}(z,\lambda).
$$
\noindent For all $R\in (0,\infty)$, we define $B_R\subseteq\Bbb{H}^3\munion\partial_\infty\Bbb{H}^3$ by:
$$
B_R=\left\{(z,\lambda)\in\Bbb{C}\times[0,\infty[\text{ s.t. }\left|z\right|^2+\lambda^2\geqslant R^2\right\}\munion\left\{\infty\right\}.
$$
\noindent For all $r\in (0,\infty)$ we define $C_r\subseteq\Bbb{H}^3\munion\partial_\infty\Bbb{H}^3$ by:
$$
C_r=\left\{(z,\lambda)\in\Bbb{C}\times[0,\infty[\text{ s.t. }\left|z\right|^2\leqslant r^2\lambda^2\right\}\munion\left\{\infty\right\}.
$$
\noindent We now define the mushroom $\opMush_{R,r}\subseteq\Bbb{H}^3\munion\partial_\infty\Bbb{H}^3$ by:
$$
\opMush_{R,r} = B_r\munion C_r.
$$
\noindent The mushrooms $\opMush_{R,r}$ converge to $\overline{\Gamma}_{0,\infty}$ in the Haussdorf topology as $R$ tends
to infinity and $r$ tends to $0$.
\medskip
\noindent Suppose that there exists $B\in (0,\infty)$ such that, for all $n$:
$$
\left|\frac{w_n}{\lambda_n}\right| < B.
$$
\noindent There exists a compact subset $L$ of $\Bbb{H}^3$ such that, for all $n$:
$$
M_np_n\in L.
$$
\noindent By the first and second structure lemmata (lemmata \procref{LemmeSolutionEstgraphe} and  \procref{LemmeMetriqueEstMince}), after taking a subsequence if necessary, we
may assume that there exists $\suite{R}{n},\suite{r}{n}\in (0,\infty)$ such that $(R_n/\lambda_n)_{n\in\Bbb{N}}$ tends to infinity and $\suite{r}{n}$ tends to $0$, and, for
all $n$:
$$
\overline{T} \subseteq \opCh_{R_n,r_n}.
$$
\noindent Consequently:
$$
M_n\overline{T} \subseteq \opCh_{R_n/\lambda_n,r_n}.
$$
\noindent It follows that $M_n\overline{T}$ converges towards $\overline{\Gamma}_{0,\infty}$ in the Haussdorf topology. For all $n$, we define the application $B_n$ by
$B_n=A_nM_n^{-1}$ and we obtain:
$$
B_n(M_n p_n)\in K.
$$
\noindent However, since $L$ and $K$ are both compact, the familly of isometries of $\Bbb{H}^3$ which send a point of
$L$ onto a point of $K$ is also compact. It follows that, after taking a further subsequence if necessary, we may assume
that there exists $B_0$ such that $\suite{B}{n}$ converges to $B_0$. By lemma \procref{LemmeConvergenceDesEnsemblesEtDesFonctions}, it follows that
$(A_n\overline{T})_{n\in\Bbb{N}}=(B_nM_n(\overline{T}))_{n\in\Bbb{N}}$ and $(A_n\overline{\Gamma}_{0,\infty})_{n\in\Bbb{N}}=(B_n\overline{\Gamma}_{0,\infty})_{n\in\Bbb{N}}$
both converge to $A_0\Gamma_{0,\infty}$ in the Haussdorf topology.
\medskip
\noindent We now suppose that no such $B$ exists. For all $n$, we define $\rho_n$ by $\rho_n=\left|w_n/\lambda_n\right|$.
After taking a subsequence if necessary, we may assume that $\suite{\rho}{n}$ tends to infinity. By the first and second structure lemmata (lemmata
\procref{LemmeSolutionEstgraphe} and  \procref{LemmeMetriqueEstMince}), after taking a further subsequence if necessary, we may assume that there
exists $\suite{R}{n},\suite{r}{n},\suite{K}{n}\in ]0,\infty[$ such that:
\medskip
\myitem{(i)} the sequences $(R_n/\lambda_n)_{n\in\Bbb{N}}$, $\suite{r}{n}$ and $\suite{K}{n}$ tend to infinity,
\medskip
\myitem{(ii)} for all $n$, $R_n/\lambda_n\geqslant K_n\rho_n$, and
\medskip
\myitem{(iii)} for all $n$:
$$\matrix
&\overline{T} \hfill&\subseteq \opCh_{R_n,r_n}\hfill\cr
\Rightarrow \hfill&M_n\overline{T} \hfill&\subseteq \opCh_{R_n/\lambda_n,r_n}.\hfill\cr
\endmatrix$$
\noindent For all $n$, we define the application $N_n$ by:
$$
N_n(z,\lambda) = (z-\frac{w_n}{\lambda_n},\lambda).
$$
\noindent For all sufficiently large $n$, we obtain:
$$
N_n\opCh_{R_m/\lambda_n,r_n} \subseteq B^c_{\rho_n/2}.
$$
\noindent Consequently, for sufficiently large $n$:
$$
N_nM_n\overline{T}_n, N_nM_n\overline{\Gamma}_{0,\infty} \subseteq B^c_{\rho_n/2}.
$$
\noindent Since $(B^c_{\rho_n/2})_{n\in\Bbb{N}}$ converges towards $\left\{\infty\right\}$ in the Haussdorf topology, it follows that\break
$(N_nM_n\overline{T}_n)_{n\in\Bbb{N}}$ and
$(N_nM_n\overline{\Gamma}_{0,\infty})_{n\in\Bbb{N}}$ also converge to $\left\{\infty\right\}$ in the Haussdorf topology. For all $n$, we define $B_n$ by $B_n=A_n(N_nM_n)^{-1}$.
By following the same reasoning as before, after taking a subsequence if necessary, we may assume that there exists $p_0\in\partial_\infty\Bbb{H}^3$ such that
$(A_n\overline{T}_n)_{n\in\Bbb{N}}$ and $(A_n\overline{\Gamma}_{0,\infty})_{n\in\Bbb{N}}$ converge towards $\left\{p_0\right\}$ in the Haussdorf topology, and the result
follows.\qed
\newhead{Ramified Coverings}
\newsubhead{Introduction}
\noindent In this section we will prove theorem \procref{PresentationChIIIRevetementsRamifiees}:
\proclaim{Theorem \procref{PresentationChIIIRevetementsRamifiees}}
\noindent Let $S$ be a Riemann surface. Let $\Cal{P}$ be a discrete subset of $S$ such that $S\setminus\Cal{P}$ is hyperbolic. Let $\varphi:S\rightarrow\hat{\Bbb{C}}$ be a ramified covering having critical points in $\Cal{P}$. Let $\kappa$ be a real number in $(0,1)$. Let
$i:S\setminus\Cal{P}\rightarrow\Bbb{H}^3$ be the unique solution to the Plateau problem $(S\setminus\Cal{P},\varphi)$ with constant Gaussian
curvature $\kappa$. Let $\hat{\Sigma}=(S\setminus\Cal{P},\hat{\mathi})$ be the Gauss lifting of $\Sigma$.
\medskip
\noindent Let $p_0$ be an arbitrary point in $\Cal{P}$. If $\varphi$ has a critical point of order $k$ at $p_0$, then $\hat{\Sigma}$ is asymptotically tubular of order
$k$ at $p_0$.
\endproclaim
\noindent We will proceed in many steps. We first prove that if $(p_n)_\ninn$ is a sequence of points in $S\setminus\Cal{P}$ which converges to $p_0$, then
$(\hat{\Sigma},p_n)_\ninn$ converges to a tube in the Cheeger/Gromov topology. We show that this tube is necessarily of order $k$. We then show how convergence in the
Cheeger/Gromov topology allows us to deduce that $(\Sigma,p_n)$ is a graph over a tube of given finite length for all sufficiently large $n$. Finally, by glueing together
such graphs, we will obtain the desired result.
\newsubhead{The Position of $\hat{\mathi}(p)$ Near to Ramification Points}
\noindent Let $S$ be a Riemann surface and let $\Cal{P}\subseteq S$ be a discrete
subset. Let $\varphi:S\rightarrow\hat{\Bbb{C}}$ be a ramified covering of $S$ over $\hat{\Bbb{C}}$ with critical points in $\Cal{P}$.
let $\kappa\in(0,1)$ be a real number.
Let $i:S\setminus\Cal{P}\rightarrow\Bbb{H}^3$ be the unique solution to the Plateau problem $(S\setminus\Cal{P},\varphi)$ with constant Gaussian curvature $\kappa$, and let us define the
immersed surface $\Sigma$ by $\Sigma=(S\setminus\Cal{P},i)$. Let $\hat{\mathi}$ be the Gauss lifting of $i$, and let $\hat{\Sigma}=(S\setminus\Cal{P},\hat{\mathi})$ be the Gauss
lifting of $\Sigma$. Let $p_0\in\Cal{P}$ be a ramification point of $\varphi$. Let $\Gamma_{0,\infty}$ be the unique geodesic in $\Bbb{H}^3$ going from $0$ to $\infty$.
\medskip
\noindent We begin by recalling the following result which gives a local description of ramified coverings near to ramification points:
\proclaim{Lemma \nextprocno}
\noindent Suppose that $\varphi(p_0)=0$. There exists a chart $(f,\Omega,\Bbb{D})$ of $S$ about $p_0$, a real number $\lambda\in (0,\infty)$ and a positive integer
$k\in\Bbb{N}$ such that the following diagram commutes:
\goodbreak
$$\commdiag{
\Omega & & \cr
\mapdown_f & \arrow(3,-2)\rt{\varphi} & \cr
\Bbb{D} & \mapright^{z\mapsto\lambda z^k} & \hat{\Bbb{C}} \cr
}$$
\endproclaim
\proclabel{LemmeTrivDesRevetementsRamifiees}
\noindent Using the boundary behaviour theorem (theorem \procref{PresentationChIIILimites}), we obtain the following result:
\proclaim{Lemma \nextprocno}
\noindent Let $p_0\in\Cal{P}$ be a ramification point of $\varphi$. If $\suite{p}{n}\in S\setminus\Cal{P}$ is a sequence which converges to $p_0$, then
$(i_n(p_n))_\ninn$ tends towards $\varphi(p_0)$ in $\Bbb{H}^3\munion\partial_\infty\Bbb{H}^3$.
\endproclaim
\proclabel{RevetementsRamifies}
\proof By the preceeding  lemma, after composing $\varphi$ with an isometry of $\Bbb{H}^3$ if necessary, we may assume that there exists
a chart $(f,U,2\Bbb{D})$ of $S\setminus\Cal{P}$ about $p_0$ and $k\in\Bbb{N}$ such that the following diagram commutes:
$$\commdiag{
U & & \cr
\mapdown_f & \arrow(3,-2)\rt{\varphi} & \cr
2\Bbb{D} & \mapright^{z\mapsto z^{k}} &\hat{\Bbb{C}} \cr
}$$
\noindent We define the compact subset $K$ of $S\setminus\Cal{P}$ by $K=f^{-1}(\left\{\left|z\right|=1\right\})$ and we define the connected component $\Omega$ of $(S\setminus\Cal{P})\setminus K$ by $\Omega=f^{-1}(\left\{0<\left|z\right|<1\right\})$. The result now follows by the boundary behaviour theorem (theorem
\procref{PresentationChIIILimites}).\qed
\newsubhead{Tubes About Geodesics}
\noindent We begin by controlling the geometry of $\Sigma$ near to $p_0$. For simplicity, we will assume that $\varphi(p_0)=0$. Let $\suite{p}{n}\in S\setminus\Cal{P}$ be a
sequence which tends to $p_0$. By identifying $\Bbb{H}^3$ with $\Bbb{C}\times (0,\infty)$, for all $n$, we define $(w_n,\lambda_n)_{n\in\Bbb{N}}\in\Bbb{H}^3$ by:
$$
(w_n,\lambda_n) = i(p_n).
$$
\noindent For all $n$, we define the isometry $A_n$ of $\Bbb{H}^3$ by:
$$
A_n(w,\lambda) = \frac{1}{\lambda_n}(w,\lambda).
$$
\noindent For all $n$, we define the immersion $i_n$ by:
$$
i_n=A_ni.
$$
\noindent For all $n$, let $\hat{\mathi}_n$ be the Gauss lifting of $i_n$. For all $n$, we define the immersed surface $\Sigma_n$ by $\Sigma_n=(S\setminus\Cal{P},i_n)$ and
we define $\hat{\Sigma}_n=(S\setminus\Cal{P},\hat{\mathi}_n)$ to be the Gauss lifting of $\Sigma$.
\medskip
\noindent We obtain the following result:
\proclaim{Lemma \nextprocno}
\noindent After extraction of a subsequence, $(i_n(p_n))_\ninn$ converges towards $(0,1)$ and $(\hat{\Sigma}_n,p_n)_\ninn$ converges to a tube about $\Gamma_{0,\infty}$.
\endproclaim
\proclabel{LemmeSurfaceTubulaireAutourDeGamma}
\proof By lemma \procref{LemmeTrivDesRevetementsRamifiees}, after applying an isometry of $\Bbb{H}^3$ if necessary, we may find a chart $(z,U,\Bbb{D})$ of $S$ about $p_0$
and $k\in\Bbb{N}$ such that the following diagram commutes:
$$\commdiag{
U & & \cr
\mapdown_z & \arrow(3,-2)\rt{\varphi} & \cr
\Bbb{D} & \mapright^{z\mapsto z^k} & \hat{\Bbb{C}} \cr
}$$
\noindent For all $l$, we define the function $\varphi_l:\Bbb{D}^*\rightarrow\hat{\Bbb{C}}$ by:
$$
\varphi_l(z)=z^l.
$$
\noindent For all $l$, let us define the immersion $f_l:\Bbb{D}^*\rightarrow\hat{\Bbb{C}}$ to be the unique solution of the Plateau problem $(\Bbb{D}^*,\varphi_l)$. For all
$l$, we define $\Sigma_l=(\Bbb{D}^*,f_l)$ and we define $\Cal{T}_l=f_l(\Bbb{D}^*)$ to be the image of $f_l$ in $\Bbb{H}^3$. Let us define $T$ by $T=\Cal{T}_1$. $\Cal{T}_l$
coincides with $T$ for all $l$. Let $\overline{T}$ be the closure of $T$ in $\Bbb{H}^3\munion\partial_\infty\Bbb{H}^3$. By Labourie's uniqueness theorem (theorem
\procref{LemmeLabNature}), the immersed surface $\Sigma_k$ is a graph above $U\setminus\left\{p_0\right\}$ in the extension of $\Sigma$. By using this to control the
behaviour of $\Sigma$ near $p_0$ in terms of $\Sigma_k$, we will be able to conclude.
\medskip
\noindent For all $n$, let us define the isometry $M_n$ of $\Bbb{H}^3$ by:
$$
M_n(w,\lambda) = \left(w - \frac{w_n}{\lambda_n},\lambda\right).
$$
\noindent In particular, for all $n$, $M_n(i_n(p_n))=(0,1)$. For all $n$ we define the immersion $j_n$ by $j_n=M_n\circ i_n$ and we denote the Gauss lifting of $j_n$ by
$\hat{\mathj}_n$. For all $n$, we define the immersed surface $\Sigma'_n$ by $\Sigma'_n=(S\setminus\Cal{P},j_n)$ and we define $\hat{\Sigma}'_n=(S\setminus\Cal{P},\hat{\mathj}_n)$
to be the Gauss lifting of $\Sigma'_n$. By Labourie's compactness theorem (theorem \procref{ThmLabCompacite}), after extracting a subsequence if necessary, we may assume that
there exists a (possibly tubular) pointed immersed surface $(\hat{\Sigma}'_0,p_0)=(S_0,\hat{\mathj}_0,p_0)$ in $U\Bbb{H}^3$ such that $(\hat{\Sigma}'_n,p_n)_\ninn$ converges towards
$(\hat{\Sigma}'_0,p_0)$ in the Cheeger/Gromov topology.
\medskip
\noindent For all $n$, we define $\overline{T}_n$ by $\overline{T}_n = M_nA_n\overline{T}$. By the third structure lemma (lemma \procref{LemmeConvergenceDeT}), after extracting a
subsequence if necessary, we may assume that there exists $\overline{T}_0\subseteq\Bbb{H}^3\munion\partial_\infty\Bbb{H}^3$ which is either the closure in
$\Bbb{H}^3\munion\partial_\infty\Bbb{H}^3$ of a geodesic in $\Bbb{H}^3$ or a point in $\partial_\infty\Bbb{H}^3$ such that $\suite{\overline{T}}{n}$ converges towards
$\overline{T}_0$ in the Haussdorf topology. At the same time $(M_n\overline{\Gamma}_{0,\infty})_\ninn$ converges to $\overline{T}_0$ in the Haussdorf topology.
\medskip
\noindent Let $q_0$ be an arbitrary point in $S_0$. Let $R\in (0,\infty)$ be a positive real number such that $q_0\in B_R(p_0)$. Let $\suite{q}{n}$ be a sequence of
points in $S\setminus\Cal{P}$ which tends to $q_0$ such that, for all $n$:
$$
q_n \in B_R(p_n) \subseteq (S\setminus\Cal{P},\hat{\mathj}_n^*g^\nu).
$$
\noindent Since, for every $n$, the surface $(S\setminus\Cal{P},\hat{\mathj}_n^*g^\nu)$ is complete, for all sufficently large $n$, we obtain:
$$
B_R(p_n) \subseteq U\setminus\left\{p_0\right\}.
$$
\noindent Consequently, for all sufficiently large $n$, the point $q_n$ is in $U\setminus\left\{p_0\right\}$. Let $\opExp:T\Bbb{H}^3\rightarrow\Bbb{H}^3$ be the exponential
mapping over $\Bbb{H}^3$. For all $n$, bearing in mind that $\hat{\mathj}_n(q_n)$ is a vector in $U\Bbb{H}^3$, we define the subset $X_n$ of
$\Bbb{H}^3$ by:
$$
X_n = \opExp(\hat{\mathj}_n(q_n)[0,\infty]) \subseteq \Bbb{H}^3\munion\partial_\infty\Bbb{H}^3.
$$
\noindent The sequence $\suite{X}{n}$ converges in the Haussdorf topology towards $\opExp(\hat{\mathj}_0(q_0)[0,\infty])$ as $n$ tends to infinity. Since $\Sigma_k$ is a graph
above $U\setminus\left\{p_0\right\}$ in the extension of $\Sigma$, it follows that for every $n$:
$$
\overline{X}_n\minter\overline{T}_n \neq \emptyset.
$$
\noindent Consequently, by common sense lemma \procref{LemmeLIntersectionDesLimites}:
$$
\opExp(\hat{\mathj}_0(q_0)[0,\infty])\minter \overline{T}_0\neq\emptyset
$$
\noindent We will begin by showing that $\hat{\Sigma}'_0$ is a tube. Let us assume the contrary in order to obtain a contradiction. There exists an immersion
$j_0:S_0\rightarrow\Bbb{H}^3$ such that $\hat{\mathj}_0$ is the Gauss lifting of $j_0$. The immersed surface $\Sigma'_0=(S_0,j_0)$ is everywhere locally convex. Let $V$ be
a sufficently small open subset of $S_0$ such that the immersed surface $(V,j_0)$ coincides with a portion of the boundary of a strictly convex subset of $\Bbb{H}^3$.
Let us define the application $\Cal{E}:V\times [0,\infty)\rightarrow\Bbb{H}^3$ and $\Cal{E}_\infty:V\rightarrow\partial_\infty\Bbb{H}^3$ by:
$$\matrix
\Cal{E}(p,t) \hfill&=\opExp(t\hat{\mathj}_0(p)), \hfill\cr
\Cal{E}_\infty(p) \hfill&=\opExp(+\infty\hat{\mathj}_0(p)). \hfill\cr
\endmatrix$$
\noindent Let us denote by $W$ the set $\Cal{E}(U\times [0,\infty))$. The application $\Cal{E}$ defines a diffeomorphism of $V\times[0,\infty)$ onto $W$. Let us define
$\pi_1:V\times [0,\infty)\rightarrow V$ to be the projection onto the first coordinate, and let us define $\pi:W\rightarrow V$ such that:
$$
\pi = \pi_1\circ\Cal{E}^{-1}.
$$
\noindent In particular $\pi$ is smooth. Let us denote by $W_\infty$ the set $\Cal{E}_\infty$. The application $\Cal{E}_\infty$ defines a homeomorphism of $V$ onto $W_\infty$.
\medskip
\noindent We now have two possibilities. Either $\overline{T}_0$ is a point in $\partial_\infty\Bbb{H}^3$, or it is the closure in $\Bbb{H}^3\munion\partial_\infty\Bbb{H}^3$
of a geodesic in $\Bbb{H}^3$. If $\overline{T}_0$ is a point $\left\{t_0\right\}$ in $\partial_\infty\Bbb{H}^3$, then, since, for all $q\in V$, the intersection of
$\opExp(\hat{\mathj}_0(q))([0,\infty))$ with $\overline{T}_0$ is non-empty, we obtain:
$$
\Cal{E}_\infty(V) = \left\{t_0\right\}.
$$
\noindent This is absurd, since $\Cal{E}_\infty$ is a homeomorphism. Let us now assume that $\overline{T}_0$ is the closure in $\Bbb{H}^3\munion\partial_\infty\Bbb{H}^3$ of
a geodesic in $\Bbb{H}^3$. Let $\gamma:\Bbb{R}\rightarrow\Bbb{H}^3$ be a parametrisation of this geodesic. Let $\suite{I}{n}$ be a collection of disjoint subintervals of
$\Bbb{R}$ such that:
$$
\gamma(\Bbb{R})\minter W = \munion_{n\in\Bbb{N}} \gamma(I_n).
$$
\noindent For every $n$, the application $\pi\circ\gamma$ is smooth over $I_n$. Consequently, if we denote by $\mu$ the two dimensional measure generated over $V$ by the metric
$\hat{\mathj}_0^*g^\nu$, we obtain:
$$
\mu(\munion_{n\in\Bbb{N}}(\pi\circ\gamma)(I_n))=0.
$$
\noindent Let us denote by $\gamma(\pm\infty)$ the boundary of the image of $\gamma$ in $\Bbb{H}^3\munion\partial_\infty\Bbb{H}^3$. Since $\Cal{E}_\infty$ is a homeomorphism,
the set $\Cal{E}_\infty^{-1}(\gamma(\pm\infty))$ consists of at most $2$ points. Consequently:
$$
\mu(\Cal{E}_\infty^{-1}(\gamma(\pm\infty)))=0.
$$
\noindent Let $q\in V$ be an arbitrary point in $V$. Since the intersection of $\opExp(\hat{\mathj}_0(q))[0,\infty])$ with $T_0$ is non-empty, $q$ must be in the union of
$\munion_{n\in\Bbb{N}}(\pi\circ\gamma)(I_n)$ with $\mu(\Cal{E}_\infty^{-1}(\gamma(\pm\infty))$. Consequently, $V$ is the union of these two sets and is thus of measure
zero. This is absurd, and it follows that $\hat{\Sigma}'_0$ is not the Gauss lifting of a k-surface, and is consequently a tube about a geodesic in $\Bbb{H}^3$. Let us
denote this geodesic by $\Gamma$, and let us denote by $\overline{\Gamma}$ the closure of $\Gamma$ in $\Bbb{H}^3\munion\partial_\infty\Bbb{H}^3$. A similar reasoning permits us
to conclude that $T_0$ coincides with $\Gamma$.
\medskip
\noindent For all $n$, $M_n\Gamma_{0,\infty}$ is the unique vertical geodesic joining $-w_n/\lambda_n$ to $\infty$. Since\break $(M_n\overline{\Gamma}_{0,\infty})_\ninn$
converges towards $\overline{T}_0=\overline{\Gamma}$ in the Haussdorf topology, and since $\overline{\Gamma}$ passes by $(0,1)$, we conclude that $\overline{\Gamma}$ is the
unique vertical geodesic passing by $(0,1)$. Consequently:
$$
\overline{\Gamma} = \overline{\Gamma}_{0,\infty}.
$$
\noindent It follows that $(w_n/\lambda_n)_\ninn$ tends towards $0$, and the first result follows. Since $(w_n/\lambda_n)_\ninn$ converges towards $0$, it follows that
the sequence $\suite{M}{n}$ converges to the identity. Consequently, the sequence of immersions $\suite{\hat{\mathi}}{n}=(M_n^{-1}\circ\hat{\mathj}_n)_\ninn$ converges towards
$\hat{\mathj}_0$. In otherwords, the sequence of immersed surfaces $(\Sigma_n,p_n)_\ninn$ converges in the Cheeger/Gromov topology towards $(\Sigma'_0,p_0)$, which is itself a tube about $\Gamma_{0,\infty}$, and the second result follows.\qed
\newsubhead{Tubes of Finite Order About Geodesics}
\noindent By continuing to identify $\Bbb{H}^3$ with the upper half space $\Bbb{C}\times\Bbb{R}^+$ and $T\Bbb{H}^3$ with $(\Bbb{C}\times\Bbb{R})_{(\Bbb{C}\times\Bbb{R}^+)}$,
we define $n_{0,1}\in U_{(0,1)}\Bbb{H}^3$ by:
$$
n_{(0,1)} = (1,0)_{(0,1)}.
$$
\noindent Let $N_{0,\infty}$ be the normal circle bundle over $\Gamma_{0,\infty}$ in $\Bbb{H}^3$. Let $N_{0,\infty}(0,1)$ be the fibre above $(0,1)$. Every subsequence of
$(\hat{\mathi}_n(p_n))_\ninn$ has a subsubsequence converging to a point in $N_{0,\infty}(0,1)$. It follows that there exists a sequence $(R_n)_\ninn$ of rotations about
$\Gamma_{0,\infty}$ such that the sequence $(R_n\circ\hat{\mathi}_n(p_n))_\ninn$ converges to $n_{(0,1)}$. By replacing $A_n$ with $R_n\circ A_n$ for all $n$, we may assume
that $(\hat{\mathi}_n(p_n))_\ninn$ converges to $n_{(0,1)}$.
\medskip
\noindent We now obtain the following stronger version of the previous result:
\proclaim{Lemma \nextprocno}
\noindent After extracting a subsequence, $(\hat{\Sigma}_n,p_n)_\ninn$ converges to a tube of order $k$ about $\Gamma_{0,\infty}$ with base point $n_{(0,1)}$.
\endproclaim
\proclabel{LemmeConvergeVersUnTubeDOrdreK}
\proof Let $N_{0,\infty}$ be the normal unitary bundle over $\Gamma_{0,\infty}$ in $U\Bbb{H}^3$. The sequence\break $(\hat{\Sigma}_n,p_n)_\ninn$ converges to a tube about
$\Gamma_{0,\infty}$ in the Cheeger/Gromov topology. Let\break $(\hat{\Sigma}_0,p_0)=(S_0,\hat{\mathi}_0,p_0)$ be the limit of this sequence. Since
$\hat{\mathi}_0:S_0\rightarrow N_{0,\infty}$ is a local isometry between two complete surfaces, there exists $m\in\Bbb{N}\munion\left\{\infty\right\}$ such that $\hat{\mathi}_0$
is an m-fold covering of $N_{0,\infty}$.We thus aim to show that $m=k$.
\medskip
\noindent As before, after applying an isometry of $\Bbb{H}^3$ if necessary, by lemma \procref{LemmeTrivDesRevetementsRamifiees}, we may find a chart $(z,U,\Bbb{D})$ of $S$ about
$p_0$ and $k\in\Bbb{N}$ such that the following diagram commutes:
$$\commdiag{
U & & \cr
\mapdown_z & \arrow(3,-2)\rt{\varphi} & \cr
\Bbb{D} & \mapright^{z\mapsto z^k} & \hat{\Bbb{C}} \cr
}$$
\noindent For all $n$, we define $z_n$ and $D_n$ by:
$$\matrix
z_n \hfill&= \lambda_n^{-1/k}z,\hfill\cr
D_n \hfill&= \left\{ z\in\Bbb{C}\text{ s.t. }0 < \left|z\right| < \lambda_n^{-1/k}\right\}.
\endmatrix$$
\noindent For all $n$, $(z_n, U\setminus\left\{p_0\right\}, D_n)$ defines a chart of $S\setminus\Cal{P}$ such that, if $\overrightarrow{n}$ is the Gauss-Minkowski mapping which
sends $U\Bbb{H}^3$ to $\partial_\infty\Bbb{H}^3$, then, for all $p\in U\setminus\left\{p_0\right\}$:
$$
\overrightarrow{n}\circ\hat{\mathi}_n(p) = z_n(p)^k.
$$
\noindent Using these charts, we will construct a sequence of pointed tubes of order $k$ about $\Gamma_{0,\infty}$ which converges to $(\hat{\Sigma}_0,p_0)$ in the
Cheeger/Gromov topology. The Haussdorf property of the Cheeger/Gromov topology will then permit us to conclude.
\medskip
\noindent We begin by constructing a number of coordinate charts that are well adapted to our problem. To begin with, we may assume that $(S_0,\hat{\mathi}_0^*g^\nu)$ is
equal either to $S^1\times\Bbb{R}$ or to $\Bbb{R}^2$, both of these spaces being furnished with the canonical Euclidean metric. Let $NN_{0,\infty}$ be the normal bundle over
$N_{0,\infty}$ in $U\Bbb{H}^3$. $NN_{0,\infty}$ is trivial and there exists a canonical vector bundle isomorphism
$\tau:(S^1\times\Bbb{R})\times\Bbb{R}^3\rightarrow NN_{0,\infty}$ which is unique up to composition with an element of $SO(3)$. For $\epsilon\in\Bbb{R}^+$ we define:
$$
N_\epsilon N_{0,\infty} = \left\{ v\in NN_{0,\infty}\text{ s.t. }\|v\|<\epsilon\right\}.
$$
\noindent Let $\opExp:NN_{0,\infty}\rightarrow U\Bbb{H}^3$ be the exponential mapping. Since $U\Bbb{H}^3$ is homogeneous, there exists $\epsilon\in\Bbb{R}^+$ such that the
restriction of $\opExp$ to $N_\epsilon N_{0,\infty}$ is a diffeomorphism onto its image. We define the mapping $\omega$ by $\omega=\opExp\circ\tau$ and we define
$\Omega\subseteq U\Bbb{H}^3$ by:
$$
\Omega = \omega((S^1\times\Bbb{R})\times B_\epsilon(0)).
$$
\noindent We may assume that $\omega$ sends the origin to $\hat{\mathi}_0(p_0)=n_{(0,1)}$. The triple $(\omega^{-1},\Omega, (S^1\times\Bbb{R})\times B_\epsilon(0))$ provides a coordinate
chart of $U\Bbb{H}^3$ which is well adapted to our problem. Let $\pi_1:(S^1\times\Bbb{R})\times B_\epsilon(0)\rightarrow (S^1\times\Bbb{R})\times\left\{0\right\}$ be the
projection onto the first factor.
\medskip
\noindent Let us denote by $\overrightarrow{n}$ the Gauss-Minkowski mapping which sends $U\Bbb{H}^3$ onto $\partial_\infty\Bbb{H}^3$. We identify $\partial_\infty\Bbb{H}^3$ with
the Riemann sphere $\hat{\Bbb{C}}$. Since the group of isometries of $\Bbb{H}^3$ which preserve $\Gamma_{0,\infty}$ acts transitively over $N_{0,\infty}$, by reducing
$\epsilon$ if necessary, we may assume that:
$$
\overrightarrow{n}(\Omega) = \Bbb{C}^* = \hat{\Bbb{C}}\setminus\left\{0,\infty\right\}.
$$
\noindent The application $\overrightarrow{n}$ defines a diffeomorphism between $N_{0,\infty}$ and $\Bbb{C}^*$. Let us denote the inverse of this mapping by $\pi_{0,\infty}$.
\medskip
\noindent Let $(\psi_n)_\ninn$ be a sequence of convergence mappings of $(\hat{\Sigma}_n,p_n)_\ninn$ with respect to $(\hat{\Sigma}_0,p_0)$. Let $R\in\Bbb{R}^+$ be a positive
real number. Let $N\in\Bbb{N}$ be such that, for all $n\geqslant N$:
\medskip
\myitem{(i)} the restriction of $\psi_n$ to $B_{R+1}(p_0)$ is a diffeomorphism onto its image,
\medskip
\myitem{(ii)} $B_{R+1/2}(p_n)$ is contained in the image of $B_{R+1}(p)$ under $\psi$, and
\medskip
\myitem{(iii)} $(i_n\circ\psi_n)(B_{R+1}(p_0))$ is contained in $\Omega$.
\medskip
\noindent For all $p\in\Bbb{N}$, let $\|\cdot\|_{C^p,R}$ be the $C^p$ norm over $B_{2R}(p_0)$. Since $\Sigma_0$ is a tube, $\hat{\mathi}_0$
takes values in $\omega((S^1\times\Bbb{R})\times\left\{0\right\})$. Consequently:
$$
\pi_1\circ(\omega^{-1}\circ\hat{\mathi}_0) = (\omega^{-1}\circ\hat{\mathi}_0).
$$
\noindent Thus:
$$
(\|\pi_1\circ(\omega^{-1}\circ\hat{\mathi}_n\circ\psi_n) - (\omega^{-1}\circ\hat{\mathi}_n\circ\psi_n)\|_{C^p,R})_\ninn\rightarrow 0.
$$
\noindent Applying $(\overrightarrow{n}\circ\omega)$, we obtain:
$$
(\|(\overrightarrow{n}\circ\omega)\circ\pi_1\circ(\omega^{-1}\circ\hat{\mathi}_n\circ\psi_n) -
(\overrightarrow{n}\circ\omega)\circ(\omega^{-1}\circ\hat{\mathi}_n\circ\psi_n)\|_{C^p,R})_\ninn\rightarrow 0.
$$
\noindent Using the definition of $\varphi_n$, we obtain:
$$
(\|(\overrightarrow{n}\circ\omega)\circ\pi_1\circ(\omega^{-1}\circ\hat{\mathi}_n\circ\psi_n) - (\varphi_n\circ\psi_n)\|_{C^p,R})_\ninn\rightarrow 0.
$$
\noindent We now apply $(\omega^{-1}\circ\pi_{0,\infty})$ to obtain:
$$
(\|(\omega^{-1}\circ\pi_{0,\infty})\circ(\overrightarrow{n}\circ\omega)\circ\pi_1\circ(\omega^{-1}\circ\hat{\mathi}_n\circ\psi_n) -
(\omega^{-1}\circ\pi_{0,\infty})\circ(\varphi_n\circ\psi_n)\|_{C^p,R})_\ninn\rightarrow 0.
$$
\noindent Using the fact that the restriction of $\pi_{0,\infty}\circ\overrightarrow{n}$ to $N_{0,\infty}$ is equal to the identity, we obtain:
$$
(\|(\pi_1\circ\omega^{-1})\circ(\hat{\mathi}_n\circ\psi_n) - (\omega^{-1}\circ\pi_{0,\infty})\circ(\varphi_n\circ\psi_n)\|_{C^p,R})_\ninn\rightarrow 0.
$$
\noindent By considering the metric $(\omega^*g^\nu)$ over $(S^1\times\Bbb{R})\times B_\epsilon(0)$, the first and last of these limits
permit us to obtain:
$$\matrix
(\|(\hat{\mathi}_n\circ\psi_n)^*(\pi_1\circ\omega^{-1})^*(\omega^*g^\nu)-(\hat{\mathi}_n\circ\psi_n)^*(\omega^{-1})^*(\omega^*g^nu)\|_{C^p,R})_\ninn\hfill&\rightarrow 0,\hfill\cr
(\|(\hat{\mathi}_n\circ\psi_n)^*(\pi_1\circ\omega^{-1})^*(\omega^*g^\nu) - (\varphi_n\circ\psi_n)^*(\omega^{-1}\circ\pi_{0,\infty})^*(\omega^*g^\nu)\|_{C^p,R})_\ninn\hfill&\rightarrow 0.\hfill\cr
\endmatrix$$
\noindent Thus:
$$
(\|(\hat{\mathi}_n\circ\psi_n)^*g^\nu - (\varphi_n\circ\psi_n)^*(\pi_{0,\infty})^*g^\nu\|_{C^p,R})_\ninn\rightarrow 0.
$$
\noindent Consequently, since $((\hat{\mathi}_n\circ\psi_n)^*g^\nu)_\ninn$ converges to $\hat{\mathi}_0^*g^\nu$, we obtain:
$$
(\|\hat{\mathi}_0^*g^\nu - (\varphi_n\circ\psi_n)^*(\pi_{0,\infty})^*g^\nu\|_{C^p,R})_\ninn\rightarrow 0.
$$
\noindent Let us define $\alpha_k:\Bbb{C}^*\rightarrow\Bbb{C}^*$ by:
$$
\alpha_k(z) = z^k.
$$
\noindent We define $g_k$ over $\Bbb{C}^*$ by:
$$
g_k = \alpha_k^*\pi_{0,\infty}^*g^\nu.
$$
\noindent Since $\varphi_n=\alpha_k\circ z_n$, we have:
$$
(\|\hat{\mathi}_0^*g^\nu - (z_n\circ\psi_n)^*g_k\|_{C^p,R})_\ninn\rightarrow 0.
$$
\noindent Since the radius of $D_n$ about $0$ tends to infinity as $n$ tends to infinity, it follows that there exists $M$ such that for $n\geqslant M$:
$$
\psi_n(B_{R+1}(p_0))\subseteq\Omega_n.
$$
\noindent Consequently, for all $n$, the function $(z_n\circ\psi_n)$ is defined and is smooth over $B_{R+1}(p_0)$. Moreover, for $n\geqslant M$, the restriction of $(z_n\circ\psi_n)$ to
$B_{R+1}(p_0)$ is a diffeomorphism onto its image. Since $R\in\Bbb{R}^+$ is arbitrary, it follows that $(z_n\circ\psi_n)_{n\in\Bbb{N}}$ defines a sequence of convergence
mappings for the sequence $(\Bbb{C}^*, \varphi_n(p_n), g_k)_\ninn$ with respect to the limit $(S_0,p_0,\hat{\mathi}_0^*g^\nu)$.
\medskip
\noindent Since $(\varphi_n(p_n))_\ninn=(\overrightarrow{n}\circ\hat{\mathi}_n(p_n))_\ninn$ converges to $\overrightarrow{n}(n_{0,1})=1$ and since the\break Cheeger/Gromov topology is Haussdorf, it follows that $(S_0,p_0,\hat{\mathi}_0^*g^\nu)$ is isometric to\break $(\Bbb{C}^*,1,g_k)$.
By considering, for example, the length of the shortest homotopically non-trivial curve in $\Bbb{C}^*$, we find that, for $k\neq k'$,
the manifolds $(\Bbb{C}^*,1,g_k)$ and $(\Bbb{C}^*,1,g_{k'})$ are not isometric. Consequently $m=k$, and the result now follows.\qed
\medskip
\noindent Since the same result holds for every subsequence of $(\hat{\Sigma}_n,p_n)_\ninn$, we obtain the following stronger version of this result:
\proclaim{Corollary \nextprocno}
\noindent $(\hat{\Sigma}_n,p_n)_\ninn$ converges to a tube of order $k$ about $\Gamma_{0,\infty}$ with base point $n_{(0,1)}$.
\endproclaim
\proof We assume the contrary. Let $(T_0,p_0)$ be the tube of order $k$ about $\Gamma_{0,\infty}$ with base point $n_{(0,1)}$. We may assume that there exists a
neighbourhood $\Omega$ of $(T_0,p_0)$ in the Cheeger/Gromov topology such that, after extraction of a subsequence, for all $n$:
$$
(\hat{\Sigma}_n,p_n)\notin\Omega
$$
\noindent However, by the preceeding result, there exists a subsequence of $(\hat{\Sigma}_n,p_n)_\ninn$ which converges towards $(T_0,p_0)$. We thus have a contradiction and
the result now follows.\qed
\medskip
\noindent Expressing this result in terms of graphs over tubes, we obtain:
\proclaim{Lemma \nextprocno}
\noindent Let $r\in\Bbb{R}^+$ be a positive real number. There exists $N\in\Bbb{N}$ such that for all $n\geqslant N$, the pointed surface $(\hat{\Sigma}_n,p_n)$ is locally
a graph over a tube about $\Gamma_{0,\infty}$ of order $k$ and of half length $r$.
\medskip
\noindent Moreover, if for all $n\geqslant N$, we denote by $\lambda_{n}$ the graph function of $(\hat{\Sigma}_n,p_n)$ over $S^1\times]-r,r[$, then $(\lambda_n)_{n\geqslant N}$
converges to $0$ in the $C^\infty_\oploc$ topology.
\endproclaim
\proof Let $\opExp:TU\Bbb{H}^3\rightarrow U\Bbb{H}^3$ be the exponential mapping over $U\Bbb{H}^3$. Let $N_{0,\infty}$ be the normal circle bundle over $\Gamma_{0,\infty}$ in
$U\Bbb{H}^3$. Let $NN_{0,\infty}$ be the normal bundle over $N_{0,\infty}$ in $TU\Bbb{H}^3$. For $\epsilon\in\Bbb{R}^+$, we define $N_\epsilon N_{0,\infty}$ by:
$$
N_\epsilon N_{0,\infty} = \left\{ v \in NN_{0,\infty} | \|v\| < \epsilon \right\}.
$$
\noindent Since $U\Bbb{H}^3$ is homogeneous, there exists $\epsilon\in\Bbb{R}^+$ such that the restriction of $\opExp$ to $N_\epsilon N_{0,\infty}$ is a diffeomorphism onto its
image. Let us define $U\subseteq U\Bbb{H}^3$ by $U=\opExp(N_\epsilon N_{0,\infty})$. Let $\pi:U\rightarrow N_{0,\infty}$ be the orthogonal projection onto $N_{0,\infty}$.
\medskip
\noindent Let $T=(S^1\times\Bbb{R},\hat{\mathi}_0)$ be a tube of order $k$ about $\Gamma_{0,\infty}$. Let $p_0$ be the origin of $T$. By the preceeding result, we may assume
that $(\hat{\Sigma}_n,p_n)_{n\in\Bbb{N}}$ converges towards $(T,p_0)$ in the Cheeger/Gromov topology. Let $(\varphi_n)_\ninn$ be a sequence of convergence mappings for
$(\hat{\Sigma}_n,p_n)_\ninn$ with respect to $(T,p_0)$.
\medskip
\noindent Let us define the application $\alpha:(S^1\times\Bbb{R})\times(S^1\times\Bbb{R})\rightarrow(S^1\times\Bbb{R})\times N_{0,\infty}$ by:
$$
\alpha(x,y) = (x,\hat{\mathi}_0(y)).
$$
\noindent Let $\Delta$ be the diagonal in $(S^1\times\Bbb{R})\times(S^1\times\Bbb{R})$:
$$
\Delta = \left\{(x,x) | x\in (S^1\times\Bbb{R}) \right\}.
$$
\noindent For all $\rho\in\Bbb{R}^+$, let $B_\rho(\Delta)$ be the tubular neighbourhood of radius $\rho$ about $\Delta$ in $(S^1\times\Bbb{R})\times(S^1\times\Bbb{R})$. Since
$\alpha$ is a local diffeomorphism, it follows by the homogeneity of $(S^1\times\Bbb{R})\times(S^1\times\Bbb{R})$ that there exists $\rho\in\Bbb{R}^+$ such that the
restriction of $\alpha$ to $B_\rho(\Delta)$ is a diffeomorphism onto its image. We will use this mapping to unravel other mappings that wrap $k$ times round $N_{0,\infty}$. We
define $V$ by:
$$
V = \alpha(B_\rho(\Delta)).
$$
\noindent Let $p_1,p_2:(S^1\times\Bbb{R})\times(S^1\times\Bbb{R})\rightarrow S^1\times\Bbb{R}$ be the projections onto the first and second factors respectively.
\medskip
\noindent Let $R>0$ be a positive real number such that $T_{2r}\subseteq B_R(p_0)\subseteq T$. Let $N_1\in\Bbb{N}$ be such that for all $n\geqslant N_1$:
\medskip
\myitem{(i)} the restriction of $\varphi_n$ to $B_{R+1}(p_0)$ is a diffeomorphism onto its image,
\medskip
\myitem{(ii)} $(i_n\circ\varphi_n)(B_{R+1}(p_0))$ is contained within $U$, and
\medskip
\myitem{(iii)} for all $x\in B_{R+1}(p_0)$:
$$
(x,\pi\circ\hat{\mathi}_n\circ\varphi_n(x))\in V.
$$
\noindent For all $n\geqslant N_1$, we define $\beta_n:B_{R+1}(p_0)\rightarrow S^1\times\Bbb{R}$ by:
$$
\beta_n(x) = p_2\circ\alpha^{-1}(x,\pi\circ\hat{\mathi}_n\circ\varphi_n(x)).
$$
\noindent We define $\beta_0:B_{R+1}(p_0)\rightarrow S^1\times\Bbb{R}$ by:
$$
\beta_0(x) = p_2\circ\alpha^{-1}(x,\pi\circ\hat{\mathi}_0(x)).
$$
\noindent We obtain:
$$
\beta_0(x) = x.
$$
\noindent Since $(\beta_n)_{n\geqslant N_1}$ converges to $\beta_0$ in the $C^\infty_\oploc$ topology, it follows by the common sense lemmata
\procref{LemmeInjectiviteLorsquOnSApprocheDuLimite} and \procref{LemmeCompacteDansLImageDeLaLimite} that there exists $N_2\geqslant N_1$ such that for all $n\geqslant N_2$:
\medskip
\myitem{(i)} the restriction of $\beta_n$ to $B_{R+1/2}(p_0)$ is a diffeomorphism onto its image, and
\medskip
\myitem{(ii)} $(S^1\times ]-3r,3r[)$ is contained in $\beta_n(B_{R+1/2}(p_0))$.
\medskip
\noindent For all $n\geqslant N_2$, we define $\psi_n$ and $\lambda_n$ by:
$$\matrix
\psi_n \hfill&= \beta_n^{-1}|_{S^1\times]-2r,2r[},\hfill\cr
\lambda_n \hfill&= \opExp^{-1}\circ\hat{\mathi}_n\circ\varphi_n\circ\psi_n.\hfill\cr
\endmatrix$$
\noindent For all $x\in S^1\times]-2r,2r[$ and for all $n\geqslant N_2$, we have:
$$\matrix
& \beta_n\circ\psi_n(x) \hfill&= x\hfill\cr
\Rightarrow\hfill& p_2\circ\alpha^{-1}(\psi_n(x),\pi\circ\hat{\mathi}_n\circ\varphi_n\circ\psi_n(x)) \hfill&= x.\hfill\cr
\endmatrix$$
\noindent Since $p_1\circ\alpha=p_1$, and thus $p_1\circ\alpha^{-1}=p_1$, we obtain:
$$\matrix
&\alpha^{-1}(\psi_n(x),\pi\circ\hat{\mathi}_n\circ\varphi_n\circ\psi_n(x)) \hfill&= (\psi_n(x),x) \hfill\cr
\Rightarrow\hfill&(\psi_n(x),\pi\circ\hat{\mathi}_n\circ\varphi_n\circ\psi_n(x)) \hfill&= (\psi_n(x), \hat{\mathi}_0(x)) \hfill\cr
\Rightarrow\hfill&\pi\circ\hat{\mathi}_n\circ\varphi_n\circ\psi_n(x) \hfill&=\hat{\mathi}_0(x). \hfill\cr
\endmatrix$$
\noindent Thus:
$$\matrix
&\pi\circ\opExp^{-1}\circ\hat{\mathi}_n\circ\varphi_n\circ\psi_n \hfill&= \hat{\mathi}_0 \hfill\cr
\Rightarrow\hfill&\pi\circ\lambda_n \hfill&= \hat{\mathi}_0.\hfill\cr
\endmatrix$$
\noindent Thus $\lambda_n$ is a section of $\hat{\mathi}_0^*NN_{0,\infty}$.
\medskip
\noindent Since $(\beta_n)_\ninn$ converges to the identity, it follows that $(\beta_n(0))_\ninn$ converges to $0$. For all $n\geqslant N_2$, let
$\phi_n:S^1\times\Bbb{R}\rightarrow S^1\times\Bbb{R}$ be the unique conformal mapping which sends $0$ to $\beta_n(0)$. There exists $N_3\geqslant N_2$ such that
for all $n\geqslant N_3$:
$$
\phi_n(S^1\times ]-r,r[)\subseteq S^1\times ]-2r,2r[.
$$
\noindent For all $n\geqslant N_3$, we define:
$$\matrix
\psi_n' \hfill&= \varphi_n\circ\psi_n\circ\phi_n, \hfill\cr
\lambda_n' \hfill&= \lambda_n\circ\phi_n, \hfill\cr
\hat{\mathi}'_{0,n}\hfill&= \hat{\mathi}_0\circ\phi_n. \hfill\cr
\endmatrix$$
\noindent For all $n\geqslant N_3$, we obtain:
\medskip
\myitem{(i)} $(S^1\times]-r,r[,\hat{\mathi}'_{0,n})$ is a tube of order $k$ about $\Gamma_{0,\infty}$,
\medskip
\myitem{(ii)} $\psi_n':(S^1\times]-r,r[,0)\rightarrow(S_n,p_n)$ is a diffeomorphism onto its image,
\medskip
\myitem{(iii)} $\lambda'_n$ is a section of $(\hat{\mathi}'_{0,n})^*NN_{0,\infty}$ over $S^1\times\Bbb{R}$, and
\medskip
\myitem{(iv)} $\hat{\mathi}_n\circ\psi_n' = \opExp\circ\lambda_n'$.
\medskip
\noindent Consequently, for all $n\geqslant N_3$, the immersed surface $(\hat{\Sigma}_n,p_n)$ is locally a graph over a tube about $\Gamma_{0,\infty}$ of order $k$ and of
half length $n$, and the first result follows.
\medskip
\noindent Moreover, we find that $(\lambda_n')_\ninn$ converges to $0$ over $S^1\times]-r,r[$ in the $C^\infty_\oploc$ topology, and the second result follows.\qed
\medskip
\noindent Since the property of being locally a graph over $\Gamma_{0,\infty}$ is invariant under isometries of $\Bbb{H}^3$ which preserve $\Gamma_{0,\infty}$, we immediately
obtain the following result:
\proclaim{Corollary \nextprocno}
\noindent Let $r,\epsilon\in\Bbb{R}^+$ be positive real numbers. There exists $N\in\Bbb{N}$ such that for $n\geqslant N$ the pointed immersed surface $(\hat{\Sigma},p_n)$
is locally a graph over a tube about $\Gamma_{0,\infty}$ of order $k$ and of half length $r$ and if $\lambda_n$ is the graph function of $(\hat{\Sigma},p_n)$
over $S^1\times]-r,r[$ then $\|\lambda_n\|<\epsilon$. Moreover $(\lambda_n)_\ninn$ tends to $0$ in the $C^\infty_\oploc$ topology.
\endproclaim
\noindent If $\epsilon$ is sufficiently small, then the graph functions and the graph diffeomorphisms are unique. Since the same result holds for every sequence of points in $S\setminus\Cal{P}$ which tends towards $p_0$, we obtain:
\proclaim{Corollary \nextprocno}
\noindent Let $r,\epsilon\in\Bbb{R}^+$ be positive real numbers. There exists an open set $\Omega$ of $p_0$ in $S$ such that if $p\in\Omega\setminus\left\{p_0\right\}$, then
$(\hat{\Sigma},p)$ is locally a graph over a tube about $\Gamma_{0,\infty}$ of order $k$ and of half length $r$ and if $\lambda_p$ is the graph function of $(\hat{\Sigma},p)$
over $S^1\times]-r,r[$ then $\|\lambda_p\|<\epsilon$. Moreover $\lambda_p$ tends to $0$ in the $C^\infty_\oploc$ topology is $p$ tends to
$p_0$.
\endproclaim
\proclabel{CorLocallyGraphsOverTubes}
\noindent By glueing these graphs together, we now obtain theorem \procref{PresentationChIIIRevetementsRamifiees}:
\proclaim{Theorem \procref{PresentationChIIIRevetementsRamifiees}}
\noindent Let $S$ be a Riemann surface. Let $\Cal{P}$ be a discrete subset of $S$ such that $S\setminus\Cal{P}$ is hyperbolic. Let $\varphi:S\rightarrow\hat{\Bbb{C}}$ be a ramified covering having critical points in $\Cal{P}$. Let $\kappa$ be a real number in $(0,1)$. Let
$i:S\setminus\Cal{P}\rightarrow\Bbb{H}^3$ be the unique solution to the Plateau problem $(S\setminus\Cal{P},\varphi)$ with constant Gaussian
curvature $\kappa$. Let $\hat{\Sigma}=(S\setminus\Cal{P},\hat{\mathi})$ be the Gauss lifting of $\Sigma$.
\medskip
\noindent Let $p_0$ be an arbitrary point in $\Cal{P}$. If $\varphi$ has a critical point of order $k$ at $p_0$, then $\hat{\Sigma}$ is asymptotically tubular of order
$k$ at $p_0$.
\endproclaim
\proof As in the proof of the preceeding lemma, let $\epsilon$ be such that the restriction of $\opExp$ to $N_\epsilon N_{0,\infty}$ is a diffeomorphism onto its
image. Let us define $U\subseteq U\Bbb{H}^3$ by $U=\opExp(N_\epsilon N_{0,\infty})$. Let $\pi:U\rightarrow N_{0,\infty}$ be the orthogonal projection onto $N_{0,\infty}$.
\medskip
\noindent Let $r$ be a positive real number. By corollary \procref{CorLocallyGraphsOverTubes}, there exists a connected neighbourhood $\Omega$ of $p_0$ in $S$ such that
if $p\in\Omega\setminus\left\{p_0\right\}$, then $(\hat{\Sigma},p)$ is locally a graph over a tube about $\Gamma_{0,\infty}$ of order $k$ and of half length $2r$ and if
$\lambda$ is the graph function of $(\hat{\Sigma},p)$ over $S^1\times(-2r,2r)$, then $\|\lambda\|<\epsilon$.
\medskip
\noindent By using foliations, we will construct a chart over an open set about $p_0$ which is well adapted to our problem. Let $\Cal{F}$ be the canonical circle foliation
of $N_{0,\infty}$ arising from its structure as a circle bundle over $\Gamma_{0,\infty}$. Let $\Cal{F}'$ be a the canonical circle foliation of $S^1\times (-2r,2r)$.
\medskip
\noindent For $p$ a point in $\Omega\setminus\left\{p_0\right\}$, let $T_p=(S^1\times(-2r,2r),\hat{\mathi}_p)$ be the tube of order $k$ and of half length $2r$ over which
$(\Sigma,p)$ is a locally a graph. Let $\varphi_p:S^1\times(-2r,2r)\rightarrow S$ be the graph diffeomorphism of $(\Sigma,p)$ over $T_p$. We define:
$$
\Omega_p = \varphi_p(S^1\times(-r,r)).
$$
\noindent Let $\lambda_p:S^1\times(-2r,2r)\rightarrow\Bbb{R}$ be the graph function of $(\Sigma,p)$ over $T_p$. Since $\|\lambda_p\|<\epsilon$, we have $\hat{\mathi}(q)\in U$
for all $q\in\Omega_p$. Moreover, by corollary \procref{CorLocallyGraphsOverTubes}, we may assume that $(\pi\circ\hat{\mathi})^*g^\nu$
defines a metric over $\Omega_p$.
\medskip
\noindent We remark that by the uniqueness of graph diffeomorphisms:
$$
q\in\Omega_p\ \Leftrightarrow\ p\in\Omega_q.
$$
\noindent We define $\hat{\Omega}$ by:
$$
\hat{\Omega} = \munion_{p\in\Omega\setminus\left\{p_0\right\}}\Omega_p.
$$
\noindent Since $\Omega$ is connected, so is $\hat{\Omega}$. For all $p$, $(\varphi_p)_*\Cal{F}'$ defines a smooth circle foliation of $\Omega_p$. This circle foliation
coincides with $\pi^*\Cal{F}$. It thus follows that $\hat{\Omega}$ is foliated by $\pi^*\Cal{F}$. Using the definition of $g^\nu$, and recalling that $\pi$ is a $k$-fold
covering map, we find that every leaf of this foliation is of length $2\pi k\nu^{-1}$ with respect to the metric $(\pi\circ\hat{\mathi})^*g^\nu$.
\medskip
\noindent Let us define $L$ by:
$$
L = \hat{\Omega}/\pi^*\Cal{F}.
$$
\noindent $L$ is a smooth connected one-dimensional manifold without boundary and is thus diffeomorphic to an open interval $I=(a,b)$ in $\Bbb{R}$. The set $\hat{\Omega}$ is thus
diffeomorphic to a smooth circle bundle over $I$. We thus obtain a diffeomorphism $\varphi_1:S^1\times I\rightarrow\hat{\Omega}$. We will show that by modifying this
diffeomorphism we obtain the desired chart.
\medskip
\noindent Let $\delta$ be an arbitrary metric over $S$ compatible with its topology. For $p\in\Omega$, we define:
$$
\Delta(p) = \minf\left\{\delta(\varphi_p(e^{i\theta},t),p_0)\text{ s.t. }(e^{i\theta},t)\in S^1\times[-r,r]\right\}.
$$
\noindent By uniqueness of graph diffeomorphisms, for all $(e^{i\theta},t)\in S^1\times(-r,r)$, we obtain:
$$
\Delta(\varphi_p(e^{i\theta},t)) = \minf\left\{\delta(\varphi_p(e^{i\phi},s),p_0)\text{ s.t. }(e^{i\phi},s)\in S^1\times[t-r,t+r]\right\}.
$$
\noindent Consequently, $\Delta$ is continuous. Since $S^1\times [-r,r]$ is compact, we find that $\Delta(p)>0$ for all $p\in\Omega\setminus\left\{p_0\right\}$.
\medskip
\noindent Let $\Omega_1$ be connected neighbourhood of $p_0$ contained in $\hat{\Omega}$. Let us define $\Delta_1$ by:
$$
\Delta_1 = \minf\left\{\Delta(p)\text{ s.t. }p\in\partial\Omega_1\right\}.
$$
\noindent Let us define $\Omega_2$ by:
$$
\Omega_2 = \left\{p\in S\text{ s.t. }\delta(p,p_0)<\Delta_1\right\}.
$$
\noindent Since $q\in\Omega_p$ if and only if $p\in\Omega_q$, we obtain:
$$
\hat{\Omega}_2\minter\partial\Omega_1 = \emptyset.
$$
\noindent Consequently:
$$
\hat{\Omega}_2 \subseteq \Omega_1.
$$
\noindent Since $\hat{\Omega}_2$ is connected and foliated by $\pi^*\Cal{F}$, there exists an open subinterval $I'\subseteq I$ such that:
$$
\hat{\Omega}_2 = \varphi_1(S^1\times I').
$$
\noindent Moreover, since $p_0$ is contained in the closure of $\hat{\Omega}_2$, it follows that the closure of $I'$ in $I$ is not compact. Consequently, we may assume that there
exists $a'\in (a,b)$ such that:
$$
I' = (a',b).
$$
\noindent Since we may choose $\Omega_1$ arbitrarily small about $p_0$, we find that $\varphi_1(e^{i\theta},t)$ tends to $p_0$ as $t$ tends to $b$.
\medskip
\noindent Let $p:N_{0,\infty}\rightarrow\Gamma_{0,\infty}$ be the canonical projection. Let us also denote by $p$ the composition
$p\circ\pi:NN_{0,\infty}\rightarrow\Gamma_{0,\infty}$. Let $t_0$ be an arbitrary point in $(a,b)$. Let $\gamma:\Bbb{R}\rightarrow\Gamma_{0,\infty}$ be a unit speed
parametrisation of $\Gamma_{0,\infty}$ such that $\gamma(t)\rightarrow 0$ as $t$ tends to $+\infty$ and:
$$
\gamma(0) = (p\circ\hat{\mathi}\circ\varphi_1)(e^{i\theta},t_0).
$$
\noindent Since $\varphi_1$ respects the foliation $\pi^*\Cal{F}$, the mapping $t\mapsto(\gamma^{-1}\circ p\circ\hat{\mathi}\circ\varphi_1)(e^{i\theta},t)$ is independant
of $\theta$ and is everywhere a local diffeomorphism. Consequently, it defines a diffeomorphism. By lemma \procref{RevetementsRamifies}, $\hat{\mathi}(p)\rightarrow 0$ as
$p\rightarrow p_0$. Consequently $(\gamma^{-1}\circ p\circ\hat{\mathi}\circ\varphi_1)(e^{i\theta},t)\rightarrow +\infty$ as $t\rightarrow b$. We may thus reparametrise
$\varphi_1$ to obtain a diffeomorphism $\varphi_2:S^1\times(a_1,\infty)\rightarrow\hat{\Omega}$ such that, for all $\theta$ and for all $t$:
$$
(p\circ\hat{\mathi}\circ\varphi_2)(e^{i\theta},t) = \gamma(t).
$$
\noindent We define the vector fields $\partial_\theta$ and $\partial_t$ over $S^1\times (-2r,2r)$ by:
$$\matrix
\partial_\theta(e^{i\theta},t) \hfill&= [\phi\mapsto (e^{i(\theta + \phi)},t)], \hfill\cr
\partial_t(e^{i\theta}, t) \hfill&= [s\mapsto (e^{i\theta},t+s)]. \hfill\cr
\endmatrix$$
\noindent For all $p$, we may orient $\varphi_p$ in such a manner that there exists $T_p\in\Bbb{R}$ such that for all $(t,e^{i\theta})$:
$$
(p\circ\hat{\mathi}\circ\varphi_p)(t,e^{i\theta}) = \gamma(t+T_p).
$$
\noindent We then define $X_p$ and $Y_p$ over $\Omega_p$ by:
$$\matrix
X_p \hfill&= (\varphi_p)_*\partial_\theta, \hfill\cr
Y_p \hfill&= (\varphi_p)_*\partial_t. \hfill\cr
\endmatrix$$
\noindent By the uniqueness of graph diffeomorphisms, for all $p,q\in\Omega$, the mapping $\varphi_p^{-1}\circ\varphi_q$ defined over $\varphi_q^{-1}(\Omega_p\minter\Omega_q)$
is an affine mapping (i.e. a rotation followed by a translation). Since $\varphi_p^{-1}\circ\varphi_q$ preserves orientation, we obtain, for all $p,q\in\Omega$:
$$\matrix
X_p|_{\Omega_p\minter\Omega_q} \hfill&= X_q|_{\Omega_p\minter\Omega_q}\hfill\cr
Y_p|_{\Omega_p\minter\Omega_q} \hfill&= Y_q|_{\Omega_p\minter\Omega_q}\hfill\cr
\endmatrix$$
\noindent We may thus define $X$ and $Y$ over the whole of $\hat{\Omega}$ such that, for all $p$:
$$
X|_{\Omega_p} = X_p,\ Y|_{\Omega_p} = Y_p.
$$
\noindent In particular $[XY]=0$. Let $\Phi$ and $\Psi$ be the flows of $X$ and $Y$. $(\Phi_t)_{t\in\Bbb{R}}$ is a flow that moves along the leaves of
the foliation $\pi^*\Cal{F}$ with speed $k\nu^{-2}$ with respect to the metric $(\pi\circ\hat{\mathi})^*g^\nu$. In particular $\Phi_t$ is defined over $\hat{\Omega}$ for all
$t\in\Bbb{R}$. Moreover, since every leaf of $\pi^*\Cal{F}$ is of length $2\pi k\nu^{-1}$, it follows that $(\Phi_t)_{t\in\Bbb{R}}$ is periodic with period $2\pi$. Let
$(G_t)_{t\in\Bbb{R}}$ be the geodesic flow along $\Gamma_{0,\infty}$ in the positive direction (i.e. towards $0$) with constant speed $k\nu^{-2}$. For all $t\geqslant 0$, the
following diagram commutes:
$$\commdiag{
\hat{\Omega} &\mapright^{\Psi_t} &\hat{\Omega} \cr
\mapdown_p& &\mapdown_p\cr
\Gamma_{0,\infty} &\mapright^{G_t} &\Gamma_{0,\infty}\cr}
$$
\noindent It follows that $\Psi_t$ is defined over $\hat{\Omega}$ for all $t\geqslant 0$. We define $\varphi_3:\Bbb{R}\times(0,\infty)\rightarrow\hat{\Omega}$
by:
$$
\varphi_3(s,t) = \Phi_{2\pi\nu^{-1}ks}\Psi_{t}(\varphi_1(0,t_0)).
$$
\noindent Since $(\Phi_t)_{t\in\Bbb{R}}$ is periodic with period $2\pi$, the mapping $\varphi_2$ quotients to an application
$\varphi:S^1\times(0,\infty)\rightarrow\hat{\Omega}$. We define $\hat{\mathj}$ and $\lambda$ by:
$$\matrix
\hat{\mathj} \hfill&= \pi\circ\hat{\mathi}\circ\varphi,\hfill\cr
\lambda \hfill&= \opExp^{-1}\circ\hat{\mathi}_0\circ\varphi.\hfill\cr
\endmatrix$$
\noindent For all $p\in\Omega$, we find that $\varphi_p^{-1}\circ\varphi|_{\varphi^{-1}(\Omega_p)}$ is the restriction of an affine transformation $\phi_p$ of $S^1\times\Bbb{R}$
to $\varphi^{-1}(\Omega_p)$. Consequently:
$$\matrix
\hat{\mathj}|_{\varphi^{-1}(\Omega_p)} \hfill&= \pi\circ\hat{\mathi}\circ\varphi|_{\varphi^{-1}(\Omega_p)}\hfill\cr
&=(\pi\circ\hat{\mathi}\circ\varphi_p)\circ(\varphi_p^{-1}\circ\varphi)|_{\varphi^{-1}(\Omega_p)}\hfill\cr
&=\hat{\mathi}_p\circ\phi_p.\hfill\cr
\endmatrix$$
\noindent We recall that $\hat{\mathi}_p$ is a locally conformal $k$-fold covering map. It thus follows that $\hat{\mathj}:S^1\times(0,\infty)\rightarrow N_{0,\infty}$ is
a locally conformal $k$-fold covering map. We thus have:
\medskip
\myitem{(i)} $(S^1\times(0,\infty),\hat{\mathj})$ defines a half tube of order $k$ about $\Gamma_{0,\infty}$,
\medskip
\myitem{(ii)} $\lambda$ is a section of $\hat{\mathj}^*NN_{0,\infty}$, and
\medskip
\myitem{(iii)} $\hat{\mathi}\circ\varphi = \opExp\circ\lambda$.
\medskip
\noindent We have thus shown that $\hat{\Sigma}$ is a graph over a half tube of order $k$ about $\Gamma_{0,\infty}$. Moreover, by corollary \procref{CorLocallyGraphsOverTubes},
and the uniqueness of graph functions, we find that for all $p\in\Bbb{R}$, $\|D^p\lambda(e^{i\theta},t)\|$ converges to $0$ as $t$ tends to $+\infty$, and the result follows.\qed
\newhead{Asymptotically Tubular Surfaces of Finite Order}
\newsubhead{Introduction}
\noindent In this section, we will prove theorem \procref{PresentationChIIISurfacesAsymptotiquementTubulaires}:
\proclaim{Theorem \procref{PresentationChIIISurfacesAsymptotiquementTubulaires}}
\noindent let $S$ be a surface and let $\Cal{P}\subseteq S$ be a discrete subset of $S$. Let $i:S\setminus\Cal{P}\rightarrow\Bbb{H}^3$ be
an immersion such that $\Sigma=(S\setminus\Cal{P},i)$ is a $k$-surface (and is thus the solution to a Plateau problem). Let
$\overrightarrow{n}:U\Bbb{H}^3\rightarrow\partial_\infty\Bbb{H}^3$ be the Gauss-Minkowski mapping which sends $U\Bbb{H}^3$ to
$\partial_\infty\Bbb{H}^3$. Let $\hat{\mathi}$ be the Gauss lifting of $i$ so that $\varphi=\overrightarrow{n}\circ\hat{\mathi}$ defines
the Plateau problem to which $i$ is the solution. Let $\Cal{H}$ be the holomorphic structure generated over $S\setminus\Cal{P}$ by
the local homeomorphism $\varphi$.
\medskip
\noindent Let $p_0$ be an arbitrary point in $\Cal{P}$, and suppose that $\Sigma$ is asymptotically tubular of order $k$ about $p_0$. Then
there exists a unique holomorphic structure $\tilde{\Cal{H}}$ over $(S\setminus\Cal{P})\munion\left\{p_0\right\}$ and a unique
holomorphic mapping $\tilde{\varphi}:(S\setminus\Cal{P})\munion\left\{p_0\right\}\rightarrow\hat{\Bbb{C}}$ such that $\tilde{\Cal{H}}$ and
$\tilde{\varphi}$ extend $\Cal{H}$ and $\varphi$ respectively. Moreover, $\tilde{\varphi}$ has a critical point of order $k$ at
$p_0$.
\endproclaim
\noindent This result will be proven in two stages. First, by using the properties of the modules of conformal rings, we obtain:
\proclaim{Lemma \nextprocno}
\noindent Let $S$ be a surface. Let $\Cal{P}$ be a discrete subset of $S$. Let $i:S\setminus\Cal{P}\rightarrow\Bbb{H}^3$ be an immersion such that the immersed surface
$\Sigma=(S\setminus\Cal{P},i)$ is a $k$-surface. Let $\hat{\Sigma}=(S,\hat{\mathi})$ be the Gauss lifting of $\Sigma$. Let $\overrightarrow{n}$ be the Gauss-Minkowski
mapping that sends $U\Bbb{H}^3$ into $\partial_\infty\Bbb{H}^3=\hat{\Bbb{C}}$. Let us define $\varphi=\overrightarrow{n}\circ\hat{\mathi}$. Let $\Cal{H}$ be the canonical
conformal structure over $\hat{\Bbb{C}}$.
\medskip
\noindent Let $p_0$ be a point in $\Cal{P}$. If $\hat{\Sigma}$ is asymptotically tubular of finite order near $p_0$, then $\varphi^*\Cal{H}$ extends to a unique conformal
structure on $(S\setminus\Cal{P})\munion\left\{p_0\right\}$.
\endproclaim
\proclabel{LemmaATSExtendConformalStructure}
\noindent Next, by showing that there exist $q_0$ such that $\varphi(p)$ tends to $q_0$ as $p$ tends to $p_0$, applying Cauchy's removeable singularity theorem, we obtain:
\proclaim{Lemma \nextprocno}
\noindent With the same hypotheses as in lemma \procref{LemmaATSExtendConformalStructure}, let $p_0$ be a point in $\Cal{P}$. If $\hat{\Sigma}$ is asymptotically
tubular of order $k$ near $p_0$, then $\varphi$ extends to a unique holomorphic function over $(S\setminus\Cal{P})\munion\left\{p_0\right\}$ having a critical point of order
$k$ at $p_0$.
\endproclaim
\proclabel{LemmaATSExtendFunction}
\noindent Theorem \procref{PresentationChIIISurfacesAsymptotiquementTubulaires} now follows as a direct corollary to these two lemmata.
\newsubhead{Conformal Rings}
\noindent In this section we will recall various properties of holomorphic rings. We define a {\emph (conformal) ring\/} to be a Riemann surface $A$ whose fundamental group is
isomorphic to $\Bbb{Z}$. For $R$ a real number greater than $1$, let us define the ring $A_R$ by:
$$
A_R = \left\{z | 1< \left|z\right|  < R \right\}.
$$
\noindent The uniformisation principal permits us to show that an arbitrary ring is biholomorphic to one of $\Bbb{C}^*$, $\Bbb{D}^*$ or $A_R$ for some $R\in (1,\infty)$.
Let $\Gamma$ be the familly of curves in $A$ which are freely homotopic to a generator of $\pi_1(A)$. For $g$ a conformal Riemannian metric over $A$, and for $\gamma\in\Gamma$
an arbitrary curve in $\Gamma$, we define $\opLen_g(\gamma)$ to be the length of $\gamma$ with respect to $g$, and we define $\Cal{L}_g(\Gamma)$ by:
$$
\Cal{L}_g(\Gamma) = \minf_{\gamma\in\Gamma}\opLen_g(\gamma).
$$
\noindent For $g$ a conformal metric over $A$, we define $\opArea_g(A)$ to be the area of $A$ with respect to $g$. We define $\opMod(A)$, the {\emph module\/} of $A$, by:
$$
\opMod(A) = \msup_{\opArea_g(A) = 1\atop g\text{ conformal}}\Cal{L}_g(\Gamma).
$$
\noindent By definition $\opMod(A)$ only depends on the conformal class of $A$. $\opMod(A)$ may be calculated in certain cases, and, in particular, we have the following
result:
\proclaim{Lemma \nextprocno}
\noindent For all $R\in (0,\infty)$:
$$\matrix
\opMod(A_R) \hfill&=\sqrt{\frac{2\pi}{\opLog(R)}},\hfill\cr
\opMod(S^1\times ]0,R[) \hfill&=\sqrt{\frac{2\pi}{R}}.\hfill\cr
\endmatrix$$
\endproclaim
\proclabel{LemmeCalculerLeModuleDUnAnneau}
\proof Let $g$ be a conformal metric of area $1$ over $A_R$. Let $g_{\opEuc}$ be the Euclidean metric over $A_R$ and let $\lambda:A_R\rightarrow (0,\infty)$ be such that:
$$
g = \lambda g_{\opEuc}.
$$
\noindent Using the Cauchy-Schwarz inequality, we obtain:
$$\matrix
\opArea_g(A_R) \hfill&= \int_1^R\int_0^{2\pi}\lambda r dr d\theta \hfill\cr
&\geqslant \int_1^R\frac{1}{2\pi r}\left(\int_0^{2\pi}r\lambda^{1/2}d\theta\right)^2dr \hfill\cr
&\geqslant \int_1^R\frac{1}{2\pi r}\Cal{L}_g(\Gamma)^2dr \hfill\cr
&=\frac{\opLog(R)}{2\pi}\Cal{L}_g(\Gamma)^2. \hfill\cr
\endmatrix$$
\noindent Since $\opArea_g(A_R)=1$, we obtain:
$$
\Cal{L}_g(\Gamma)^2 \leqslant \frac{2\pi}{\opLog(R)}.
$$
\noindent We obtain equality if and only if $\lambda=Kr^{-2}$ for some normalising factor $K\in (0,\infty)$, and an explicit calculation of $K$ permits us to obtain the
first result. The second result follows by a similar reasoning.\qed
\medskip
\noindent The following lemma permits us to compare the modules of two rings of which one is contained inside the other:
\proclaim{Lemma \nextprocno}
\noindent Let $A_1$ and $A_2$ be two rings. Let $i:A_1\rightarrow A_2$ be an embedding. If $i_*\pi_1(A_1) = \pi_1(A_2)$, then:
$$
\opMod(A_2) \leqslant \opMod(A_1).
$$
\endproclaim
\proclabel{LemmeComparerDeuxAnneaux}
\proof Let $\Gamma_2$ be the familly of curves in $A_2$ which are freely homotopic to a generator of $\pi_1(A_2)$. By the proof of the preceeding lemma, there exists a conformal
metric $g$ over $A_2$ such that:
$$\matrix
\opArea_g(A_2) \hfill&= 1, \hfill\cr
\Cal{L}_g(\Gamma_2) \hfill&= \opMod(A_2). \hfill\cr
\endmatrix$$
\noindent Let $\Gamma_1$ be the familly of curves in $A_1$ which are freely homotopic to a generator of $\pi_1(A_1)$. The mapping $i_*$ sends $\Gamma_1$ into $\Gamma_2$. We thus obtain:
$$
\Cal{L}_g(i_*\Gamma_1) \geqslant \Cal{L}_g(\Gamma_2).
$$
\noindent Let us define $h$ by:
$$
h = \frac{1}{\opArea_{i^*g}(A_1)}i^*g.
$$
\noindent Since $\opArea_g(A_1)\leqslant \opArea_g(A_2)=1$, it follows that $h\geqslant i^*g$ and consequently:
$$
\Cal{L}_h(\Gamma_1) \geqslant \Cal{L}_{i^*g}(\Gamma_1) \geqslant \Cal{L}_g(\Gamma_2) = \opMod(A_2).
$$
\noindent Since $\opArea_h(A_1)=1$, the result follows.\qed
\medskip
\noindent In particular, we obtain:
\proclaim{Corollary \nextprocno}
\noindent $\opMod(A)=0$ if and only if $A$ is conformally equivalent to $\Bbb{D}^*$ or to $\Bbb{C}^*$.
\endproclaim
\proclabel{CorClassificationDesAnneauxParModule}
\proof Let $R\in ]0,\infty[$ be a positive real number. We have:
$$
\left\{\frac{1}{R} < \left|z\right| < 1\right\} \subseteq \Bbb{D}^* \subseteq \Bbb{C}^*.
$$
\noindent Using the previous result, we thus obtain:
$$
\opMod(\Bbb{C}^*) \leqslant \opMod(\Bbb{D}^*) \leqslant \sqrt{\frac{2\pi}{\opLog(R)}}.
$$
\noindent By letting $R$ tends to infinity, we obtain:
$$
\opMod(\Bbb{C}^*) = \opMod(\Bbb{D}^*) = 0.
$$
\noindent The converse follows by the uniformisation principal.\qed
\newsubhead{Extending the Complex Structure}
\noindent Let $S$ be a surface and let $\Cal{P}$ be a discrete subset of $S$. Let $i:S\setminus\Cal{P}\rightarrow\Bbb{H}^3$ be an immersion such that the immersed surface
$\Sigma=(S\setminus\Cal{P},i)$ is a k-surface. Let $\hat{\mathi}$ be the Gauss lifting of $i$ and let us denote by $\overrightarrow{n}$ the Gauss-Minkowski mapping that sends
$U\Bbb{H}^3$ to $\partial_\infty\Bbb{H}^3=\hat{\Bbb{C}}$. Let $\Cal{H}$ be the canonical conformal structure over $\hat{\Bbb{C}}$. Let us define $\varphi$ by
$\varphi=\overrightarrow{n}\circ\hat{\mathi}$. Since $\varphi$ is a local homeomorphism $\varphi^*\Cal{H}$ defines a conformal structure over $S\setminus\Cal{P}$. In this
section, we will prove the lemma \procref{LemmaATSExtendConformalStructure}.
\medskip
\noindent Let $p_0$ be an arbitrary point in $\Cal{P}$. We suppose that $\hat{\Sigma}$ is asymptotically tubular of finite order about $p_0$. We obtain the following result:
\proclaim{Lemma \nextprocno}
\noindent For every sufficiently small neighbourhood $U$ of $p_0$ in $S$ which is homeomorphic to a disc, the Riemann surface $(U\setminus\left\{p_0\right\},\varphi^*\Cal{H})$ is conformally
equivalent to $\Bbb{D}^*$.
\endproclaim
\proclabel{LemmaATSPuncturedHolomorphicDisc}
\proof Let $\Cal{H}'$ be the conformal structure generated over $S\setminus\Cal{P}$ by the metric $\hat{\mathi}^*g^\nu$ and the canonical orientation of $S$. By
lemma \procref{LemmaCSQCEquivalence}, $\Cal{H}'$ is k-quasiconformally equivalent to $\varphi^*\Cal{H}$. It thus suffices to show that $(U\setminus\left\{p_0\right\},\Cal{H}')$
is conformally equivalent to $\Bbb{D}^*$.
\medskip
\noindent Let $\Gamma_{0,\infty}$ be the unique geodesic in $\Bbb{H}^3$ joining $0$ to $\infty$. We may assume that $\hat{\Sigma}$ is asymptotically tubular about
$\Gamma_{0,\infty}$. Let $N_{0,\infty}$ be the normal circle bundle over $\Gamma_{0,\infty}$ in $U\Bbb{H}^3$. Let $\opExp:TU\Bbb{H}^3\rightarrow U\Bbb{H}^3$ be the
exponential mapping over $U\Bbb{H}^3$.
\medskip
\noindent Let $T=(S^1\times(0,\infty),\hat{\mathj})$ be a half tube of order $k$ about $\Gamma_{0,\infty}$ such that there exists:
\medskip
\myitem{(i)} a neighbourhood $\Omega$ of $p_0$ in $S$,
\medskip
\myitem{(ii)} a graph diffeomorphism $\phi:S^1\times(0,\infty)\rightarrow\Omega\setminus\left\{p_0\right\}$, and
\medskip
\myitem{(iii)} a graph function $\lambda\in\Gamma(\hat{\mathj}^*NN_{0,\infty})$,
\medskip
\noindent such that:
\medskip
\myitem{(i)} $\varphi(e^{i\theta},t)\rightarrow p_0$ as $t$ tends to $+\infty$,
\medskip
\myitem{(ii)} for all $p\in\Bbb{N}$, $\|D^p\lambda(e^{i\theta},t)\|$ tends to $0$ as $t$ tends to $+\infty$, and
\medskip
\myitem{(iii)} $\hat{\mathi}\circ\phi = \opExp\circ\lambda$.
\medskip
\noindent For $R,T>0$, we define the set $A_{R,T}$ by:
$$
A_{R,T} = S^1\times(T,T+R).
$$
\noindent Since $\|D^1\lambda(e^{i\theta},t)\|$ tends to $0$ as $t$ tends to $+\infty$, we have:
$$
\left|\phi^*\hat{\mathi}^*g^\nu - \hat{\mathj}^*g^\nu\right|(e^{i\theta},t)\rightarrow 0\text{ as }t\rightarrow\infty.
$$
\noindent Consequently, if we denote by $d_{R,T}$ the complex dilatation of the metric $\phi^*\hat{\mathi}^*g^\nu$ relative to $\hat{\mathj}^*g^\nu$ over $A_{R,T}$, we find that
for all $R$:
$$
d_{R,T}\rightarrow 0\text{ as }T\rightarrow+\infty.
$$
\noindent Thus, by the translation invariance of $\hat{\mathj}^*g^\nu$, for all $R$ we obtain:
$$
\opMod(A_{R,T},\phi^*\hat{\mathi}^*g^\nu)\rightarrow\opMod(A_{0,R},\hat{\mathj}^*g^\nu)\text{ as }T\rightarrow+\infty.
$$
\noindent It follows by lemma \procref{LemmeComparerDeuxAnneaux} that, for all $R$:
$$
\opMod(S^1\times(0,\infty),\phi^*\hat{\mathi}^* g^\nu)\leqslant \opMod(A_{0,R},\hat{\mathj}^*g^\nu).
$$
\noindent Thus:
$$\matrix
\opMod(S^1\times(0,\infty),\phi^*\hat{\mathi}^*g^\nu)\hfill&\leqslant \opMod(S^1\times(0,\infty),\hat{\mathj}^*g^\nu)\hfill\cr
&=0\hfill\cr
\endmatrix$$
\noindent Consequently, by corollary \procref{CorClassificationDesAnneauxParModule}, $(\Omega\setminus\left\{p_0\right\},\hat{\mathi}^*g^\nu)$ is biholomorphic either
to $\Bbb{C}^*$ or to $\Bbb{D}^*$. Thus, by reducing $\Omega$ if necessary, we obtain the desired result.\qed
\medskip
\noindent We now obtain lemma \procref{LemmaATSExtendConformalStructure} as a corollary to this result:
\proclaim{Lemma \procref{LemmaATSExtendConformalStructure}}
\noindent Let $S$ be a surface. Let $\Cal{P}$ be a discrete subset of $S$. Let $i:S\setminus\Cal{P}\rightarrow\Bbb{H}^3$ be an immersion such that the immersed surface
$\Sigma=(S\setminus\Cal{P},i)$ is a $k$-surface. Let $\hat{\Sigma}=(S,\hat{\mathi})$ be the Gauss lifting of $\Sigma$. Let $\overrightarrow{n}$ be the Gauss-Minkowski
mapping that sends $U\Bbb{H}^3$ into $\partial_\infty\Bbb{H}^3=\hat{\Bbb{C}}$. Let us define $\varphi=\overrightarrow{n}\circ\hat{\mathi}$. Let $\Cal{H}$ be the canonical
conformal structure over $\hat{\Bbb{C}}$.
\medskip
\noindent Let $p_0$ be a point in $\Cal{P}$. If $\hat{\Sigma}$ is asymptotically tubular of finite order near $p_0$, then $\varphi^*\Cal{H}$ extends to a unique conformal
structure on $(S\setminus\Cal{P})\munion\left\{p_0\right\}$.
\endproclaim
\proof Let $U$ be a neighbourhood of $p_0$ in $S$ such that $U\setminus\left\{p_0\right\}$ is biholomorphic to $\Bbb{D}^*$ and let
$\alpha:U\setminus\left\{p_0\right\}\rightarrow\Bbb{D}^*$ be this biholomorphism. Let $\Omega$ be a neighbourhood of zero in $\Bbb{D}$ and let $\gamma$ be a simple closed
curve in $\Omega$ such that:
$$
0\in \opInt(\gamma)\setminus\left\{0\right\}
$$
\noindent Let us define the curve $\tilde{\gamma}$ by $\tilde{\gamma}=\alpha^{-1}\circ\gamma$. It is a simple closed curve in $U\setminus\left\{p_0\right\}$. It follows that the
complement of $\tilde{\gamma}$ in $U$ consists of two connected components $U_1$ and $U_2$. We may assume that $p_0\in U_1$. Since $\alpha$ is a homeomorphism, it sends
$U_1\setminus\left\{p_0\right\}$ either onto $\opInt(\gamma)\setminus\left\{0\right\}=\Omega$ or onto $\opExt(\gamma)\minter\Bbb{D}$. However, by the preceeding lemma:
$$
\opMod(U_1\setminus\left\{p_0\right\})=0.
$$
\noindent Consequently, the set $U_1\setminus\left\{p_0\right\}$ is not biholomorphic to $\opExt(\gamma)\minter\Bbb{D}$, and so:
$$
\alpha(U_1\setminus\left\{p_0\right\}) = \opInt(\gamma)\setminus\left\{0\right\}.
$$
\noindent It follows that $\alpha(p)$ tends to zero as $p$ tends to $p_0$ and we may thus extends $\alpha$ to a continuous mapping over $U$ by defining:
$$
\alpha(p_0) = 0
$$
\noindent Since $\alpha$ is bijective, by the princical of invariance of domains, it is a homeomorphism. We thus obtain a holomorphic chart $(\alpha, U, \Bbb{D})$ of
$(S\setminus\Cal{P})\munion\left\{p_0\right\}$ about $p_0$ which extends the conformal structure of $S\setminus\Cal{P}$, and we thus obtain existence. Uniqueness follows from the Cauchy removeable singularity theorem.\qed
\newsubhead{Extending the Holomorphic Function}
\noindent We continue to work with the construction of the previous section. We now obtain the following result:
\proclaim{Lemma \nextprocno}
\noindent If $\hat{\Sigma}$ is asymptotically tubular of finite order about $p_0$, then there exists a point $q_0\in\hat{\Bbb{C}}$ such that $\varphi(p)$ tends to
$q_0$ as $p$ tends to $p_0$.
\endproclaim
\proof Let $\Gamma_{0,\infty}$ be the geodesic joining $0$ to infinity. We may assume that $\Sigma$ is asymptotically tubular about $\Gamma_{0,\infty}$. Let $N_{0,\infty}$
be the normal circle bundle of $\Gamma_{0,\infty}$ in $U\Bbb{H}^3$. Let $\pi:N_{0,\infty}\rightarrow\Gamma_{0,\infty}$ be the canonical projection. Let
$\gamma:\Bbb{R}\rightarrow\Gamma_{0,\infty}$ be a unit speed parametrisation of $\Gamma_{0,\infty}$ such that $\gamma(t)$ tends to $0$ as $t$ tends to $+\infty$. By identifying
$\Bbb{H}^3$ with $\Bbb{C}\times(0,\infty)$, we define $h:\Bbb{R}\rightarrow (0,\infty)$ such that, for all $t$:
$$
\gamma(t) = (0,h(t)).
$$
\noindent For all $t\in\Bbb{R}$, we define $A_t\in\opIsom(\Bbb{H}^3)$ by:
$$
A_t(z,s) = \left(\frac{z}{h(t)},\frac{s}{h(t)}\right).
$$
\noindent We also denote by $A_t$ the actions of $A_t$ on $U\Bbb{H}^3$ and $\partial_\infty\Bbb{H}^3$.
\medskip
\noindent Let $T=(S^1\times(0,\infty),\hat{\mathj})$, $\Omega$, $\phi$ and $\lambda$ be as in the proof of lemma \procref{LemmaATSPuncturedHolomorphicDisc}. For all
$t\in\Bbb{R}$, we have:
$$\matrix
&\pi\circ A_t\circ\hat{\mathj}(e^{i\theta},t) \hfill&= (0,1) \hfill\cr
\Rightarrow\hfill &\left|\overrightarrow{n}\circ A_t\circ\hat{\mathj}(e^{i\theta},t)\right| \hfill&= 1\hfill\cr
\endmatrix$$
\noindent Since $\lambda(e^{i\theta},t)$ tends to $0$ as $t$ tends to $+\infty$, we have:
$$\matrix
&\left|\overrightarrow{n}\circ A_t\circ\hat{\mathj}(e^{i\theta},t)-\overrightarrow{n}\circ A_t\circ\hat{\mathi}\circ\phi(e^{i\theta},t)\right|\rightarrow 0
\text{ as } t\rightarrow +\infty,\hfill\cr
\Rightarrow\hfill&\left|\overrightarrow{n}\circ A_t\circ\hat{\mathi}\circ\phi(e^{i\theta},t)\right|\rightarrow 1\text{ as } t\rightarrow +\infty,\hfill\cr
\Rightarrow\hfill&\left|A_t\circ\overrightarrow{n}\circ\hat{\mathi}\circ\phi(e^{i\theta},t)\right|\rightarrow 1\text{ as } t\rightarrow +\infty,\hfill\cr
\Rightarrow\hfill&\left|\overrightarrow{n}\circ\hat{\mathi}\circ\phi(e^{i\theta},t)\right|\rightarrow 0\text{ as } t\rightarrow +\infty,\hfill\cr
\Rightarrow\hfill&\varphi(p)\rightarrow 0\text{ as }p\rightarrow p_0.\hfill\cr
\endmatrix$$
\noindent The result now follows.\qed
\medskip
\noindent We now obtain lemma \procref{LemmaATSExtendFunction} as a corollary to this result:
\proclaim{Lemma \procref{LemmaATSExtendFunction}}
\noindent With the same hypotheses as in lemma \procref{LemmaATSExtendConformalStructure}, let $p_0$ be a point in $\Cal{P}$. If $\hat{\Sigma}$ is asymptotically
tubular of order $k$ near $p_0$, then $\varphi$ extends to a unique holomorphic function over $(S\setminus\Cal{P})\munion\left\{p_0\right\}$ having a critical point of order
$k$ at $p_0$.
\endproclaim
\proof It follows by Cauchy's removeable singularity theorem and the previous lemma that $\varphi$ extends to a unique holomorphic function over
$(S\setminus\Cal{P})\munion\left\{p_0\right\}$.
\medskip
\noindent Using the same reasoning and notation as in the preceeding lemma, we find that there exists $T>0$ such that for $t\geqslant T$,
and for all $\theta$:
$$\matrix
&\left|(A_t\circ\varphi\circ\phi)(e^{i\theta},t) - (A_t\circ\overrightarrow{n}\circ\hat{\mathj})(e^{i\theta},t)\right|\hfill&<1\hfill\cr
\Rightarrow\hfill&\left| (A_t\circ\varphi\circ\phi)(e^{i\theta},t) - e^{ik\theta}\right|\hfill&<1\hfill\cr
\endmatrix$$
\noindent It follows that, for $t\geqslant T$, the curve $\theta\mapsto(\varphi\circ\phi)(e^i\theta,t)$ is homotopic in $\Bbb{C}^*$
to the curve $\theta\mapsto e^{ik\theta}$, and it thus turns $k$ times round the origin.
\medskip
\noindent By lemma \procref{LemmaATSExtendConformalStructure}, we may assume that $(\phi(S^1\times(0,[T,\infty)),\varphi^*\Cal{H})$ is biholomorphic to $\Bbb{D}^*$.
Let $\alpha:\phi(S^1\times[T,\infty))\rightarrow\Bbb{D}^*$ be a biholomorphic mapping. For $t\geqslant T$, since $\alpha$ is a homeomorphism,
the curve $\theta\mapsto (\alpha\circ\phi)(e^{i\theta},t)$ turns once around the origin. The result now follows by considering $(\varphi\circ\alpha^{-1})$.\qed
\newhead{Bibliography}
{\leftskip = 5ex \parindent = -5ex
\leavevmode\hbox to 4ex{\hfil\cite{BallGromSch}}\hskip 1ex{Ballman W., Gromov M., Schroeder V., {\sl Manifolds of nonpositive curvature}, Progress in Mathematics, {\bf 61}, Birkh\"auser, Boston, (1985)}
\medskip
\leavevmode\hbox to 4ex{\hfil\cite{Grom}}\hskip 1ex{Gromov M., Foliated plateau problem, part I : Minimal varieties, {\sl GAFA} {\bf 1}, no. 1, (1991), 14--79}%
\medskip
\leavevmode\hbox to 4ex{\hfil\cite{LabB}}\hskip 1ex{Labourie F., Probl\`emes de Monge-Amp\`ere, courbes holomorphes et laminations, {\sl GAFA} {\bf 7}, no. 3, (1997), 496--534}
\medskip
\leavevmode\hbox to 4ex{\hfil\cite{LabA}}\hskip 1ex{Labourie F., Un lemme de Morse pour les surfaces convexes, {\sl Invent. Math.} {\bf 141} (2000), 239--297}
\medskip
\leavevmode\hbox to 4ex{\hfil\cite{LV}}\hskip 1ex{Lehto O., Virtanen K. I., {\sl Quasiconformal mappings in the plane}, Die Grund\-lehren der mathemathischen Wissenschaften, {\bf 126}, Springer-Verlag, New York-Heidelberg, (1973)}
\medskip
\leavevmode\hbox to 4ex{\hfil\cite{Muller}}\hskip 1ex{Muller M.P., Gromov's Schwarz lemma as an estimate of the gradient for holomorphic curves, In: {\sl Holomorphic curves in symplectic geometry},
Progress in Mathematics, {\bf 117}, Birkh\"auser, Basel, (1994)}
\medskip
\leavevmode\hbox to 4ex{\hfil\cite{RosSpruck}}\hskip 1ex{Rosenberg H., Spruck J. On the existence of convex hyperspheres of constant Gauss curvature in hyperbolic space, {\sl J. Diff. Geom.} {\bf 40} (1994), no. 2,
379--409} 
\medskip
\leavevmode\hbox to 4ex{\hfil\cite{SmiA}}\hskip 1ex{Smith G., Special Legendrian structures and Weingarten problems, Preprint, Orsay (2005)}%
\medskip
\leavevmode\hbox to 4ex{\hfil\cite{SmiB}}\hskip 1ex{Smith G., Hyperbolic Plateau problems, Preprint, Orsay (2005)}%
\medskip
\leavevmode\hbox to 4ex{\hfil\cite{SmiE}}\hskip 1ex{Smith G., Th\`ese de doctorat, Paris (2004)}%
\par
}%
\enddocument

%% file: references.tex
\global\def\_@citation@BallGromSch{1}
\global\def\_@citation@Grom{2}
\global\def\_@citation@LabB{3}
\global\def\_@citation@LabA{4}
\global\def\_@citation@LV{5}
\global\def\_@citation@Muller{6}
\global\def\_@citation@RosSpruck{7}
\global\def\_@citation@SmiA{8}
\global\def\_@citation@SmiB{9}
\global\def\_@citation@SmiE{10}
\global\def\_@proc@PresentationChIIILimites{1.1}
\global\def\_@proc@PresentationChIIIRevetementsRamifiees{1.2}
\global\def\_@proc@PresentationChIIISurfacesAsymptotiquementTubulaires{1.3}
\global\def\_@proc@LemmaGeometricMaximumPrincipal{2.1}
\global\def\_@proc@LemmeConvergenceDesEnsemblesEtDesFonctions{2.3}
\global\def\_@proc@LemmeLIntersectionDesLimites{2.4}
\global\def\_@proc@LemmeLimitesDInverses{2.5}
\global\def\_@proc@LemmeInjectiviteDeLaLimite{2.6}
\global\def\_@proc@LemmeInjectiviteLorsquOnSApprocheDuLimite{2.7}
\global\def\_@proc@LemmeCompacteDansLImageDeLaLimite{2.8}
\global\def\_@head@HeadDifferentialGeometryInTTM{3}
\global\def\_@proc@LemmaCSQCEquivalence{3.2}
\global\def\_@proc@ThmSmithExistence{4.1}
\global\def\_@proc@LemmeLabNature{4.2}
\global\def\_@proc@ThmLabCompacite{4.3}
\global\def\_@proc@LemmeSolutionEstgraphe{6.1}
\global\def\_@proc@LemmeMetriqueEstMince{6.2}
\global\def\_@proc@LemmeConvergenceDeT{6.3}
\global\def\_@subhead@SubHeadGraph{6.2}
\global\def\_@proc@Limites{6.4}
\global\def\_@proc@LemmeBehaviourOfIAtLimites{6.5}
\global\def\_@fig@FigureConvexiteStricteDeSigma{6.1}
\global\def\_@subhead@SubHeadBehaviourOfMetric{6.3}
\global\def\_@proc@CorLaFormeDeLagraphe{6.10}
\global\def\_@subhead@SubHeadHaussdorfConvergence{6.4}
\global\def\_@proc@LemmeTrivDesRevetementsRamifiees{7.1}
\global\def\_@proc@RevetementsRamifies{7.2}
\global\def\_@proc@LemmeSurfaceTubulaireAutourDeGamma{7.3}
\global\def\_@proc@LemmeConvergeVersUnTubeDOrdreK{7.4}
\global\def\_@proc@CorLocallyGraphsOverTubes{7.8}
\global\def\_@proc@LemmaATSExtendConformalStructure{8.1}
\global\def\_@proc@LemmaATSExtendFunction{8.2}
\global\def\_@proc@LemmeCalculerLeModuleDUnAnneau{8.3}
\global\def\_@proc@LemmeComparerDeuxAnneaux{8.4}
\global\def\_@proc@CorClassificationDesAnneauxParModule{8.5}
\global\def\_@proc@LemmaATSPuncturedHolomorphicDisc{8.6}